\pdfoutput=1
\documentclass[a4paper,11pt]{article}
\usepackage[english]{babel}

\usepackage{hyperref}
\hypersetup{
colorlinks=true,
linkcolor=blue,
filecolor=magenta,      
urlcolor=cyan,
citecolor=[RGB]{0,204,0},
}
\usepackage{resizegather}
\usepackage[activate={true, nocompatibility}, final, tracking=true, kerning=true, spacing=true, factor=1100, stretch=10, shrink=10]{microtype}

\usepackage[left=2.9cm,top=3cm,right=2.9cm,bottom=3cm]{geometry}

\usepackage{graphicx,wrapfig}

\usepackage{amsfonts,amsthm,amsmath,amssymb,enumerate,enumitem,mathtools,tikz,stmaryrd}
\setlist[enumerate,1]{label={(\arabic*)}}

\expandafter\def\csname opt@stmaryrd.sty\endcsname
{only,shortleftarrow,shortrightarrow}
\usepackage{tikz-cd}
\usepackage[all]{xy}
\usepackage{adjustbox}
\usepackage{extpfeil}
\usepackage{comment}
\usepackage{mathabx}
\usepackage{float}
\usepackage{caption}
\usepackage[toc]{appendix}
\usepackage{dsfont}
\usetikzlibrary{matrix}
\usepackage[nameinlink, capitalize]{cleveref}
\usepackage[normalem]{ulem}
\usepackage{sidecap}
\usepackage{mathrsfs,mathscinet}

\usepackage{quiver}

\usepackage{color} %CLASHES with hyperref
\definecolor{darkgreen}{rgb}{0.0, 0.7, 0.0}
\newenvironment{??}{\noindent \color{darkgreen}{\bf ???:} \footnotesize}{}
\definecolor{cyan}{cmyk}{1,0,0,0}

\newcommand{\cc}{\color{cyan}}

\newcommand{\bdg}{\begin{dg}}

\theoremstyle{definition}
\newtheorem{thm}{Theorem}[section]
\newtheorem{coroll}[thm]{Corollary}
\newtheorem{construction}[thm]{Construction}
\newtheorem{defn}[thm]{Definition}
\newtheorem{lemma}[thm]{Lemma}

\newtheorem{prop}[thm]{Proposition}
\newtheorem{question}[thm]{Question}
\newtheorem{remark}[thm]{Remark}

\newtheorem*{thm*}{Main Theorem}
\newtheorem{notn}[thm]{Notation}

\newtheorem{example}[thm]{Example}
\newtheorem{setup}[thm]{Setup}
\tikzset{
  symbol/.style={
    draw=none,
    every to/.append style={
      edge node={node [sloped, allow upside down, auto=false]{$#1$}}}
  }
}
\microtypecontext{spacing=nonfrench}

\DeclareMathOperator{\RHom}{RHom}
\DeclareMathOperator{\Filt}{Filt}

\newcommand{\codim}{\mathrm{codim}}

\newcommand{\cO}{\mathcal{O}}

\newcommand{\cC}{\mathcal{C}}

\newcommand{\cI}{\mathcal{I}}

\newcommand{\cH}{\mathcal{H}}

\newcommand{\cN}{\mathcal{N}}

\newcommand{\cU}{\mathcal{U}}
\newcommand{\cW}{\mathcal{W}}
\newcommand{\cX}{\mathcal{X}}

\newcommand{\cL}{\mathcal{L}}

\newcommand{\bG}{\mathbb{G}}

\renewcommand{\P}{\mathbb{P}}

\newcommand{\bZ}{\mathbb{Z}}

\newcommand{\Extstack}{\mathscr{E}xt}
\newcommand{\sU}{\mathscr{U}}
\newcommand{\sX}{\mathscr{X}}
\newcommand{\sY}{\mathscr{Y}}
\newcommand{\sZ}{\mathscr{Z}}

\newcommand{\hot}{{\mathrm{h.o.t.}}}
\newcommand{\lot}{{\mathrm{l.o.t.}}}

\newcommand{\rnc}{{\mathrm{nc}\text{-}}}
\newcommand{\rtf}{{\mathrm{tf}}}
\newcommand{\ST}{\overline{\mathrm{ST}}}
\newcommand{\univ}{{\mathrm{univ}}}

\newcommand\blank{\underline{\ \ }}

\newcommand{\on}{\operatorname}

\newcommand{\ev}{\mathrm{ev}}

\DeclareMathOperator{\ch}{ch}

\DeclareMathOperator{\Coker}{Coker}

\newcommand{\Ext}{ \on{Ext}}
\DeclareMathOperator{\Fil}{Fil}
\DeclareMathOperator{\Fit}{Fit}

\DeclareMathOperator{\gr}{gr}
\newcommand{\Hilb}{\on{Hilb}}

\newcommand{\Hom}{ \on{Hom}}
\DeclareMathOperator{\id}{id}

\newcommand{\Map}{ \on{Hom}}

\newcommand{\Quot}{\on{Quot}}

\DeclareMathOperator{\rk}{rk}
\newcommand{\Spec}{\on{Spec}}

\DeclareMathOperator{\Supp}{Supp}
\newcommand{\Td}{\on{Td}}
\newcommand{\TF}{\on{TF}}
\newcommand{\Tors}{\on{Tors}}
\DeclareMathOperator{\triv}{triv}

\newcommand{\Pic}{\on{Pic}}

\newcommand{\SCoh}{\mathscr{C}oh}
\newcommand{\SPic}{\mathscr{P}ic}
\newcommand{\Tot}{\mathrm{Tot}}

\usepackage{color} %CLASHES with hyperref
\definecolor{darkgreen}{rgb}{0.0, 0.7, 0.0}

\newcommand{\nps}{\mathcal{U}}

\begin{document}
\title{\textbf{Moduli of non-pure sheaves and\\
contractions of Hilbert schemes}}
\author{Andres Fernandez Herrero and Svetlana Makarova}
\date{}

\maketitle

\begin{abstract}
We define certain stacks of rank one sheaves on a smooth projective variety, and show that they admit proper good moduli spaces.
We offer several applications to contractions of subschemes inside Hilbert schemes of points.
We construct a surgery diagram via non-GIT wall-crossing, and use the interpretation of the surgery as a fine moduli of sheaves to prove instances of Kawamata's DK-hypothesis in this setting.
\end{abstract} 

\tableofcontents

\section{Introduction}
This paper is concerned with the existence of contractions of subschemes inside a given variety. To be more precise, we are motivated by the following:

\begin{question}
    Let $X$ be a smooth projective variety, and let $j:Y \hookrightarrow X$ be a closed subscheme equipped with a surjective morphism $\pi: Y \to B$.
    We are interested in conditions that guarantee the existence of a proper algebraic space $\overline{X}$ equipped with a surjective morphism $f: X \to \overline{X}$ that contracts $Y$ to $B$.
    I.e., we would like to produce a closed immersion $B \hookrightarrow \overline{X}$ such that the restriction $f: f^{-1}(\overline{X}\setminus B) \to \overline{X} \setminus B$ is an isomorphism and such that the following diagram commutes
 \[
\begin{tikzcd}
 Y \ar[r, "j"]
  \ar[d, "\pi"] & X \ar[d, "f"] \\ B \ar[r, symbol = \hookrightarrow] &
  \overline{X}
\end{tikzcd}
\]
\end{question}

The contraction theorem in birational geometry \cite[Thm. 3.7 (3)]{kollar-mori} provides a way to construct contractions of $X$ associated to $\omega_X$-negative extremal faces in the cone of curves. However, we are potentially interested in situations where the contraction cannot be constructed with a semiample line bundle on $X$. A more general study of contractions was conducted in Artin's paper \cite{artin_algebraization_formal_ii}, which provides certain conditions \cite[Thm. 6.2]{artin_algebraization_formal_ii} guaranteeing the existence of $\overline{X}$ as an algebraic space.
In practice, verifying the hypotheses of \textit{loc.~cit.} usually requires a good understanding of the formal completion of $Y \subset X$. 
Artin's methods seem particularly powerful under strict negativity assumptions on the normal sheaf $\cN_{Y/X}:=j^*(\cI_Y)^{\vee}$ (for example, when the restriction of $\cN_{Y/X}$ to each $\pi$-fiber is an anti-ample vector bundle). Some examples of results under such strict negativity assumptions are provided by \cite[Cor. 6.12]{artin_algebraization_formal_ii} when $Y \subset X$ is a divisor (see also \cite{fujiki-nakano} and \cite{grauert_contractions} for analogs in the analytic context). 

We aim to study the more delicate situation where we allow the restriction of $\cN_{Y/X}$ to the $\pi$-fibers to be non-positive (for example, if we require the restrictions to be anti-nef vector bundles).
In this case it is necessary to specify in addition that $\pi: Y \to B$, viewed as a family of subschemes of $X$, satisfies a (set-theoretic) versality condition: indeed, if this family of subschemes deformed further, then we would be forced to contract those deformations as well.
More precisely, we need to require that the closed image of the induced morphism $B \to \Hilb(X)$ to the Hilbert scheme is also open; we denote this image by $\alpha$. Roughly speaking, $\alpha$ is specifying the data of the families of subvarieties that we want to contract.

We approach this question by viewing $X$ as the good moduli space of the stack of ideal sheaves of a single point on $X$, and adding certain rank one sheaves containing torsion.
We call $\nps_{\alpha,1}$ the resulting stack of sheaves.
More generally, we define an enlargement $\nps_{\alpha, \ell}$ of the stack of ideal sheaves of $\ell$ points, which realizes contractions of the Hilbert scheme $\Hilb^{\ell}(X)$ of $\ell$ points on $X$.

\paragraph{\textbf{Moduli of non-pure sheaves.}}
For the setup of our moduli problem, we fix a finite union of connected components $\alpha \subset \Hilb(X)$ such that every point $[W \subset X] \in \alpha$ satisfies that $W$ is geometrically integral, and either $W$ is geometrically normal or $W$ is a divisor. We also fix a positive integer $\ell>0$. Our results hold subject to a technical condition on the pair $(\alpha, \ell)$ (see condition (LR) in \Cref{defn: condition dagger}), which in practice amounts to a bound on the positivity of normal bundles of elements $[W \subset X] \in \alpha$.

\begin{figure}[h]
    \centering
    \includegraphics[width=0.3\linewidth]{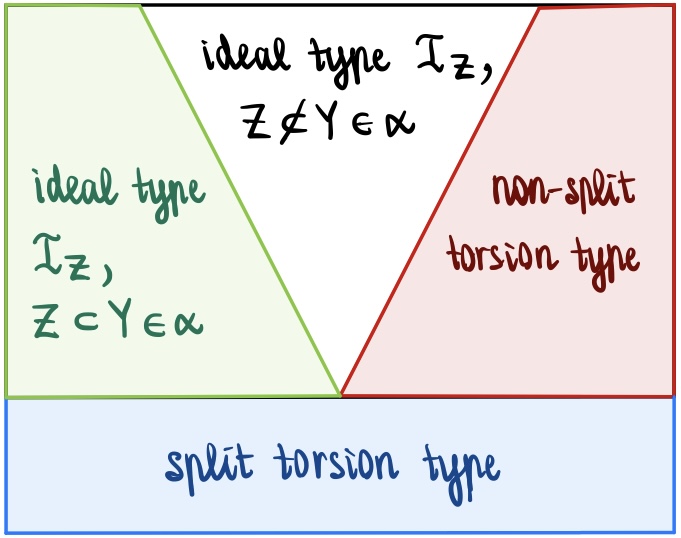}
    \captionof{figure}{Stack $\nps_{\alpha,\ell}$.}
    \label{fig: stack U_alpha l}
\end{figure}

Our stack $\nps_{\alpha, \ell}$ parameterizes coherent sheaves $F$ on $X$ of one of the following forms:

\begin{itemize}
    \item $F$ is the ideal sheaf $\cI_Z$ of some length $\ell$ finite subscheme $Z \subset X$, or
    \item $F$ fits into a short exact sequence $0 \to i_*(G) \to F \to \cI_W \to 0$,
    where $[i:W' \hookrightarrow X]$ and $[W \subset X]$ belong to the same connected component $\alpha_j \subset \alpha$, and $G$ is a rank one torsion-free sheaf on $W'$ with Hilbert polynomial $P_{\cO_{W'}} - \ell$ (here $P_{\cO_{W'}}$ denotes the Hilbert polynomial of the coherent sheaf $\cO_{W'}$).
    The locus of such sheaves is depicted in blue and red in \cref{fig: stack U_alpha l}, depending on whether the short exact sequence splits or not, respectively.
\end{itemize}

See \Cref{defn: stack U_alpha} for a more formal definition of $\nps_{\alpha, \ell}$ as a substack of a closed substack $\SCoh(X)_{\cO_X, P_{\ell}}$ of the stack of coherent sheaves on $X$.
The first main result of this paper is the following:
\begin{thm}[see \cref{thm: main theorem on U_alpha} for a more general statement]
\label{thm: main theorem intro 1}
    Assume that the pair $(\alpha, \ell)$ satisfies condition (LR) as in \cref{defn: condition dagger} and that $\alpha$ is of finite type.
    Then the substack $\nps_{\alpha,\ell} \subset \SCoh(X)_{\cO_X, P_{\ell}}$ is open, admits a proper good moduli space denoted $\Hilb^{\ell}(X)_{\alpha}$, and there is a morphism $f:\Hilb^\ell(X) \to \Hilb^\ell(X)_\alpha$.
\end{thm}

\begin{remark}
    It can be checked that two closed points $[Z \subset X], [Z' \subset X] \in \Hilb^{\ell}(X)$ have the same image in $\Hilb^{\ell}(X)_{\alpha}$ if and only if there exists some $[W \subset X] \in \alpha$ such that $Z, Z' \subset W$ and the ideal sheaves of $Z$ and $Z'$ inside $W$ are isomorphic. Hence, we have good control on the locus that gets contracted by $f$.
\end{remark}

The proper algebraic space $\Hilb^{\ell}(X)_{\alpha}$ is not a scheme in general (see \Cref{example: mukai flop contractions}).
In particular, the techniques of Geometric Invariant Theory cannot be used to prove the existence of this moduli of sheaves.
Instead, we employ the intrinsic approach developed by Alper, Halpern-Leistner and Heinloth \cite{AHLH} to show our main theorem in Sections \ref{section: openness and boundedness of the stack} and \ref{section: existence of moduli space}.

In \cref{fig: stack U_alpha l}, if one removes the red and blue loci of sheaves with torsion, one recovers the stack whose good moduli space is the Hilbert scheme $\Hilb^{\ell}(X)$.
One may then wonder if the status of the green locus of ideal sheaves and the red locus of non-split sheaves with torsion is symmetric;
namely, what would happen if we considered the substack $\nps_{\alpha,\ell}^{\circ} \subset \nps_{\alpha,\ell}$ obtained by removing the blue and green loci from $\nps_{\alpha,\ell}$.
To that end, we prove the following result.

\begin{thm}[see \cref{thm: existence of good moduli space other side of the wall}]
\label{thm: main theorem intro 2}
    Assume that $(\alpha, \ell)$ satisfies condition (LR) as in \cref{defn: condition dagger} and that $\alpha$ is of finite type.
    Then the substack $\nps_{\alpha,\ell}^\circ \subset \nps_{\alpha,\ell}$ is open, admits a proper good moduli space denoted $\Hilb^{\ell}(X)_{\alpha}^\circ$, and there is an induced morphism $\Hilb^\ell(X)_\alpha^\circ \to \Hilb^\ell(X)_\alpha$.
\end{thm}

\begin{figure}[h]
    \centering
    \includegraphics[width=0.9\linewidth]{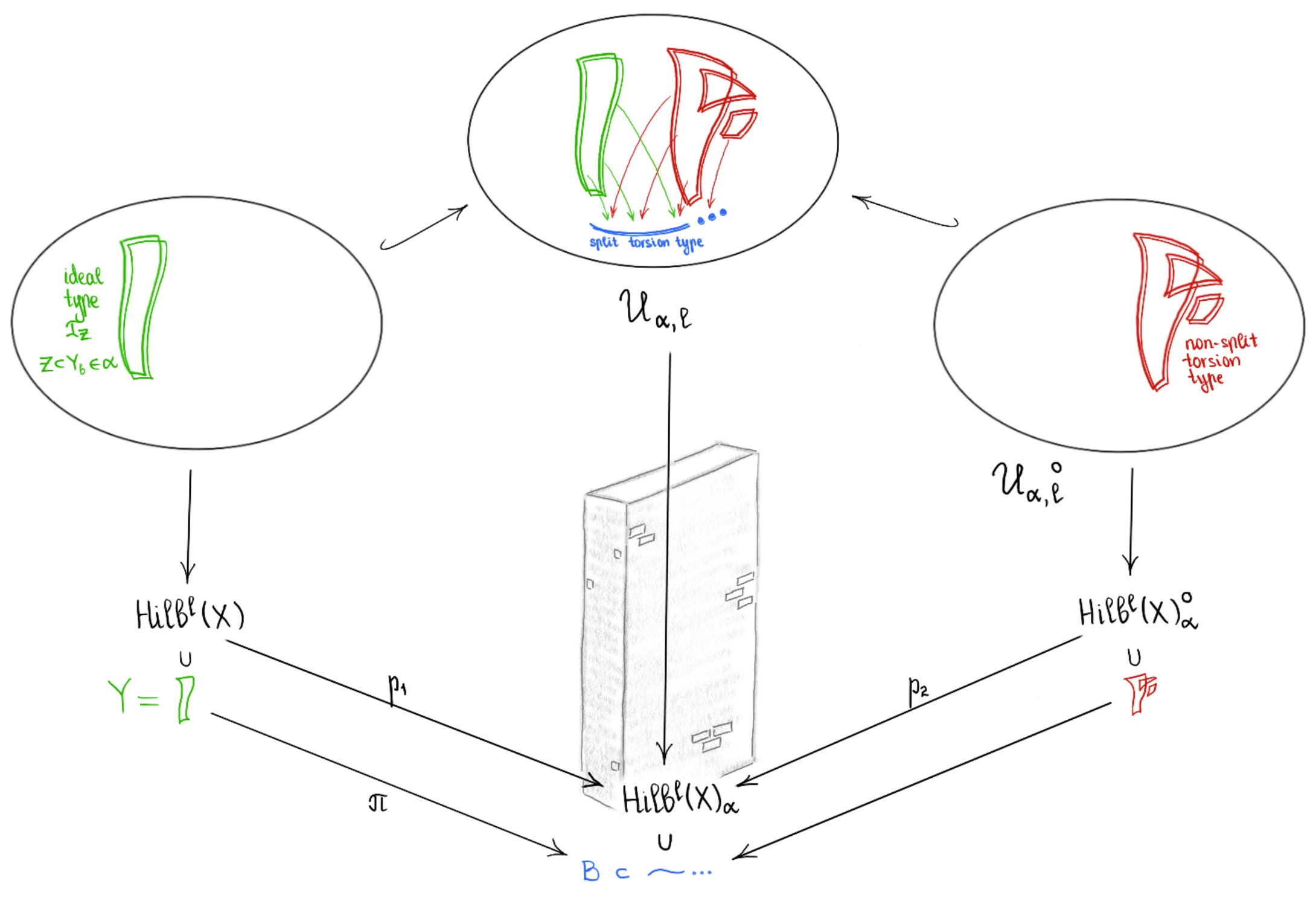}
    \caption{Surgery diagram (wall-crossing).}
    \label{fig: surgery diagram}
\end{figure}

We use \cref{fig: surgery diagram} to visualize the above theorem.
In the picture, the locus of sheaves of ideal type (green on the left, see \cref{defn: stack U_alpha}) and the locus of sheaves of non-split torsion type (red on the right, see \cref{remark: split and nonsplit torsion type}) behave like loci of some semistable objects that get interchanged under wall-crossing.
The blue locus of sheaves of split torsion type consists of some closed points of \(\cU_{\alpha,\ell}\) (by \cref{prop: closed points}) and lies in the closure of both green and red loci.
The arrows indicate the existence of specializations, which in our case can be realized as degenerations of non-split short exact sequences to the split ones.
Thus, the blue locus behaves like the locus of polystable sheaves at the wall that contracts semistable objects on both the green and the red sides.
We point out that on the level of good moduli spaces, the green locus does not necessarily surject onto the blue locus.
In the rest of the paper, we will refer to this kind of picture as a surgery diagram.

\paragraph{\textbf{Applications to contractions.}}
After possibly passing to subspaces of $\Hilb^\ell(X)_\alpha$ and $\Hilb^\ell(X)_\alpha^\circ$, and possibly also their normalizations, the construction recovers classically known contractions, flips and flops in certain cases, and provides new examples.

\begin{example}[Contracting rational curves]
    Consider the toy setting when $\ell=1$, $B$ is a point and $Y \subset X$ is a smooth curve.
    In this case, \cite{jimenez_contraction} provides criteria that characterize when $Y$ can be contracted, which rely on a good understanding of the formal model around $Y$.
    A conceptually different set of hypotheses that guarantee the existence of a contraction of a rational curve $\mathbb{P}^1 \cong Y \subset X$ inside a smooth threefold $X$ appears in the work of Reid \cite[Part II]{mmp_threefolds_reid} (see also related works \cite{laufer_exceptional} and \cite{pinkham_birational_factorization}).
    The result in \cite{mmp_threefolds_reid} states that $Y$ can be contracted to a point if the following are satisfied:
    \begin{itemize}
        \item (Non-positivity). The normal bundle $\cN_{Y/X}$ is isomorphic to either $\cO_{\mathbb{P}^1}(-1)^{\oplus 2}$ or $\cO_{\mathbb{P}^1} \oplus \cO_{\mathbb{P}^1}(-2)$.
        \item (Set-theoretic versality). The corresponding point $[Y \subset X] \in \Hilb(X)$ is open and closed in the Hilbert scheme.
    \end{itemize}
    The proof in \cite{mmp_threefolds_reid} constructs the contraction $\overline{X}$ explicitly by classifying the possible local models around $Y$.

    A corollary of our main theorems for $\ell=1$ fully recovers Reid's contraction result, and extends it to higher dimensions provided that $\P^1_k \cong Y \hookrightarrow X$ satisfies non-positivity for its normal bundle and set-theoretic versality (see \cref{cor: contraction of P^1 that don't move}).
    In this case \cref{fig: surgery diagram} specializes to \cref{fig: surgery diagram for contracting P1}, where, if we denote by $p \in \overline X$ the image of the contracted $\P^1$,
    the preimage $\mathbb{P}^1_k \cong Y = p_1^{-1}(p)$ in $X$ is replaced after the surgery with the projectivization of the vector space $H^1(Y, \cN_{Y/X}(-1))$ (see \Cref{example: surgery of non-movable rational curves}). 
    \begin{figure}[h]
        \centering
        \includegraphics[width=0.9\linewidth]{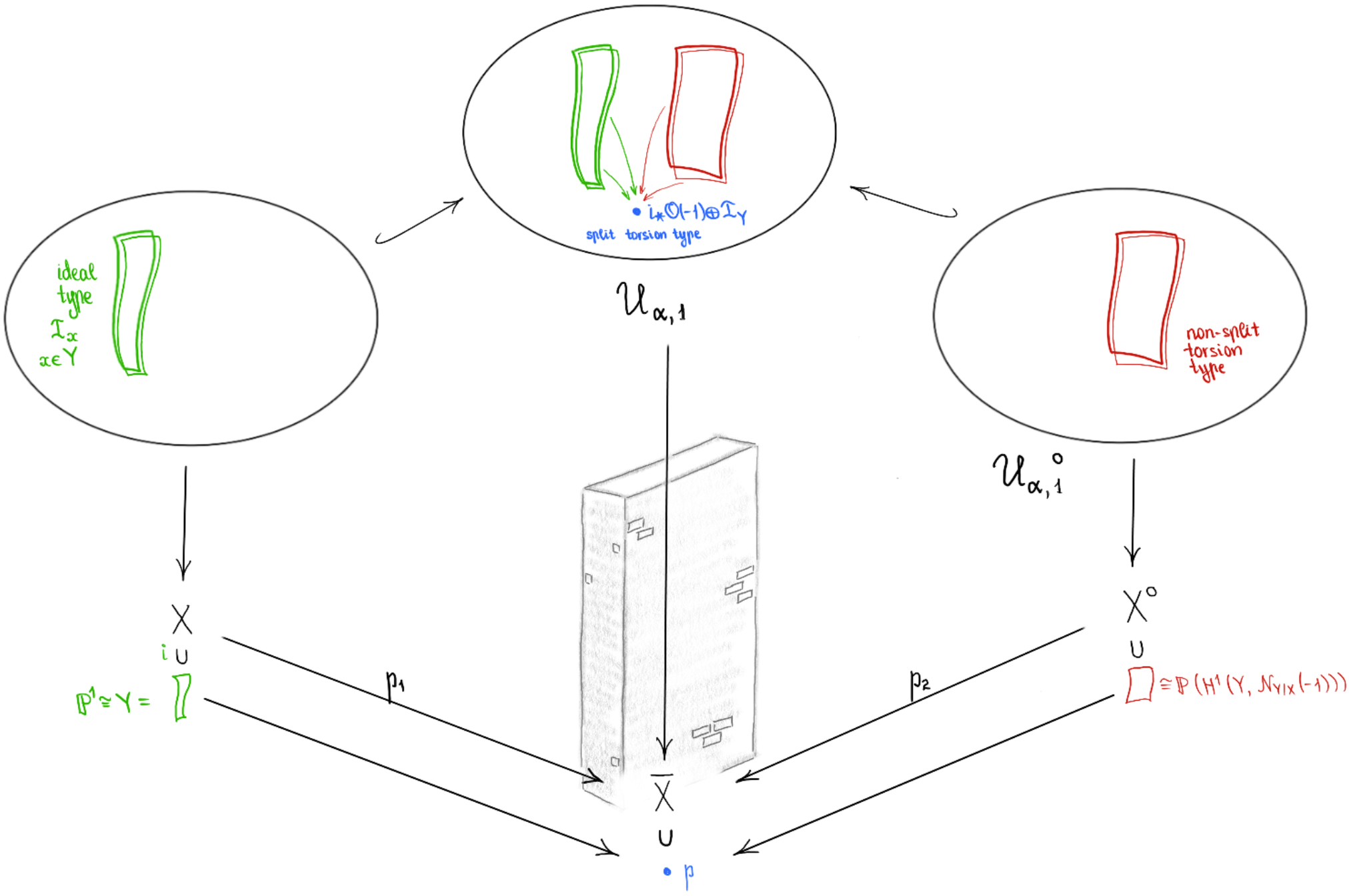}
        \caption{Surgery diagram for contracting $\P^1_k$.}
        \label{fig: surgery diagram for contracting P1}
    \end{figure}

    One may also use our results in this paper to study contractions of non-movable rational curves whose conormal bundle is not necessarily nef, see \cref{example: beyond nonpositive normal bundle};
    or contractions of families of rational curves $Y \to B$ where $B$ may be positive-dimensional, see \cref{subsection: contraction of projective fibrations}. 
\end{example}

\begin{example}[The DK hypothesis]
    The construction in the previous example realizes $X^{\circ}$ as a fine moduli of sheaves in $X$, which is related (a posteriori after the construction of the contraction) to the moduli of perverse ideal sheaves defined by Bridgeland in \cite[\S 5]{bridgeland_flops} whenever $\overline{X}$ has rational singularities.
    If $X^{\circ}$ is smooth and $\deg(i^*(\omega_X)) \leq 0$, we show that there is a fully faithful embedding of derived categories of coherent sheaves $D^b(X^{\circ}) \hookrightarrow D^b(X)$, which is an equivalence if $\deg(i^*(\omega_X))=0$ (see \cref{thm: derived categories of surgeries} for a more general statement).
\end{example}

\begin{example}[Abel-Jacobi contractions, see \cref{subsection: Abel-Jacobi contractions}]
    Let $X$ be a smooth projective surface, and let $\alpha \subset \Hilb(X)$ be a connected component of the Hilbert scheme which parameterizes integral subcurves of $X$. The intersection $C\cdot C'$ of any two elements $[C\subset X], [C' \subset X] \in \alpha$ remains constant as we vary $C$ and $C'$; we denote this number by $I(\alpha)$. We denote the universal family by $\cC \to \alpha$.
    
    Given a positive integer $\ell>I(\alpha)$, we can form the relative Hilbert scheme $\Hilb^{\ell}(\cC/\alpha)$, and the relative compactified Jacobian $\overline{\Pic}(\cC/\alpha)$ of rank one torsion-free sheaves on the family of integral curves $\cC \to \alpha$. There is an Abel-Jacobi morphism $AJ_{\ell}: \Hilb^{\ell}(\cC/\alpha) \to \overline{\Pic}(\cC/\alpha)$ which sends a length $\ell$ subscheme $Z$ to its ideal sheaf. We denote by $B$ the image of $AJ_{\ell}$.
    
    There is a natural closed injective morphism $Y:= \Hilb^{\ell}(\cC/\alpha) \to \Hilb^{\ell}(X)$ whose image consists of subschemes $[Z \subset X] \in \Hilb^{\ell}(X)$ that are contained in some curve $[C \subset X] \in \alpha$. \Cref{prop: abel jacobi contractions} asserts that the fibers of $AJ_{\ell}: Y \to B$ can be contracted, in the sense that there is a commutative diagram of proper algebraic spaces
    \[
        \begin{tikzcd}
         Y \ar[d, "AJ_{\ell}"]
          \ar[r] & \Hilb^{\ell}(X) \ar[d, "f"]\\ B \ar[r, symbol = \hookrightarrow] & Y_{\ell}
        \end{tikzcd}
    \]
    such that $f$ is an isomorphism over the complement of $B \subset Y_{\ell}$.
\end{example}

\paragraph{\textbf{Comparison with existing moduli of sheaves in the literature.}}
Unlike the classical moduli of Gieseker semistable sheaves constructed by Simpson \cite{Simpson-repnI} and Langer \cite{langer-boundedness-general}, our moduli problem in \Cref{thm: main theorem intro 1} includes non-pure coherent sheaves.
In this context, the Gieseker semistable locus is the open substack of $\nps_{\alpha, \ell}$ parameterizing ideal sheaves $\cI_{Z}$ of length $\ell$ subschemes $[Z \subset X] \in \Hilb^{\ell}(X)$.
So $\nps_{\alpha, \ell}$ is strictly larger than the Gieseker semistable locus.

On the other hand, unlike the moduli of complexes subject to a Bridgeland stability condition \cite{bridgeland_stability}, our stack $\nps_{\alpha, \ell}$ always parameterizes coherent sheaves (as opposed to other objects in the derived category). We provide some examples where $\nps_{\alpha, \ell}$ can be recovered as a stack of Bridgeland semistable objects (\Cref{example: contractions of hilbert scheme of projective plane}), and also provide an example where the algebraic space $\Hilb^{\ell}(X)_{\alpha}$ is provably not recoverable as a moduli of Bridgeland semistable objects (\Cref{example: mukai flop contractions}).

Finally, in the context of \Cref{subsection: contraction of projective fibrations} for $\ell=1$, the moduli of sheaves $\nps_{\alpha, 1}^{\circ}$ described in \Cref{thm: existence of good moduli space other side of the wall} and \Cref{example: surgery of P ell fibrations} can be related to the moduli of perverse ideal sheaves constructed by Bridgeland in \cite[\S 5]{bridgeland_flops} whenever the corresponding contraction has rational singularities.

\paragraph{\textbf{Notation and conventions.}}
We work over an arbitrary ground field $k$. Given two morphisms of schemes $Y \to S$ and $T \to S$, we may sometimes use the notation $Y_T: = Y\times_S T$ to denote the base-change.
In this case, given a quasicoherent sheaf $F$ on $Y$, we denote by $F_T$ its pullback to $Y_T$. 

For the entirety of the paper, we work with a fixed smooth, projective, geometrically connected variety $X$ over $k$. We assume throughout that $\dim(X) \geq 2$. We also fix the choice of an ample line bundle $\cO_X(1)$. For convenience, we include here some common notation in alphabetical order:
\begin{description}[style=multiline, leftmargin=17ex, labelwidth=16ex, align=left, itemsep=0pt]
    \item[$\SCoh(X)$]
    Stack of coherent sheaves on $X$.
    \item[$\SCoh^{\rtf}(X)_{\cO_X, P_{\ell}}$]
    Stack of torsion-free coherent sheaves with trivial determinant and Hilbert polynomial $P_\ell$.
    \item[$\Hilb^\ell(X)$]
    The Hilbert scheme of finite length $\ell$ subschemes of $X$.
    \item[$\Hilb^\ell(X)_\alpha$]
    The good moduli space of $\cU_{\alpha,\ell}$ (see \cref{thm: existence of proper}).
    \item[$\Hilb^\ell(X)^{\rnc\alpha}$]
    The open subscheme of $\Hilb^\ell(X)$ on which the morphism $\Hilb^\ell(X) \to \Hilb^\ell(X)_\alpha$ is an isomorphism (see \cref{notation: nc-alpha}).
    \item[$\Hilb^\ell(X)_\alpha^\circ$]
    The good moduli space of $\cU_{\alpha,\ell}^\circ$ (see \cref{thm: existence of good moduli space other side of the wall}).
    \item[$\mathfrak{o}$] 
    The unique closed point of $\Theta_R$ or $\ST_R$ for a DVR $R$, also the unique closed point of $\cX_K$ for a field $K$.
    \item[$P_F$] 
    Hilbert polynomial of a sheaf $F$ on $X_K$ with respect to $\cO_X(1)$ for any field extension $K \supset k$.
    \item[$P_{\ell}$]
    Given a positive integer $\ell$, we set $P_{\ell}(t):= P_{\cO_X}(t) - \ell$.
    \item[$\ST_R$]
    Given a discrete valuation ring $R$ with uniformizer $\pi$, we define $\ST_R := \left[ \Spec(R[s,t]/(st-\pi)) / \bG_m \right]$, where $\mathbb G_m$ acts with weight $1$ on $s$ and weight $-1$ on $t$. 
    \item[$\cU_{\alpha,\ell}$]
    A certain enlargement of $\SCoh^{\rtf}(X)_{\cO_X, P_{\ell}}$ (see \cref{defn: stack U_alpha}).
    \item[$\cU_{\alpha,\ell}^\circ$]
    A certain open substack of $\cU_{\alpha,\ell}$ (see \cref{defn: other side of the wall}).
    \item[$\cX$]
    Define $\cX := \left[ \Spec(\bZ[s,t]/st) / \bG_m \right]$, where $\mathbb G_m$ acts with weight $1$ on $s$ and weight $-1$ on $t$.
    \item[$\Theta_R$]
    Given a ring $R$, define $\Theta_R := [\mathbb A^1_R / \mathbb G_m]$, where $\mathbb G_m$ acts with weight $-1$. 
\end{description}

\paragraph{\textbf{Acknowledgements.}} We are particularly grateful to Yoonjoo Kim, who greatly influenced and inspired the applications of the main result. Special thanks are also due to Daniel Halpern-Leistner, Mirko Mauri,  and James McKernan, who suggested some of the applications to derived categories, contractions, and the construction of flips/flops, respectively.
We also thank Anand Deopurkar, Robert Lazarsfeld, Emanuele Macr\`i, Joaqu\'in Moraga, Tony Pantev, Franco Rota and Chenyang Xu for useful discussions related to the contents of this paper.

We acknowledge the MSRVP program at the Mathematical Sciences Institute, which supported a visit of A.F.H. to Australian National University.
We also thank the Summer Research Institute in Algebraic Geometry 2025 (SRI) and Colorado State University, where the SRI was held, for providing a helpful platform for collaboration.

\paragraph{\textbf{AI Declaration.}} No AI input was used for any of the mathematical arguments or the writing of the first version of the article.

After the original submission, the only use of AI throughout the revision process was asking a free version of Claude to give us a list of typos after we received the first referee report. All of the revision replies and modifications were done by us with no AI input.

% ---
\section{Statement of the main result}

In this section, we start by introducing notation for some standard moduli stacks and giving a proof of the folklore result that the stack of torsion-free coherent sheaves with trivial determinant and Hilbert polynomial of an ideal sheaf of points is a trivial $\bG_m$-gerbe over the corresponding Hilbert scheme of points (\cref{prop: hilbert scheme as moduli space of sheaves}).
We then define the stack $\cU_{\alpha,\ell}$ which plays a prominent role in the rest of the paper and state the first main theorem (\cref{thm: main theorem on U_alpha}).

% -
\subsection{Moduli of coherent sheaves} \label{subsection: moduli of coherent sheaves}

We denote by $\SCoh(X)$ the stack of coherent sheaves on $X$. This is an algebraic stack \cite[\href{https://stacks.math.columbia.edu/tag/09DS}{09DS}]{stacks-project} with affine diagonal \cite[\href{https://stacks.math.columbia.edu/tag/08K9}{08K9}]{stacks-project} and locally of finite type over $k$ \cite[\href{https://stacks.math.columbia.edu/tag/08KD}{08KD}, \href{https://stacks.math.columbia.edu/tag/0CMY}{0CMY}]{stacks-project}. For any choice of polynomial $P(t)$, we denote by $\SCoh(X)_P \subset \SCoh(X)$ the open and closed substack of $\SCoh(X)$ parameterizing coherent sheaves with Hilbert polynomial $P(t)$ with respect to the fixed polarization $\cO_X(1)$.

There is a determinant morphism $\det: \SCoh(X) \to \SPic(X)$ to the Picard stack of $X$ \cite[pg. 36,37]{huybrechts.lehn}.
Given a line bundle $L$ on $X$, we denote by $\SCoh(X)_{L}$ the closed substack of $\SCoh(X)$ that fits into the following Cartesian diagram
 \[
\begin{tikzcd}
 \SCoh(X)_{L} \ar[r]
  \ar[d] & \SCoh(X) \ar[d, "\det"] \\ B\mathbb{G}_m \ar[r, "{[L]}"] &
  \SPic(X),
\end{tikzcd}
\]
where the bottom horizontal arrow is the closed immersion of the residual gerbe of $\SPic(X)$ at $L$.
Alternatively, we may define $\SCoh(X)_{L}$ as the subfunctor of $\SCoh(X)$ that sends a $k$-scheme $T$ to the subgroupoid of $T$-families of sheaves $F$ in $\SCoh(X)(T)$ such that there is some line bundle $M$ on $T$ and an isomorphism of line bundles $\det(F) \cong L \boxtimes M$ on $X \times T$. Note that we only require the existence of $M$ and the isomorphism; they are not part of the data.

Given a fixed Hilbert polynomial $P(t)$, we set $\SCoh(X)_{L,P} := \SCoh(X)_L \cap \SCoh(X)_P$. The main goal of this paper is to study certain open substacks of $\SCoh(X)_{L, P_L - \ell}$ that admit good moduli spaces. Since there is an isomorphism $(-) \otimes L: \SCoh(X)_{\cO_X, P_{\cO_X} -\ell} \xrightarrow{\sim} \SCoh(X)_{L, P_L - \ell}$, we will mainly restrict ourselves to the case where $L= \cO_X$.

\begin{notn}
    Given a positive integer $\ell>0$, we set $P_{\ell}(t):= P_{\cO_X}(t) - \ell$.
\end{notn}

There is an open substack $\SCoh^{\rtf}(X)_{\cO_X, P_{\ell}} \subset \SCoh(X)_{\cO_X, P_{\ell}}$ parameterizing torsion-free sheaves \cite[Thm. 12.2.1(iii)]{egaiv}. Given an integer $\ell>0$, we denote by $\Hilb^{\ell}(X)$ the Hilbert scheme of zero-dimensional length $\ell$ subschemes of $X$. There is a morphism $f: \Hilb^{\ell}(X) \to \SCoh^{\rtf}(X)_{\cO_X, P_{\ell}}$ that sends a subscheme $Z$ to its ideal sheaf $\cI_Z$. The following statement is certainly well-known; we provide a proof since we couldn't locate a suitable reference. 

\begin{prop}
\label{prop: hilbert scheme as moduli space of sheaves}
    Suppose that $\dim(X)\geq 2$. Then there is an isomorphism $\psi: \Hilb^{\ell}(X) \times B\mathbb{G}_m  \xrightarrow{\sim} \SCoh^{\rtf}(X)_{\cO_X, P_{\ell}}$
    that fits into the following commutative diagram
 \[
\begin{tikzcd}
 \Hilb^{\ell}(X) \ar[d, "{(\id, \triv)}"]
  \ar[dr, "f"] & \\ \Hilb^{\ell}(X) \times B\mathbb{G}_m \ar[r, "\psi"] &
   \SCoh^{\rtf}(X)_{\cO_X, P_{\ell}},
\end{tikzcd}
\]
where $\triv: \Hilb^{\ell}(X) \to \Spec(k) \to B\mathbb{G}_m$ classifies the trivial line bundle. In particular, $\SCoh^{\rtf}(X)_{\cO_X, P_{\ell}}$ has good moduli space $\Hilb^{\ell}(X)$.
\end{prop}

\begin{proof}
    Recall that for any $k$-scheme $T$ the groupoid $\Hilb^{\ell}(X) \times B\mathbb{G}_m(T)$ parameterizes pairs $([W\subset X_T], \cL)$, where $[W\subset X_T] \in \Hilb^{\ell}(X)(T)$ is a $T$-flat family of subschemes and $\cL$ is a line bundle on $T$. Let $\psi: \Hilb^{\ell}(X) \times B\mathbb{G}_m  \to \SCoh^{\rtf}(X)_{\cO_X, P_{\ell}}$ be the morphism that at the level of $T$-points sends a pair $([W \subset X_T], \cL)$ to the $T$-family of sheaves $\cI_W \otimes \pi^*(\cL)$, where $\pi: X_T \to T$ is the structure morphism. It is evident that the diagram in the statement of the proposition is commutative, the only thing left to show is that $\psi$ is an isomorphism. We prove this by exhibiting an inverse morphism $\theta$.

    Let $F \in \SCoh^{\rtf}(X)_{\cO_X, P_{\ell}}$ be a $T$-family of torsion-free sheaves.
    Since $F$ is a family in $\SCoh(X)_{\cO_X}$, we have that $\det(F)$ is the pullback of a line bundle from $T$.
    Hence, by the projection formula, we have that $\pi_*(\det(F))$ is a line bundle on $T$ and the counit $g: \pi^* \pi_*(\det(F)) \xrightarrow{\sim} \det(F)$ is an isomorphism.
    We have an isomorphism $(F^{\vee})^{\vee} \cong \det(F)$, and its formation commutes with base-change on $T$, see \Cref{app: properties of determinants}.
    Consider the natural morphism $h:F \to (F^{\vee})^{\vee} \cong \det(F)$. For all $t \in T$, the base-change $h_t: F_t \to (F^{\vee})^{\vee}_t \cong ((F_t)^{\vee})^{\vee}$ is generically an isomorphism, and hence injective since $F_t$ is torsion-free.
    It follows by the slicing criterion for flatness \cite[\href{https://stacks.math.columbia.edu/tag/00ME}{00ME}]{stacks-project} that $h:F \to \det(F)$ is injective, and the cokernel $\det(F)/F$ is $T$-flat. Consider the composition 
    \[ i:F \otimes \pi^*\pi_*(\det(F))^{\vee} \xrightarrow{h \otimes id} \det(F) \otimes \pi^*\pi_*(\det(F))^{\vee} \xrightarrow{\sim} \cO_{X_T},\]
    where the last isomorphism is induced from the counit $g: \pi^* \pi_*(\det(F)) \xrightarrow{\sim} \det(F)$. The morphism $i$ is injective with $T$-flat cokernel, and it exhibits $F \otimes \pi^*\pi_*(\det(F))^{\vee}$ as the ideal sheaf of a $T$-flat family of subschemes $[W \subset X_T] \in \Hilb^{\ell}(X)(T)$. We set $\theta(F):= ([W \subset X_T], \pi_*(\det(F))^{\vee}) \in \Hilb^{\ell}(X) \times B\mathbb{G}_m(T)$. Then it follows from construction that $\theta$ is an inverse of the morphism $\psi$, as desired.
    \qedhere
\end{proof}

% -
\subsection{The substack \texorpdfstring{$\nps_{\alpha,\ell}$}{U alpha ell}}

\begin{setup} \label{setup: alpha}
    We fix an open and closed subscheme $\alpha \subset \Hilb(X)$ of the Hilbert scheme of $X$, and express it as a disjoint union $\alpha = \bigsqcup_{j \in I} \alpha_j$ of its connected components. We assume that for all algebraically closed field extensions $K \supset k$ and all subschemes $[Z \subset X_K] \in \alpha(K)$, we have either:
    \begin{enumerate}
        \item $Z$ is integral and normal, or
        \item $Z$ is integral of codimension $1$. 
    \end{enumerate}
\end{setup}

\begin{example}[Non-movable rational curve]
\label{example: nm rat curve 1}
    One helpful example to keep in mind is when we have a non-movable \(\mathbb P^1\) inside of \(X\), i.e. an open and closed point \(\alpha = \{[i:\mathbb P^1 \hookrightarrow X]\}\) in the Hilbert scheme.
    This example will be treated in full detail in \cref{subsec: contractions of rigid P1s}.
\end{example}

\begin{defn}
\label{defn: stack U_alpha}
Fix a natural number $\ell > 0$ and recall that $P_{\ell} = P_{\cO_X} - \ell$.
We define $\nps_{\alpha,\ell} \subset \SCoh(X)_{\cO_X, P_{\ell}}$ to be the subfunctor that sends a $k$-scheme $T$ to the groupoid of families $F: T \to \SCoh(X)_{\cO_X, P_{\ell}}$ such that for all geometric points $t:\Spec(K) \to T$ the corresponding sheaf $F_t$ on $X_K$ satisfies one of the following:
\begin{enumerate}
    \item
    \label{item: ideal type}
    \emph{(Ideal type)}
    $F_t \cong \cI_Z$ for a length $\ell$ finite subscheme $Z\subset X_K$, or
    \item
    \label{item: torsion type}
    \emph{(Torsion type)} $F_t$ fits into a short exact sequence
    \[ 0 \to i_* G \to F_t \to \cI_W \to 0,\]
    where $W \in \alpha_j(K)$ is a subscheme belonging to some connected component $\alpha_j \subset \alpha$, we have that $[i: Y \hookrightarrow X_K] \in \alpha_j(K)$ belongs to the same component, and $G$ is a rank $1$ torsion-free sheaf on $Y$.
\end{enumerate}
\end{defn}

We note that the locus of sheaves of ideal type in \(\cU_{\alpha,\ell}\) is open in $\SCoh(X)_{\cO_X, P_{\ell}}$, because it is equivalent to being torsion-free (see \Cref{prop: hilbert scheme as moduli space of sheaves}), which is an open condition \cite[Thm. 12.2.1(iii)]{egaiv}.

\begin{remark}
    To justify that \Cref{defn: stack U_alpha} \ref{item: torsion type} is reasonable as stated, we note that any two subschemes $W, W' \subset X_K$ such that $\cI_W \cong \cI_{W'}$ lie in the same connected component of $\Hilb(X)$. Indeed, if we consider the $K$-vector space $V = \Hom_{X_K}(\cI_W, \cO_{X_K})$ as an affine space over $K$, then there is a morphism $V \setminus \{0\} \to \Hilb(X)$ sending a nonzero homomorphism $\varphi: \cI_W \hookrightarrow \cO_{X_K}$ to the closed subscheme with structure sheaf given by the cokernel $\Coker(\varphi)$.
    In this case $[W \subset X_K]$ and $[W' \subset X_K]$ belong to the image of the connected scheme $V$, and so they lie in the same connected component of $\Hilb(X)$.
\end{remark}

\begin{defn} \label{remark: split and nonsplit torsion type}
    For some of our arguments it will be useful to specify whether, for torsion type sheaves as in \cref{defn: stack U_alpha} \ref{item: torsion type}, their corresponding short exact sequence splits or not.
    We will refer to these subcases as the \emph{split torsion type} and the \emph{nonsplit torsion type}, respectively.
\end{defn} 

\noindent
\begin{minipage}[t]{0.64\textwidth}
\begin{example}[Non-movable rational curve]
\label{example: nm rat curve 2}
    In the example of a non-movable rational curve \(\alpha = \{[i:\mathbb P^1 \hookrightarrow X]\}\) as in \cref{example: nm rat curve 1}, we can describe the torsion type sheaves that get added to \(\cU_{\alpha,1}(K)\).
    A rank 1 torsion-free sheaf \(G\) on \(\mathbb P^1_K\) that would produce a torsion type sheaf \(F\) in \(\cU_{\alpha,1}(K)\) should have Hilbert polynomial \(P_G(t) = P_1(t) - P_{\cI_{\mathbb P^1}}(t) = P_{\cO_{\mathbb P^1}}(t)-1\), and thus the only such is \(G \cong \cO_{\mathbb P^1}(-1)\).
    Therefore, sheaves \(F \in \cU_{\alpha,1}(K)\) of non-split torsion type (the red locus in \cref{figure: P1 non-symmetric}) are described by lines in \(\Ext^1_{X_K} (\cI_{\mathbb P^1} , i_* \cO_{\mathbb P^1} (-1))\).
    There is only one sheaf of split torsion type (the blue locus in \cref{figure: P1 non-symmetric}) up to isomorphism, namely \(\cI_{\mathbb P^1} \oplus i_*\cO_{\mathbb P^1} (-1)\).
\end{example}
\end{minipage}
\hfill
\begin{minipage}[t]{0.32\textwidth}
    \centering
    \vspace{0pt}
    \includegraphics[width=0.9\linewidth]{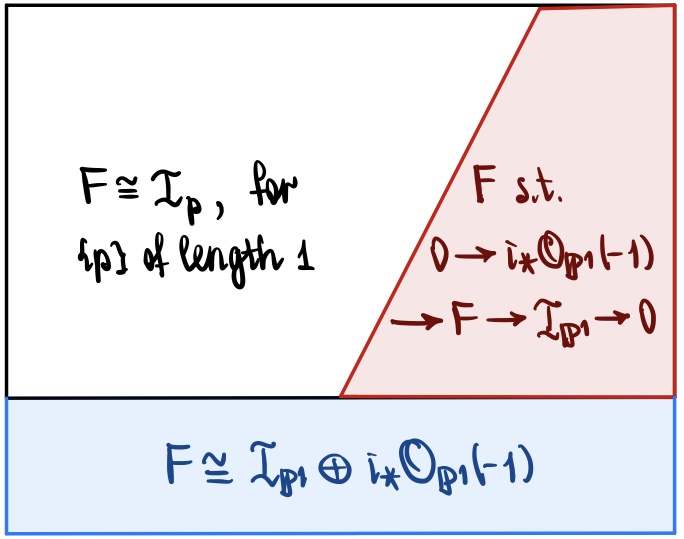}
    \captionof{figure}{The stack \(\cU_{\{[\mathbb P^1]\},1}\).}
    \label{figure: P1 non-symmetric}
\end{minipage}

\begin{defn} \label{defn: condition dagger}
    Let $\alpha = \bigsqcup_{j \in I} \alpha_j \subset \Hilb(X)$ be as in \cref{setup: alpha}, and let $\ell>0$ be a positive integer.
    We say that the pair $(\alpha, \ell)$ \emph{satisfies condition (LR)} if the following hold for all field extensions $K \supset k$:
    \begin{enumerate}
        \item For all $[W \subset X_K], [W' \subset X_K] \in \alpha(K)$, any nonzero homomorphism $f: \cI_{W'} \to \cI_{W}$ is an isomorphism;
        \item
        For all subschemes
        $[W \subset X_K], [i: Y \hookrightarrow X_K] \in \alpha(K)$,
        and all rank one torsion-free sheaves $G$ on $Y$ with Hilbert polynomial $P_{\cO_Y} - \ell$, we have $\Hom(\cI_W, i_*(G))=0$.
    \end{enumerate} 
\end{defn}

Condition (LR) in particular ensures that the closed points of the stack \(\cU_{\alpha,\ell}\) have linearly reductive stabilizers (see \cref{prop: local reductivity}), hence the choice of the label.

\begin{remark} \label{remark: condition 1 automatic dimension}
    Note that if all integral subschemes $[W \subset X_K] \in \alpha(K)$ have the same codimension, then part (1) in \cref{defn: condition dagger} is automatically satisfied. This will be the case in all of our applications.
\end{remark}

\begin{example}[Non-movable rational curve]
    Concluding \cref{example: nm rat curve 2} with \(\alpha = \{[i:\mathbb P^1 \hookrightarrow X]\}\), we observe that condition (LR) is equivalent to asking that \(\Hom(\cI_W, i_*(G)) \cong H^0(\mathbb P^1_K, \cN_{\mathbb P^1_K / X_K}(-1))\) vanishes, where \(\cN_{\mathbb P^1_K / X_K}\) is the normal bundle.
    Since all vector bundles on \(\mathbb P^1_K\) split into line bundles, this is the same as asking that the normal sheaf \(\cN_{\mathbb P^1_K / X_K} \cong \bigoplus_{i=1}^{\dim(X)-1} \cO(a_i)\), with all \(a_i \leq 0\).
    Furthermore, since \(\mathbb P^1\) is non-movable, at least one \(a_i\) must be strictly negative by deformation theory.
    We treat this example in bigger generality, allowing \(\alpha\) to be a family of rational curves and \(\ell \geq 1\), in \cref{thm: contraction of projective space bundles}.
\end{example}

The following is the main result of this article.
\begin{thm} \label{thm: main theorem on U_alpha}
    Let $\alpha$ and $\ell>0$ be as in \cref{setup: alpha}.
    Then the following hold:
    \begin{enumerate}
        \item
        \label{item: main theorem: openness}
        The substack $\nps_{\alpha,\ell} \subset \SCoh(X)_{\cO_X, P_{\ell}}$ is an open substack.
        \item
        \label{item: main theorem: existence of gms}
        If $\alpha$ is of finite type and $(\alpha, \ell)$ satisfies condition (LR) in the sense of \cref{defn: condition dagger}, then $\nps_{\alpha,\ell}$ admits a proper good moduli space $\Hilb^{\ell}(X)_{\alpha}$.
        \item 
        \label{item: main theorem: morhism from Hilbert scheme}
        In the setting of part \ref{item: main theorem: existence of gms}, there is a morphism $\Hilb^\ell(X) \to \Hilb^\ell(X)_\alpha$.
    \end{enumerate}
\end{thm}

\begin{proof}
Openness of $\nps_{\alpha, \ell}$ is proven in \cref{thm: openness general case}.
The existence of a good moduli space and its properness is proven in \cref{thm: existence of proper}.
Part \ref{item: main theorem: morhism from Hilbert scheme} follows from the inclusion $\SCoh^{\rtf}(X)_{\cO_X,P_\ell} \subset \cU_{\alpha,\ell}$, the fact that $\Hilb^\ell(X)$ is the good moduli space for $\SCoh^{\rtf}(X)_{\cO_X,P_\ell}$ (\cref{prop: hilbert scheme as moduli space of sheaves}), and universality of good moduli spaces \cite[Thm. 6.6]{alper-good-moduli}.
\end{proof}

% ---
\section{Openness and boundedness of the stack} \label{section: openness and boundedness of the stack}

In order to ensure that the stack $\cU_{\alpha,\ell}$ is algebraic, we show that it is open inside the algebraic stack $\SCoh(X)_{\cO_X,P_\ell}$.
For that, we start by proving that it is constructible (\cref{prop: constructibility of the stack}), then move on to proving that it is closed under generalization (\cref{thm: openness general case}).
Furthermore, we show that the stack $\cU_{\alpha,\ell}$ is quasi-compact under certain conditions on $(\alpha,\ell)$ (\cref{prop: boundedness}).

% -
\subsection{Constructibility}

\begin{prop}
\label{prop: constructibility of the stack}
    Let $\alpha$ be as in \cref{setup: alpha}. Let $|\nps_{\alpha, \ell}| \subset |\SCoh(X)_{\cO_X, P_{\ell}}|$ denote the subspace of the underlying topological space of $\SCoh(X)_{\cO_X, P_{\ell}}$ that consists of geometric points lying in the substack $\nps_{\alpha, \ell}$ (\cref{defn: stack U_alpha}). Then, $|\nps_{\alpha, \ell}|$ is a locally constructible subspace of the locally Noetherian topological space $|\SCoh(X)_{\cO_X, P_{\ell}}|$.
\end{prop}

\begin{proof}
    It suffices to show that for all Noetherian $k$-schemes $T$ and all families of sheaves $F: T \to \SCoh(X)_{\cO_X, P_{\ell}}$, the locus of points $t \in T$ such that $F_t$ belongs to $\nps_{\alpha, \ell}$ is constructible. By Noetherian induction, it suffices to restrict to the case when $T$ is integral with generic point $\eta$ and show that both of the following hold:
    \begin{enumerate}
        \item[(i)] if $F_{\eta}$ is contained in $\nps_{\alpha, \ell}$, then there is an open subscheme $U \subset T$ such that $F_u$ is contained in $\nps_{\alpha, \ell}$ for all $u \in U$, and
        \item[(ii)] if $F_{\eta}$ is not contained in $\nps_{\alpha, \ell}$, then there is an open subscheme $U \subset T$ such that $F_u$ is not contained in $\nps_{\alpha, \ell}$ for all $u \in U$.
    \end{enumerate} 
    
    Let's prove part (i). Suppose that $F_{\eta}$ is contained in $\nps_{\alpha, \ell}$. If $F_{\eta}$ is torsion-free, then by \cite[Thm. 12.2.1(iii)]{egaiv} there exists some open $U \subset T$ such that $F_u$ is torsion-free for all $u \in U$.
    Since $\SCoh^{\rtf}(X)_{\cO_X, P_{\ell}}$ is contained in $\nps_{\alpha, \ell}$ (see \Cref{prop: hilbert scheme as moduli space of sheaves}), it follows that $F_u$ is in $\nps_{\alpha, \ell}$ for all $u \in U$. 
    
    Otherwise, if $F_{\eta}$ is of torsion type, then there exist geometrically integral subschemes $[i: Y \hookrightarrow X_{\eta}], [W \subset X_{\eta}] \in \alpha_j(\eta)$ and a rank one torsion-free sheaf $G$ on $Y$ such that $F_{\eta}$ fits into a short exact sequence
    \[ 0 \to i_*(G) \to F_{\eta} \to \cI_W \to 0.\]
    By spreading out, there is an open subscheme $U \subset T$ and closed subschemes $[\widetilde{i}: \widetilde{Y} \hookrightarrow X_U], [\widetilde{W} \subset X_U] \in \alpha_j(U)$ such that $\widetilde{Y}_{\eta} = Y$ and $\widetilde{W}_{\eta} = W$. After possibly shrinking $U$, we may spread out $G$ to a $U$-family of rank one torsion-free sheaves $\widetilde{G}$ on $\widetilde{Y} \to U$. Moreover, after further shrinking $U$, we may spread out the morphism $F_{\eta} \twoheadrightarrow\cI_W$ to a surjection $F_U \twoheadrightarrow \cI_{\widetilde{W}}$ of sheaves on $X_U$. If we denote the kernel of this surjection by $K$, then after shrinking $U$ we may spread out the isomorphism of generic fibers $(\widetilde{i}_*(\widetilde{G}))_{\eta} = i_*(G) \cong K_{\eta}$ to an isomorphism $\widetilde{i}_*(\widetilde{G}) \cong K$. We conclude that $F_U$ fits into a short exact sequence
    \[ 0 \to \widetilde{i}_*(\widetilde{G}) \to F_U \to \cI_{\widetilde{W}} \to 0,\]
    and it follows that $F_u$ belongs to $\nps_{\alpha, \ell}$ for all $u \in U$. This concludes the proof of (i).

    We are left to show (ii). Suppose that $F_{\eta}$ is not in $\nps_{\alpha, \ell}$. This means that $F_{\eta}$ is not torsion-free. Consider the torsion filtration
    \begin{equation} \label{equation: proof of constructibility torsion filtration}
        0 \to i_*(G) \to F_{\eta} \to Q \to 0,
    \end{equation}
    where $i: Y \hookrightarrow X_{\eta}$ is the scheme-theoretic support of the maximal torsion subsheaf of $F_{\eta}$, and $Q$ denotes the rank one torsion-free quotient.
    Note that $\det(Q)^{\vee} \cong \det(i_*(G))$ is effective (see \cref{lemma: det is an effective Cartier divisor}), and hence $\det(Q)$ is the ideal sheaf of a subscheme of $X_{\eta}$.
    Using the inclusion $Q \hookrightarrow Q^{\vee \vee} \cong \det(Q)$ (see \cref{lemma: det(rank 1 tf) = double dual}), we conclude that $Q= \cI_W$ is the ideal sheaf of a closed subscheme $W \subset X_{\eta}$. Over a small enough open subscheme $U \subset T$, we may spread out the short exact sequence \eqref{equation: proof of constructibility torsion filtration} to
    \[ 0 \to \widetilde{i}_*(\widetilde{G}) \to F_U \to \cI_{\widetilde{W}} \to 0,\]
    where $\widetilde{Y} \hookrightarrow X_U$ and $\widetilde{W} \hookrightarrow X_U$ are $U$-flat families of subschemes, and $\widetilde{G}$ is a $U$-flat family of sheaves on $\widetilde{Y} \to U$. Note that the torsion-freeness of the family of sheaves $\cI_{\widetilde{W}}$ implies that for all $u \in U$ the restriction of this short exact sequence to $X_u$ recovers the torsion filtration of $F_u$. We may further assume, after perhaps shrinking $U$, that $\widetilde{Y}$ agrees with the scheme-theoretic support of $\widetilde{i}_*(\widetilde{G})$. Hence for all $u \in U$ the subscheme $\widetilde{Y}_u \subset X_u$ agrees set-theoretically with the support of the torsion $\widetilde{i}_*(\widetilde{G})_u$ of $F_u$.
    Since $F_{\eta}$ is not in $\nps_{\alpha, \ell}$ one of the following must hold:
    \begin{enumerate}
        \item[(a)] $[W \subset X_{\eta}]$ does not belong to $\alpha(\eta)$.
        \item[(b)] $[W \subset X_{\eta}] \in \alpha_j(\eta)$ but $[i: Y \hookrightarrow X_{\eta}]$ does not belong to $\alpha_j(\eta)$.
        \item[(c)] $G$ is not torsion-free of rank one on $Y$.
    \end{enumerate}
    Case (a). Suppose that $W \subset X_{\eta}$ does not belong to $\alpha(\eta)$. Consider the morphism $U \xrightarrow{\widetilde{W}} \Hilb(X)$. Since $\alpha$ is open and closed and $U$ is irreducible, it follows that for all $u \in U$ we have that $\widetilde{W}_u$ is not in $\alpha(u)$. Note that any two subschemes with isomorphic ideal sheaves belong to the same component of the Hilbert scheme of $X$, and so it follows that the torsion-free quotient $\cI_{\widetilde{W}_u}$ of $F_u$ is not isomorphic to the ideal sheaf of any element in $\alpha(u)$. Hence $F_u$ is not in $\nps_{\alpha, \ell}$ for all $u \in U$, as desired.
    
    Case (b). Suppose that $[W \subset X_{\eta}] \in \alpha_j(\eta)$ but $Y$ does not belong to $\alpha_j(\eta)$. If $Y$ is not geometrically irreducible, then by constructibility of the geometrically irreducible locus \cite[\href{https://stacks.math.columbia.edu/tag/055B}{055B}]{stacks-project}, there is a sufficiently small open subset of points $u$ such that the set-theoretic support $\widetilde{Y}_u$ of the torsion of $F_u$ is not geometrically irreducible, and hence $F_u$ does not belong to $\nps_{\alpha,\ell}$. 
    
    If $Y$ is not geometrically reduced, then after passing to a finite flat cover of $U$ we may assume that the reduced subscheme $r: Y_{red} \hookrightarrow Y$ is geometrically reduced, and we have that the natural morphism $G \to r_*r^*(G)$ is not an isomorphism, since the scheme-theoretic support of $G$ is $Y$. After spreading out the closed immersion $\widetilde{r}: \widetilde{Y}_{red} \hookrightarrow \widetilde{Y}$ and perhaps shrinking $U$, we see that for all $u \in U$ the natural morphism $\widetilde{G}_u \to (r_u)_*r_u^*(\widetilde{G}_u)$ is not an isomorphism. We conclude that the scheme-theoretic support of the torsion $(\widetilde{i})_*(\widetilde{G})_u \subset F_u$ is not geometrically reduced, and hence $F_u$ is not in $\nps_{\alpha, \ell}$ for all $u \in U$. 
    
    We are left to consider the case when $Y$ is geometrically integral.
    By constructibility of the geometrically integral locus \cite[\href{https://stacks.math.columbia.edu/tag/0579}{0579}, \href{https://stacks.math.columbia.edu/tag/055B}{055B}]{stacks-project}, we may assume after shrinking $U$ that $\widetilde{Y}_u$ is geometrically integral for all $u \in U$ . Since $\alpha_j \subset \Hilb(X)$ is closed, there is some open subset of points $u$ such that the geometrically integral subscheme $\widetilde{Y}_u \subset X_u$ is not in $\alpha_j(u)$. We conclude that the irreducible set-theoretic support of the torsion $\widetilde{i}_*(\widetilde{G})_u$ of $F_u$ is not in $\alpha_j$, which means that $F_u$ cannot be in $\nps_{\alpha, \ell}$.

    Case (c). Assume that $Y$ belongs to $\alpha_j(\eta)$, but that $G$ is not torsion-free of rank one in $Y$. As in the previous paragraph, we may assume that $\widetilde{Y}$ is a family of geometrically integral subschemes by \cite[\href{https://stacks.math.columbia.edu/tag/0579}{0579}, \href{https://stacks.math.columbia.edu/tag/055B}{055B}]{stacks-project}. If $G$ has torsion, then by spreading out the torsion filtration we see that $\widetilde{G}_u$ has torsion for all $u$ in an open of $T$, thus implying that $F_u$ is not in $\nps_{\alpha, \ell}$ for all such $u$. The other possibility is that $G$ does not have rank one. This is equivalent to the first Fitting support $\Fit_1(G) \subset Y$ being strictly smaller than the (integral) support $Y$ of $G$. After shrinking $U$, we may assume that for all $u \in U$ the dimension of $\Fit_1(\widetilde{G})_u = \Fit_1(\widetilde{G}_u)$ is strictly smaller than $\dim(\widetilde{Y}_u) = \dim(Y)$. We conclude that for all $u \in U$ we have that $\widetilde{G}_u$ does not have rank one, and hence $F_u$ does not belong to $\nps_{\alpha, \ell}$.
\end{proof}

\subsection{Openness}

To conclude the proof of openness of the substack $\nps_{\alpha, \ell} \subset \SCoh(X)_{\cO_X, P_{\ell}}$, in this subsection we show that it is closed under generalization.

\begin{construction}
\label{construction: limit in Quot}
    We fix $\alpha \subset \Hilb(X)$ as in \cref{setup: alpha}.
    Let $R$ be a discrete valuation ring, and let $\eta$ (resp. $s$) denote the generic (resp. special) point of $\Spec(R)$.
    Let $F$ be an $R$-flat family of coherent sheaves on $X_R$ with trivial determinant such that the special fiber $F_s \in \nps_{\alpha,\ell}(s)$ and is of torsion type, i.e. it fits into a short exact sequence
    \begin{equation}
        0 \to i_{s,*} (G_s) \to F_s \to \cI_{W_s} \to 0
    \end{equation}
    for some $[i_s : Y_s \to X_s], [W_s \subset X_s] \in \alpha_j(s)$, and a rank one torsion-free sheaf $G_s$ on $Y_s$.
    Assume that $F_\eta$ is not torsion-free.
    Take the torsion sequence
    \begin{equation}
        \label{eq: construction: torsion sequence for F_eta}
        0 \to \Tors(F_\eta) \to F_\eta \to \TF(F_\eta) \to 0.
    \end{equation}
    Consider the unique extension $T \subset F$ of the subsheaf $\Tors(F_\eta) \subset F_{\eta}$ fitting into a short exact sequence
    \begin{equation}
        \label{eq: construction: limit in the Quot scheme}
        0 \to T \to F \to E \to 0
    \end{equation}
    such that $T$ and $E$ are $R$-flat (that is, $T$ is the limit in the corresponding Quot scheme).
\end{construction}

\begin{lemma}
\label{lemma: T is pure}
    Assume the assumptions and notation from \cref{construction: limit in Quot}.
    Then $T$ is a family of pure sheaves on $X_R \to \Spec(R)$, and $T_s \cong i_{s,*}(A_s)$ for some rank one torsion-free subsheaf $A_s \subset G_s$.
\end{lemma}

\begin{proof}
    Restrict \eqref{eq: construction: limit in the Quot scheme} to the fiber over the closed point $X_s$ and compare the result with the torsion sequence for $F_s$:
    % https://q.uiver.app/#q=WzAsMTAsWzAsMCwiMCJdLFsyLDAsIlRfcyJdLFs0LDAsIkZfcyJdLFs2LDAsIkVfcyJdLFs4LDAsIjAiXSxbMCwyLCIwIl0sWzIsMiwiXFxUb3JzKEZfcykiXSxbNCwyLCJGX3MiXSxbNiwyLCJcXFRGKEZfcykiXSxbOCwyLCIwIl0sWzAsMV0sWzEsMl0sWzIsM10sWzMsNF0sWzUsNl0sWzYsN10sWzcsOF0sWzgsOV0sWzEsNiwiIiwxLHsic3R5bGUiOnsidGFpbCI6eyJuYW1lIjoiaG9vayIsInNpZGUiOiJ0b3AifX19XSxbMiw3LCI9Il0sWzMsOCwiIiwxLHsic3R5bGUiOnsiaGVhZCI6eyJuYW1lIjoiZXBpIn19fV1d
    \begin{equation}
    \label{eq: diagram of torsion sequences}
    \begin{tikzcd}[ampersand replacement = 
    \&]
    	0 \&\& {T_s} \&\& {F_s} \&\& {E_s} \&\& 0 \\
    	\\
    	0 \&\& {i_{s,*}(G_s)} \&\& {F_s} \&\& {\cI_{W_s}} \&\& 0
    	\arrow[from=1-1, to=1-3]
    	\arrow[from=1-3, to=1-5]
    	\arrow[hook, from=1-3, to=3-3]
    	\arrow[from=1-5, to=1-7]
    	\arrow["{=}", from=1-5, to=3-5]
    	\arrow[from=1-7, to=1-9]
    	\arrow[two heads, from=1-7, to=3-7]
    	\arrow[from=3-1, to=3-3]
    	\arrow[from=3-3, to=3-5]
    	\arrow[from=3-5, to=3-7]
    	\arrow[from=3-7, to=3-9]
    \end{tikzcd}
    \end{equation}
    Note that $T_s$ is a torsion sheaf and it is contained in $\Tors(F_s)\cong i_{s,*}(G_s)$.
    Then it follows from $T_s \subset i_{s,*} (G_s)$ that $T_s$ is a pushforward as well, i.e. $T_s \cong i_{s,*}(A_s)$ for some rank one torsion-free subsheaf $A_s \subset G_s$.
    In particular, $T_s$ is pure.
    Since being pure is open in families, we conclude the desired statement.
\end{proof}

\begin{lemma}
\label{lemma: E eta is an ideal sheaf}
    With the assumptions and notation from \cref{construction: limit in Quot}, the following hold:
    \begin{enumerate}
        \item 
        \label{item: T s is the max torsion of F s}
        The special fiber $T_s$ of $T$ is the maximal torsion subsheaf of $F_s$, so $E_s$ is torsion-free.
        \item 
        \label{item: E eta is the ideal of W eta from alpha j}
        $E_{\eta}$ is isomorphic to an ideal sheaf $\cI_{V}$ for some $V \in \alpha_j(\eta)$.
    \end{enumerate}
\end{lemma}

\begin{proof}
    The proof of this lemma proceeds differently depending on whether elements of $\alpha_j$ are divisorial.
    In both cases, we shall use the following short exact sequence that we obtain from \cref{lemma: T is pure} using the third isomorphism theorem:
    \begin{equation}
    \label{eq: third iso theorem}
        0 \to i_{s,*}(G_s/A_s) \to E_s \to \cI_{W_s} \to 0.
    \end{equation}
    Since $\rk(A_s) = \rk(G_s) = 1$ on $Y_s$, we have that $\codim\left( \Supp \left( i_{s,*}(G_s/A_s) \right)\right) > \codim(Y_s) = \codim(W_s)$.

    \paragraph{\textbf{Divisorial case.}}
    Suppose that $\codim (W_s) = \codim (Y_s) = 1$.
    We subdivide the proof into several steps.

  \noindent \textbf{Step 1. $E_\eta$ is isomorphic to some ideal sheaf $\cI_V$.}

        By purity of $T_\eta$ from \cref{lemma: T is pure}, the scheme-theoretic support
        $Y' := \Supp(T_\eta)$ is pure of codimension $1$.     
        Since $F$ is a family of sheaves with trivial determinant,
        we get from multiplicativity of the determinant in short exact sequences (\cref{app: properties of determinants}) and \cref{construction: limit in Quot}
        that $\det (E_\eta) \cong \det (T_\eta)^\vee \subset \cO_{X_\eta} (-Y')$.
        Consequently, since $E_\eta$ is torsion-free of rank one, we have an injection $E_\eta \hookrightarrow \cO_{X_\eta} (-Y')$, and therefore we can realize this sheaf as $E_\eta \cong \cI_{V}$ for some closed subscheme $V \subset X_\eta$ with $V \supset Y'$.

        \noindent \textbf{Step 2. $P_{E_s} = P_{\cI_{W_s}}$ and $V$ is divisorial.}
        
        Using $V \subset X_\eta$ as in the previous step, we get the following chain of inequalities of Hilbert polynomials:
        \[
            P_{\cI_V}
            = P_{E_\eta}
            = P_{E_s} \geq
            P_{\cI_{W_s}},
        \]
        where the last inequality follows from surjectivity of $E_s \to \cI_{W_s}$ in \eqref{eq: diagram of torsion sequences}.
        On the other hand, the inclusion $E_\eta \subset \det(E_\eta)$ implies that 
        \[ P_{\cI_V} = P_{E_\eta} \leq P_{\det(E_\eta)} = P_{\det (E_s)} = P_{\cI_{W_s}},\]
        where the last equality follows from \eqref{eq: third iso theorem} and the fact that $\Supp \left( i_{s,*}(G_s/A_s) \right)$ has codimension at least $2$.
        We now have inequalities in both directions, which implies that $P_{E_s} = P_{\cI_{W_s}}$. Furthermore, we must also have $P_{\cI_V} = P_{\det(E_\eta)}$, and thus $\cI_V \cong E_\eta = \det(E_\eta)$ is a line bundle and $V$ is divisorial.
 
        \noindent \textbf{Step 3. $A_s = G_s$, hence $T_s$ is the maximal torsion subsheaf of $F_s$.}

        By the short exact sequences in \eqref{eq: diagram of torsion sequences} and the equality $P_{E_s} = P_{\cI_{W_s}}$ from Step 2, we get $P_{A_s}= P_{G_s}$. This implies that the inclusion $A_s \subseteq G_s$ from \cref{lemma: T is pure} is an equality, hence $T_s = i_{s,*}(G_s) = \Tors(F_s)$. We have thus proven part \ref{item: T s is the max torsion of F s} of this lemma.

        \noindent \textbf{Step 4. $V$ and $W_s$ are in the same component of $\Hilb(X)$.}

        Since $V$ is divisorial, we can realize it as the vanishing locus of some section $\sigma : \cO_{X_\eta} \to \cO_{X_\eta} (V) = \det (E_\eta)^{\vee}$.
        Multiplying $\sigma$ by a sufficiently high power of the uniformizer of $R$, we get a section of $\det (E)^{\vee}$ that possibly vanishes over the special point $s$; dividing the latter section by another power of the uniformizer yields a section $\varphi : \cO_{X_R} \to \det(E)^{\vee}$ that does not vanish over $s$.
        The vanishing locus of this section is the flat closure $\overline V \subset X_R$.

        From \eqref{eq: third iso theorem} we get that $\det(E_s)^{\vee} \cong \cO(W_s)$, so the restriction of $\varphi$ to $s$ yields a section of $\cO(W_s)$ whose vanishing locus is $\left( \overline V \right)_s$.
        This implies that $W_s$ and $\left( \overline V \right)_s$ are in the same component $\alpha_j$ of the Hilbert scheme of $X$.
        Consequently, $V \in \alpha_j(\eta)$, which concludes the proof of part \ref{item: E eta is the ideal of W eta from alpha j} in the divisorial case.

    \paragraph{\textbf{Case of codimension at least $2$.}}
    Suppose that $\codim (W_s) = \codim (Y_s) \geq 2$.
    Again, we proceed in steps.

    \noindent \textbf{Step 1. $E_\eta$ is isomorphic to some ideal sheaf $\cI_V$.}

        In this case, the support of $T_\eta$ has codimension at least $2$, hence we have $\det(E_\eta) \cong \det(F_\eta) \otimes \det(T_\eta)^\vee \cong \cO_{X_\eta}$,
        by \cref{lemma: det is an effective Cartier divisor}.
        Since $E_\eta$ is torsion-free and of rank $1$, we have $E_\eta \subset \det(E_\eta)$, and we conclude that $E_\eta\cong \cI_V$ for some closed subscheme $V \subset X_\eta$.

        \noindent \textbf{Step 2. $W_s$ is contained in $\left(\overline V\right)_s$.}

        Using $V \subset X_\eta$ as in the previous step, we can consider its closure $\overline V$ in $X_R$, which is an $R$-flat closed subscheme.
        Consider the maximal open subset $U \subset X_R$ where $E$ is locally free, and set $Z := X_R \setminus U$ with the reduced subscheme structure.
        By \cite[Lemma \href{https://stacks.math.columbia.edu/tag/0CZR}{0CZR}]{stacks-project}, the fiber $Z_\eta$ is precisely the locus where $E_\eta \cong \cI_V$ is not locally free, hence $Z_\eta = V$ set-theoretically.
        Taking the closure of $V$ in $X_R$ yields $\overline V \subset Z$, because $Z$ is closed, hence we get $\left(\overline V\right)_s \subset Z_s$ at the level of sets.
        
        Similarly applying \cite[Lemma \href{https://stacks.math.columbia.edu/tag/0CZR}{0CZR}]{stacks-project} and using \eqref{eq: third iso theorem}, we get the equality $Z_s = \Supp(G_s/A_s) \cup W_s$ as sets. 
        Combining the above, we get $\left(\overline V\right)_s \subset \Supp(G_s/A_s) \cup W_s$.
        Note that, by \eqref{eq: third iso theorem}, we have
        \begin{equation} \label{eqn: containment inequality hilbert}
            P_{\cI_V} = P_{E_{\eta}} = P_{E_s} = P_{\cI_{W_s}} + P_{ i_{s,*}(G_s/A_s)}.
        \end{equation}
        Since $\codim\left( \Supp \left( i_{s,*}(G_s/A_s) \right)\right) > \codim(W_s)$, the equality \eqref{eqn: containment inequality hilbert} above implies that $\codim(V) = \codim(W_s)$. Since $\codim\left(\left(\overline V\right)_s\right) = \codim(V) = \codim(W_s)$, while $\Supp(G_s/A_s)$ has strictly larger codimension, the containment $\left(\overline V\right)_s \subset \Supp(G_s/A_s) \cup W_s$ and the integrality of $W_s$ force that $W_s \subset \left(\overline V\right)_s$ at the level of schemes.

        \noindent
        \textbf{Step 3. $\left(\overline V\right)_s = W_s$, hence $V$ is contained in $\alpha_j(\eta)$.}

        By \eqref{eqn: containment inequality hilbert}, we have $P_{\cI_{\left(\overline V\right)_s}} = P_{\cI_V} \geq P_{\cI_{W_s}}$, or equivalently,
        $P_{\left(\overline V\right)_s} \leq P_{W_s}$.
        On the other hand, we have just proved that $W_s \subset \left(\overline V\right)_s$. This forces $\left(\overline V\right)_s = W_s$.
        Since $\overline V$ is an $R$-flat family whose special fiber is in $\alpha_j$, we conclude that its generic fiber $V$ is in $\alpha_j$ as well.

        \noindent \textbf{Step 4. $A_s = G_s$, hence $T_s$ is the maximal torsion subsheaf of $F_s$.}
        
        The proof of $A_s = G_s$ is the same as in Step 3 of the divisorial case. \qedhere
\end{proof}

\begin{lemma}
\label{lemma: pure sheaves with support in beta}
    Let $\alpha \subset \Hilb(X)$ as in \Cref{setup: alpha}. Let $T$ be an $R$-flat family of pure sheaves on $X_R \to \Spec(R)$. Suppose that the special fiber $T_s$ is of the form $i_*(G)$ for some $[i : Y \hookrightarrow X_s] \in \alpha_j(s)$ and some rank one torsion-free sheaf $G$ on $Y$. Then, the generic fiber $T_{\eta}$ is of the form $(i')_*(G')$ for some $[i' : Y' \hookrightarrow X_{\eta}] \in \alpha_j(\eta)$ and some rank one torsion-free sheaf $G'$ on $Y'$.
\end{lemma}
\begin{proof}
Since the special fiber $T_s= i_*(G)$ is the pushforward of a rank $1$ torsion-free sheaf on the geometrically integral subscheme $Y \subset X_s$, it follows that the coherent sheaf $T_s$ on $X_s$ is geometrically integral in the sense of \cite[Defn. 4.6.19]{egaiv_2}. By \cite[Thm. 12.2.1(x)]{egaiv}, the generic fiber $T_{\eta}$ is also geometrically integral.
It follows that the scheme-theoretic support $i': \Supp(T_{\eta}) \hookrightarrow X_{\eta}$ is geometrically integral and $T_{\eta} = (i')_*(G')$ for some rank $1$ torsion-free sheaf $G'$ on $Y' := \Supp(T_{\eta})$.
To conclude, it suffices to show that the integral subscheme $Y'$ belongs to the same connected component $\alpha_j$ of \(\alpha \subset \Hilb(X)\) as $Y$.
If we set $\widetilde{Y} := \overline{Y'}\subset X_R$, then $[\widetilde{Y} \hookrightarrow X_R]$ corresponds to a morphism $\Spec(R) \to \Hilb(X)$ that maps the generic point $\eta$ to $[Y' \subset X_{\eta}]$. To conclude, it suffices to show that the special fiber $\widetilde{Y}_s$ recovers $Y$.
Note that we have $\widetilde{Y} \subset \Supp(T)_{red}$, and passing to special fibers we have that $\widetilde{Y}_s$ is contained in $Y = \Supp(T_s)$ set-theoretically. Since $Y$ is irreducible of the same dimension as $\widetilde{Y}_s$, it follows that the desired equality $Y = \widetilde{Y}_s$ holds set-theoretically. To conclude, we just need to show that $\widetilde{Y}_s$ is reduced.

We treat the two cases from \cref{setup: alpha} separately, depending on whether the integral subscheme $[i: Y \hookrightarrow X_s]$ is normal or a divisor.

\noindent \textbf{Case 1.}
Suppose that $Y$ is normal. 
Let $\Fit_0(T) \subset X_R$ denote the zeroth Fitting support of $T$, whose formation commutes with arbitrary base-change.
Since $Y' = \Supp(T_{\eta}) \subset \Fit_0(T_{\eta}) = \Fit_0(T)_{\eta}$, it follows that $\widetilde{Y} = \overline{Y'} \subset \Fit_0(T)$.
Passing to special fibers, we have that $\widetilde{Y}_s \subset \Fit_0(T_s)$.
Note that the sheaf $T_s = i_*(G)$ is isomorphic to $\cO_Y$ over an open dense subset of its support.
Hence, we have that $\Fit_0(T_s)$, which set-theoretically agrees with $Y$, also agrees with $Y$ scheme-theoretically over an open dense subset.
Hence, the same holds for $\widetilde{Y}_s \subset \Fit_0(T_s)$, and so $\widetilde{Y}_s$ is generically reduced. Since the reduction $(\widetilde{Y}_s)_{red} =Y$ is normal, it follows by Hironaka's lemma (see \cite[Cor. 3]{kollar_flatness}) that $\widetilde{Y}_s$ is reduced.
This concludes the proof in the case when $Y$ is geometrically normal.

\noindent \textbf{Case 2.}
Suppose that $Y \subset X_s$ is a divisor. Consider the determinant $D:=\det(T)$, whose formation commutes with restriction to fibers over \(\Spec(R)\).
We have that $D$ is an effective Cartier divisor in $X_R$ (see \Cref{lemma: det is an effective Cartier divisor}), and its special fiber $D_s = \det(i_*(G))$ agrees with the integral divisor $Y \subset X_s$.
Similarly, the generic fiber $D_{\eta} = \det((i')_*(G'))$ agrees with the integral divisor $Y' \subset X_{\eta}$. Since $D_{\eta} = Y'$, we have $\widetilde{Y} = \overline{Y'} \subset D$. Taking special fibers, we have $\widetilde{Y}_s \subset D_s = Y$. Since $(\widetilde{Y}_s)_{red} = Y$, this forces $\widetilde{Y}_s = Y$, as desired.
\end{proof}

\begin{thm} \label{thm: openness general case}
    Let $\alpha \subset \Hilb(X)$ as in \cref{setup: alpha}.
    Then, the subfunctor $\nps_{\alpha, \ell}\subset \SCoh(X)_{\cO_X, P_{\ell}}$ is represented by an open substack of $\SCoh(X)_{\cO_X, P_{\ell}}$.
\end{thm}

\begin{proof}
    To check that $\nps_{\alpha, \ell}$ is open in $\SCoh(X)_{\cO_X, P_{\ell}}$, since the latter is locally of finite type over a field, it is enough to check that for every family $f : B \to \SCoh(X)_{\cO_X, P_{\ell}}$ parameterized by a Noetherian scheme $B$, the locus $U\subset B$ of points that land in $\nps_{\alpha, \ell}$ is open.
    By \cref{prop: constructibility of the stack}, the subset $U$ is constructible.  By \cite[\href{https://stacks.math.columbia.edu/tag/004X}{004X}, \href{https://stacks.math.columbia.edu/tag/0542}{0542}]{stacks-project}, to conclude the proof of openness it suffices to show that $U \subset B$ is closed under generalization.

    Assume that we have a point $b\in B$ and some other point $b'$ in its closure such that $f(b') \in \nps_{\alpha, \ell}$.
    Then the inclusion $b' \in \overline{\{b\}}$ can be realized by a morphism $\Spec R \to B$ for some discrete valuation ring $R$.
    Let $\eta$ (resp. $s$) denote the generic (resp. special) point of $\Spec(R)$.
    The composition of $\Spec R \to B$ with $f$ produces a flat family of sheaves $F$ on $X_R$ with trivial determinant such that $F_s \in \nps_{\alpha, \ell}(s)$.
    We consider two cases depending on whether $F_{\eta}$ is torsion-free or not.
    
    If $F_\eta$ is torsion-free, then it is necessarily of rank $1$, because the rank of $F_s$ is $1$.
    Hence it admits an embedding into its reflexive hull $F_\eta \to \det(F_\eta) \cong \cO_{X_\eta}$ (see \cref{lemma: det(rank 1 tf) = double dual}), which means that it is an ideal sheaf.
    By flatness, we get that $P_{F_\eta} = P_{\ell}$, hence $F_\eta$ must be the ideal sheaf of a length $\ell$ subscheme of $X_\eta$.

    If $F_\eta$ has torsion, then $F_s$ must be of torsion type, because being torsion-free is preserved by generalization.
    This puts us in the situation of \cref{construction: limit in Quot}.
    By \cref{lemma: E eta is an ideal sheaf} and \cref{lemma: pure sheaves with support in beta}, $F_\eta$ is of torsion type in $\nps_{\alpha, \ell}(\eta)$.
\end{proof}

\subsection{Boundedness}
\begin{prop} \label{prop: boundedness}
    Let $\alpha \subset \Hilb(X)$ as in \cref{setup: alpha}. Suppose that $\alpha$ is of finite type over $k$. Then, the open substack $\nps_{\alpha, \ell} \subset \SCoh(X)_{\cO_X, P_{\ell}}$ is of finite type over $k$.
\end{prop}
\begin{proof}
It suffices to find a finite type stack $\mathscr{V}$ with a morphism $\mathscr{V} \to \SCoh(X)_{\cO_X, P_{\ell}}$ such that $\nps_{\alpha, \ell}$ is contained in the image (at the level of geometric points). Since the open substack $\SCoh^{\rtf}(X)_{\cO_X, P_{\ell}}$ is of finite type (for example, by \cref{prop: hilbert scheme as moduli space of sheaves}), it suffices to find a morphism $\mathscr{V} \to \SCoh(X)_{\cO_X, P_{\ell}}$ such that all sheaves of torsion type in $\nps_{\alpha, \ell}$ are contained in the image.

For each of the finitely many connected components $\alpha_j \subset \alpha$, consider the universal subscheme $Y_j \hookrightarrow X_{\alpha_j}$ over $\alpha_j$. Note that the fibers of $Y_j \to \alpha_j$ are geometrically integral, and they all have the same Hilbert polynomial, which we denote by $P_j$. Let $\overline{\SPic}(Y_j/\alpha_j)_{P_j -\ell}$ be the relative stack of families of rank one torsion-free sheaves on $Y_j \to \alpha_j$ with Hilbert polynomial $P_j - \ell$, where we use the polarization on $Y_j \to \alpha_j$ induced from the inclusion $Y_j \hookrightarrow X_{\alpha_j}$ to define the Hilbert polynomial. This stack is of finite type (for example, these are automatically Gieseker semistable, and so one may use the main result in \cite{langer-boundedness-general}). 

Set $H:= \bigsqcup_j \overline{\SPic}(Y_j/\alpha_j)_{P_j-\ell} \times \alpha_j$. We denote by $W \subset X_H$ the pullback of the universal subscheme $\bigsqcup_{j} Y_j \to \bigsqcup_j \alpha_j$ under the second projections 
\[H = \bigsqcup_j \overline{\SPic}(Y_j/\alpha_j)_{P_j-\ell} \times \alpha_j \xrightarrow{pr_2} \bigsqcup_j \alpha_j.\]
Similarly, we denote by $i: Y \hookrightarrow X_H$ the pullback of the universal $\bigsqcup_{j} Y_j \to \bigsqcup_j \alpha_j$ under the first projections 
\[
    H = \bigsqcup_j \overline{\SPic}(Y_j/\alpha_j)_{P_j-\ell} \times \alpha_j \xrightarrow{pr_1} \bigsqcup_j \overline{\SPic}(Y_j/\alpha_j)_{P_j-\ell} \to \bigsqcup_j \alpha_j.
\]
By definition we have an $H$-flat family of rank one torsion-free sheaves $G$ on $Y$ coming from the universal sheaf on $\bigsqcup_j \overline{\SPic}(Y_j/\alpha_j)_{P_j-\ell}$.
Let $\Extstack \to H$ denote the extension stack that sends a scheme $T \to H$ to the groupoid of extensions of sheaves on $X_T$
\[ 0 \to i_*(G)|_{X_T} \to F \to \cI_W|_{X_T} \to 0.\]
Note that $\Extstack$ is represented by an algebraic stack of finite type over $H$. Indeed, this can be checked smooth locally on $H$. Using \cite[\href{https://stacks.math.columbia.edu/tag/0A1J}{0A1J}]{stacks-project} for the $H$-flat sheaves $\cI_W$ and $i_*(G)$ on $X_H$, we see that after passing to a smooth cover $H' \to H$ by an affine scheme, we may find a finite complex of vector bundles $V^{\bullet}$ that universally computes $\Ext^{\bullet}( \cI_W, i_*(G))$. If we denote by $\Tot(V^i) \to H'$ the total space of the vector bundles appearing in the complex, then the differential induces morphisms of relatively affine $H'$-schemes $\Tot(V^0) \xrightarrow{\partial} \Tot(V^1) \xrightarrow{\partial} \Tot(V^2)$. Let $Z \subset \Tot(V^1)$ denote the preimage of the $0$-section $0: H' \to \Tot(V^2)$ under the differential $\Tot(V^1) \xrightarrow{\partial} \Tot(V^2)$. We may view $\Tot(V^0)$ as a smooth $H'$-group scheme under addition, and then $\Tot(V^0)$ acts on $Z$ at the level of scheme-valued points by $t \cdot z := \partial(t) +z$. The quotient stack $Z/\Tot(V^0)$ is evidently of finite type over $H'$, and it follows from construction that it represents the base-change $\Extstack_{H'}$.

Since $H$ is of finite type, it follows that $\Extstack$ is of finite type. The universal extension $F_{\univ}$ on $X_{\Extstack}$ induces a morphism $\Extstack \to \SCoh(X)$, and it follows from construction that the image contains all sheaves of torsion type in $\nps_{\alpha, \ell}$. Then, setting $\mathscr{V} \subset \Extstack$ to be the preimage of the closed substack $\SCoh(X)_{\cO_X, P_{\ell}} \subset \SCoh(X)$ we obtain the desired morphism $\mathscr{V} \to \SCoh(X)_{\cO_X, P_{\ell}}$.
\end{proof}

% ---
\section{Existence of the moduli space} \label{section: existence of moduli space}

This section is the technical heart of the paper, where we verify the local properties of $\cU_{\alpha,\ell}$ from \cite[Thm. 5.4]{AHLH} that guarantee the existence of a good moduli space.
We start by classifying the closed points of the stack and describing their groups of automorphisms (\cref{prop: closed points}, \cref{prop: local reductivity}).
We then prove that $\cU_{\alpha,\ell}$ is $\Theta$-reductive (\cref{prop: theta reductivity}).
To show S-completeness, we invoke a result of Zhang (\cite[Thm. 4.15]{zhang2024axis_completeness}, see \cref{lemma: axis implies S}), which allows us to verify the more intuitive property of axis-completeness instead (\cref{prop: S completeness}).
We continue by showing the valuative criterion for properness (\cref{prop: existence valuative criterion}).

The section culminates with the existence of a good moduli space for $\cU_{\alpha,\ell}$, denoted $\Hilb^\ell(X)_\alpha$ (\cref{thm: existence of proper}), and its main properties:
we determine the locus $\Hilb^\ell(X)^{\rnc\alpha} \subset \Hilb^\ell(X)$ that is unchanged under the morphism $\Hilb^\ell(X) \to \Hilb^\ell(X)_\alpha$ (\cref{notation: nc-alpha}, \cref{prop: locus where the modified Hilbert scheme is an iso}),
and show a condition that ensures that $\Hilb^\ell(X) \to \Hilb^\ell(X)_\alpha$ is surjective (\cref{prop: condition for surjectivity from Hilbert scheme}).

% -
\subsection{Closed points and local reductivity}

\begin{prop}
\label{prop: closed points}
    Let $\alpha \subset \Hilb(X)$ as in \cref{setup: alpha} and $\ell>0$. Suppose that $(\alpha, \ell)$ satisfies condition (LR) as in \cref{defn: condition dagger}. For any algebraic field extension $K \supset k$, a sheaf $F \in \nps_{\alpha, \ell}(K)$ is a closed point of $\nps_{\alpha, \ell}$ if and only if one of the following is satisfied:
    \begin{enumerate}
        \item $F \cong \cI_Z$ is of ideal type and for all geometric points $[Y \subset X_{\overline{k}}] \in \alpha(\overline{k})$ we have $Z_{\overline{k}} \not\subset Y$ scheme-theoretically.
        \item $F$ is of split torsion type.
    \end{enumerate}
\end{prop}
\begin{proof}
    After base-changing to the algebraic closure of $k$, we may assume without loss of generality that $k = K =\overline{k}$. We need to show that a $k$-point of $\nps_{\alpha, \ell}$ is closed if and only if it is of the form (1) or (2).
    Suppose first that $F \in \nps_{\alpha, \ell}(k)$ is not of the form (1) or (2). We show that $F$ is not closed by considering two cases:
    \begin{itemize}
        \item Assume that $F \cong \cI_Z$ for some length-$\ell$ subscheme $Z \subset X$ which is contained in some $[i: Y \hookrightarrow X] \in \alpha(k)$. Then there is a short exact sequence of sheaves
    \[ 0 \to \cI_Y \to F \to i_*(G) \to 0, \]
    where $G$ denotes the ideal sheaf of $Z\subset Y$. We may view $\cI_Y \subset F$ as a weighted filtration of $F$, and then the Rees construction \cite[Prop. 1.0.1]{halpernleistner2018structure} yields a morphism $\psi: [\mathbb{A}^1_{k}/\mathbb{G}_m] \to \nps_{\alpha, \ell}$ such that $\psi(1) \cong F$ and $\psi(0) \cong i_*(G) \oplus \cI_Y$. This shows that $F$ is not a closed point.

    \item Suppose that $F$ is of nonsplit torsion type. This means that $F$ fits into a short exact sequence
    \[ 0 \to i_*(G) \to F \to \cI_W \to 0\]
    which does not split. By considering $i_*(G) \subset F$ as a weighted filtration of $F$, the Rees construction yields a degeneration of $F$ to the sheaf $i_*(G) \oplus \cI_W$, thus showing that $F$ is not closed.
    \end{itemize}

    We are left to show that $F$ is closed whenever it is of the form (1) or (2).
    Suppose first that $F$ is of the form $\cI_Z$ as in (1).
    Assume for the sake of contradiction $F$ is not a closed point.
    It follows that there is a discrete valuation ring $R$ with generic point $\eta$ and special point $s \cong \Spec(k)$, and a morphism $\Spec(R) \to \nps_{\alpha, \ell}$ corresponding to a sheaf $\widetilde{F}$ on $X_R$ such that the generic fiber $\widetilde{F}_{\eta}$ is isomorphic to the base-change $F_{\eta}$ and such that $\widetilde{F}_s \not\cong F$.
    Note that the special fiber $\widetilde{F}_s$ of such a degeneration must be of torsion type, since every $k$-point of the open substack $\SCoh^{\rtf}(X)_{\cO_X, P_{\ell}} \subset \nps_{\alpha, \ell}$ is closed by \cref{prop: hilbert scheme as moduli space of sheaves}.
    By a Rees degeneration argument as above, we may further specialize to a sheaf of split torsion type and assume without loss of generality that $\widetilde{F}_s \cong i_*(G) \oplus \cI_W$. By upper semicontinuity of the dimension of $\Hom(-, \cI_Z)$ in families \cite[Cor. 7.7.8.2]{egaiii_2}, we must have that $\Hom(i_*(G) \oplus \cI_W, \cI_Z)$ is nonzero. Since $\cI_Z$ is torsion-free, we have $\Hom(i_*(G), \cI_Z) =0$, and it follows that $\Hom(\cI_W, \cI_Z) \neq 0$. If we choose a nonzero morphism $f: \cI_W \to \cI_Z$, then $f$ is a monomorphism and the composition $\cI_W \hookrightarrow \cI_Z \hookrightarrow \cO_X$ corresponds to an element $[W' \subset X]$ of $\alpha(k)$ with $\cI_{W'} \cong \cI_{W}$ and $Z \subset W'$. This contradicts the assumption that $Z$ is not contained in a subscheme in $\alpha$.

    We are left to consider the case when $F$ is of split torsion type, so we may write it in the form $F = i_*(G) \oplus \cI_W$.
    Assume for the sake of contradiction that $F$ is not a closed point of $\nps_{\alpha, \ell}$. As in the previous paragraph, there exists a family $\widetilde{F}: \Spec(R) \to \nps_{\alpha, \ell}$ such that $\widetilde{F}_{\eta} \cong F_{\eta}$ and such that $\widetilde{F}_s \cong i'_*(G') \oplus \cI_{W'}$ is of split torsion type with $\widetilde{F}_s \not\cong F$.
    By upper semicontinuity of the dimension of $\Hom(-, \cI_W)$, we must have that $\Hom(i'_*(G') \oplus \cI_{W'}, \cI_W) \neq 0$.
    Since $\cI_W$ is torsion-free, it follows that $\Hom(i'_*(G'), \cI_W) =0$, and hence we must have a nonzero morphism $\cI_{W'} \to \cI_W$. By condition (LR), this forces $\cI_W \cong \cI_{W'}$, and so $[i: Y \hookrightarrow X]$ and $[i': Y' \hookrightarrow X]$ are integral subschemes in the same component of $\Hilb(X)$, and $P_{i_*(G)} = P_{(i')_*(G')}$.
    By upper-semicontinuity of $\Hom(i_*(G), -)$, we must have $\Hom(i_*(G), i'_*(G') \oplus \cI_{W'}) \neq 0$.
    Since $\cI_{W'}$ is torsion-free, we have $\Hom(i_*(G), \cI_{W'}) =0$, and hence it follows that $\Hom(i_*(G), i'_*(G')) \neq 0$. Using adjunction and the fact that $Y,Y' \subset X$ are integral of the same dimension, this implies that $i=i'$ and $G \subset G'$.
    Furthermore, since $P_{i_*(G)}= P_{(i')_*(G')}$, this inclusion must be an equality $G = G'$.
    We have thus shown that $\cI_W \cong \cI_{W'}$ and $i_*(G) \cong i'_*(G')$, which proves that $F \cong \widetilde{F}_s$, a contradiction.
\end{proof}

\begin{prop}
\label{prop: local reductivity}
    Suppose that $(\alpha, \ell)$ satisfies condition (LR) as in \cref{defn: condition dagger}. Then, the stack $\nps_{\alpha, \ell}$ is locally reductive in the sense of \cite[Def 2.5]{AHLH}.
\end{prop}

\begin{proof}
    The closed points of $\nps_{\alpha, \ell}$ are of the form (1) or (2) as described in \cref{prop: closed points}. Note that a point as in \cref{prop: closed points}(1) has stabilizer isomorphic to $\mathbb{G}_{m}$, whereas the assumption (LR) implies that a point as in \cref{prop: closed points}(2) has stabilizer isomorphic to $\mathbb{G}_{m}^2$. In particular, we see that the stabilizers of all closed points of $\nps_{\alpha, \ell}$ are linearly reductive. By the Luna \'etale slice theorem for stacks (see e.g. \cite[Thm. 2.2]{AHLH}), in order to show that $\nps_{\alpha, \ell}$ is locally reductive it is enough to show that every geometric point $x \in \nps_{\alpha, \ell}(K)$ specializes to a closed point. Since $\nps_{\alpha, \ell}$ is locally of finite type over $k$, any such geometric point $x$ specializes to a point defined over the algebraic closure $\overline{k}$, and so we may assume that $x \in \nps_{\alpha, \ell}(\overline{k})$.
    
    Suppose first that $x$ corresponds to a sheaf $F \cong \cI_Z$ of ideal type for some length $\ell$ subscheme $Z \subset X_{\overline{k}}$. We distinguish two subcases:
    \begin{itemize}
        \item If $Z$ is contained in some $[i: Y \hookrightarrow X_{\overline{k}}]$, then the Rees degeneration as in the proof of \cref{prop: closed points} exhibits a specialization of $x$ to a sheaf of split torsion type $i_*(G) \oplus \cI_Y$, which is a closed point of $\nps_{\alpha, \ell}$ by \cref{prop: closed points}. 

    \item If $Z$ is not contained in any $[Y \subset X_{\overline{k}}] \in \alpha(\overline{k})$, then it follows that $x$ is already a closed point of $\nps_{\alpha, \ell}$ by \cref{prop: closed points}.
    \end{itemize}

    We are left to deal with the case when $x$ corresponds to a sheaf $F$ of torsion type, which is of the form
    \[ 0 \to i_*(G) \to F \to \cI_W \to 0.\]
     We may view $i_*(G) \subset F$ as a weighted filtration of $F$, and then the Rees construction yields a degeneration $\psi: [\mathbb{A}^1_{\overline{k}}/\mathbb{G}_m] \to \nps_{\alpha, \ell}$ such that $\psi(1) \cong F$ and $\psi(0) \cong i_*(G) \oplus \cI_W$. Since $i_*(G) \oplus \cI_W$ is of split torsion type, we conclude by \cref{prop: closed points} that $x$ specializes to a closed point of $\nps_{\alpha, \ell}$.
\end{proof}

% -
\subsection{\texorpdfstring{$\Theta$}{Theta}-reductivity}

We begin by recalling a couple of standard lemmas that will be used in the proof of $\Theta$-reductivity.

\begin{lemma} \label{lemma: EGA lemma}
    Let $f: Y \to T$ be a proper morphism of Noetherian schemes. Let $F,I$ be coherent sheaves on $Y$ such that $F$ is $T$-flat and $I$ is isomorphic to the cokernel of a morphism $\varphi: E \to H$ of locally finitely generated locally free sheaves on $Y$. Then there exists a coherent sheaf $N$ on $T$ such that for all quasicoherent sheaves $M$ on $T$ we have an isomorphism of quasicoherent sheaves
    \[ f_*(\cH om(I, F \otimes f^*(M))) \xrightarrow{\sim} \cH om(N, M)\]
    which is natural in the variable $M$.
\end{lemma}
\begin{proof}
    This is the content of \cite[Cor. 7.7.8]{egaiii_2}.
\end{proof}

\begin{lemma}
\label{lemma: global hom from fiberwise homs}
    Let $R$ be a discrete valuation ring with fraction field $K$ and residue field $\kappa$, and let $\eta$ and $s$ be the corresponding open and closed points of $\Spec(R)$, respectively.
    Let $I$ and $F$ be $R$-flat coherent sheaves on $X_R$.
    Assume that $\dim \Hom_{X_\eta}(I_\eta,F_\eta) = \dim \Hom_{X_s}(I_s,F_s) = n$ for some integer $n$.
    Then $\Hom_{X_R}(I,F)$ is a free $R$-module of rank $n$ whose restriction to $\eta$ (resp. $s$) is isomorphic to $\Hom_{X_\eta}(I_\eta,F_\eta)$ (resp. $\Hom_{X_s}(I_s,F_s)$); in particular, there exists a homomorphism $I \to F$ that is $R$-fiberwise nonzero.
\end{lemma}

\begin{proof}
    By \Cref{lemma: EGA lemma} with $T=\Spec(R)$, there is a finitely generated $R$-module $N$ such that for all $R$-modules $M$, there is a natural isomorphism
    $\Hom_{X_R}(I,F\otimes_R M) \cong \Hom_R (N,M)$.
    Setting $M=K$ (resp. $M=\kappa$) calculates that the dimension of $\Hom_R(N,K) \cong (N_K)^\vee$ (resp. $(N_\kappa)^\vee$) is $n$. Hence $N$ is a finitely generated module over $R$ of constant rank $n$, and so we have $N \cong R^{\oplus n}$. 
    Thus $\Hom_{X_R}(I,F) \cong \Hom_R(N,R)\cong R^{\oplus n}$. Any basis element in $\Hom_{X_R}(I,F) \cong R^{\oplus n}$ yields a homomorphism $I \to F$ which is $R$-fiberwise nonzero.
\end{proof}

\begin{lemma}
\label{lemma: Hom from I W' to non-split F is zero}
    Let $\alpha$ be as in \cref{setup: alpha}.
    Suppose that $(\alpha, \ell)$ satisfies condition (LR) as in \cref{defn: condition dagger}.
    Let $k \subset K$ be a field extension and let $F \in \nps_{\alpha,\ell}(K)$ be a sheaf on $X_K$ of torsion type that fits into the short exact sequence
    \[
    0 \to i_*(G) \to F \to \cI_W \to 0
    \]
    as in \cref{defn: stack U_alpha}, with $W\in \alpha(K)$.
    Let $\cI_{W'}$ be an ideal sheaf with $W' \in \alpha(K)$.
    If $\Hom_{X_K} (\cI_{W'},F) \neq 0$, then the following hold:
    \begin{enumerate}
        \item 
        the short exact sequence splits (i.e., $F$ is of split torsion type);
        \item 
        $\cI_W \cong \cI_{W'}$;
        \item 
        $\dim \Hom(\cI_W,F) = 1$, and if $f : \cI_{W'} \to F$ is nonzero, then
        $f$ factors through the direct summand $\cI_W$.
    \end{enumerate}
\end{lemma}

\begin{proof}
    Let $f : \cI_{W'} \to F$ be a nonzero morphism.
    If $\cI_W \not\cong \cI_{W'}$, then by condition (LR), we have $\Hom(\cI_{W'},\cI_W)=0$, which implies that $f$ factors through $i_*(G)$, but
    this contradicts condition (LR).
    We conclude that $\cI_W \cong \cI_{W'}$.
    Noting $\Hom(\cI_W,\cI_W) \cong K$, we conclude that $f$ gives a splitting of the short exact sequence.
    By condition (LR), $\Hom(\cI_W,i_*(G))=0$, which yields the claims in part (3).
\end{proof}

\begin{lemma}
\label{lemma: Hom from I W' to ideal F is at most dimension 1}
    Let $\alpha$ be as in \cref{setup: alpha}.
    Suppose that $(\alpha, \ell)$ satisfies condition (LR) as in \cref{defn: condition dagger}.
    Let $k \subset K$ be a field extension and let $F \in \nps_{\alpha,\ell}(K)$ be a sheaf of ideal type. For all $[W \subset X_K] \in \alpha(K)$, we have $\dim \Hom(\cI_W, F)  \leq 1$.
\end{lemma}
\begin{proof}
    Suppose for the sake of contradiction that there are two $K$-linearly independent homomorphisms $\psi_1, \psi_2: \cI_W \to F$. The cokernel of $\psi_1$ is a sheaf of the form $i_*(G)$, where $[i: Y \hookrightarrow X_K] \in \alpha(K)$ in the same component of $\alpha$ as $W$, and $G$ is a rank one torsion-free sheaf with Hilbert polynomial $P_{\cO_W}-\ell = P_{\cO_Y}-\ell$ on $Y$. Then, the composition $\cI_W \xrightarrow{\psi_2} F \twoheadrightarrow i_*(G)$ is nonzero, thus contradicting condition (LR).
\end{proof}

For the next definition and proposition, we use the standard notation $\Theta_R := [\mathbb{A}^1_R/\mathbb{G}_m]$ for any ring $R$, where $\mathbb{G}_m$ acts linearly on $\mathbb{A}^1_R = \Spec(R[t])$ so that the coordinate $t$ has weight $-1$. If $R$ is a discrete valuation ring, we denote by $\mathfrak{o} \subset \Theta_R$ 
the residual gerbe of the unique closed point of the stack $\Theta_R$ (corresponding to the origin of the special fiber of $\mathbb{A}^1_R \to \Spec(R)$). 
\begin{defn}[{\cite[Def. 3.10]{AHLH}}]
\label{defn: Theta reductivity}
    Let $\mathscr{X} \to \sY$ be a morphism of algebraic stacks locally of finite type over $k$.
    We say that $\mathscr{X} \to \sY$ is \emph{$\Theta$-reductive} if it satisfies the unique right lifting property with respect to $\Theta_R \setminus \mathfrak{o} \to \Theta_R$, for any discrete valuation ring $R$.
    In other words, any commutative diagram of solid arrows as below admits a filling with the dotted arrow, and the filling is unique up to a unique isomorphism.
    % https://q.uiver.app/#q=WzAsNCxbMCwwLCJcXGNYXFxzZXRtaW51cyAwIl0sWzIsMCwiXFxzWSJdLFsyLDIsIlxcc1oiXSxbMCwyLCJcXGNYIl0sWzAsMV0sWzMsMl0sWzAsM10sWzEsMl0sWzMsMSwiIiwxLHsic3R5bGUiOnsiYm9keSI6eyJuYW1lIjoiZGFzaGVkIn19fV1d
    \[\begin{tikzcd}
	{\Theta_R\setminus \mathfrak{o}} && \sX \\
	\\
	\Theta_R && \sY
	\arrow[from=1-1, to=1-3]
	\arrow[from=1-1, to=3-1]
	\arrow[from=1-3, to=3-3]
	\arrow[dashed, from=3-1, to=1-3]
	\arrow[from=3-1, to=3-3]
    \end{tikzcd}\]
    We say that a stack $\sX$ is $\Theta$-reductive if the structure morphism $\sX \to \Spec(k)$ is $\Theta$-reductive.
\end{defn}

\begin{prop} \label{prop: theta reductivity}
    Let $\alpha$ be as in \cref{setup: alpha}.
    Suppose that $(\alpha, \ell)$ satisfies condition (LR) as in \cref{defn: condition dagger}.
    Then, the stack $\nps_{\alpha, \ell}$ is $\Theta$-reductive.
\end{prop}

\begin{proof}
    Let $\Theta_R \setminus \mathfrak{o} \to \nps_{\alpha,\ell} \subset \SCoh(X)$ be a morphism.
    Since $\SCoh(X)$ is $\Theta$-reductive by \cite[Lem. 7.17, Prop. 3.18]{AHLH}, there is a unique extension $\Theta_R \to \SCoh(X)$.
    This translates, by \cite[Cor. 7.13]{AHLH}, into having an $R$ flat family $F$ of sheaves in $\nps_{\alpha,\ell}(R)$ endowed with a decreasing filtration $\cdots \supset \Fil_0(F) \supset \Fil_1(F) \supset \cdots$ such that the following hold:
    \begin{enumerate}
        \item \label{item 1 filtration} $\Fil_{-N}(F) = F$ and $\Fil_N(F) = 0$ for some $N\gg 0$,
        \item \label{item 2 filtration} the associated graded object $\gr_\bullet(F)$ is $R$-flat,
        \item \label{item 3 filtation} $\gr_\bullet(F_{\eta})$ is in $\nps_{\alpha,\ell}(\eta)$.
    \end{enumerate}
    It remains to check that $\gr_\bullet(F_s)$ is in $\nps_{\alpha,\ell}(s)$.

    If $\gr_\bullet(F_\eta)$ is concentrated in one degree, then the filtration of $F_\eta$ is shifted trivial.
    One can extend this filtration to a shifted trivial filtration of $F_s$; by uniqueness of the extension to $\SCoh(X)$, this coincides with the filtration given by $\Theta_R \to \SCoh(X)$.
    But the associated graded $\gr_\bullet(F_s)$ of a shifted trivial filtration is just $F_s$, which is in $\nps_{\alpha,\ell}$, as desired.

    If $\gr_\bullet(F_\eta)$ is concentrated in more than one degree, then it must be decomposable and hence of the split torsion type.
    In other words, we can express $\gr_\bullet(F_\eta) \cong i_{\eta,*}(G_\eta) \oplus \cI_{W_\eta}$ for some $[W_\eta\subset X_{\eta}],  [i_{\eta} : Y_\eta \to X_\eta] \in \alpha_j(\eta)$, and a rank $1$ torsion-free sheaf $G_\eta$ on $Y_\eta$.
    Both of the summands $i_{\eta,*}(G_\eta)$ and $\cI_{W_{\eta}}$ are indecomposable, since both are (pushforwards of) rank $1$ torsion-free sheaves on integral schemes. Hence the graded object $\gr_\bullet(F_\eta)$ is concentrated in two degrees, say $0$ and $1$ after appropriate shifts, and we have a short exact sequence
    \begin{equation}
    \label{eq: SES of one-step filtration}
        0 \to \gr_1(F_\eta) \to F_\eta \to \gr_0(F_\eta) \to 0.
    \end{equation}
    We will now consider two cases.

    \paragraph{\textbf{Case: $\Fil_1(F_\eta) = \gr_1(F_\eta) \cong \cI_{W_\eta}$.}}
    Take the closure $W := \overline{W_\eta}$ in $X_R$.
    Note that we have $W_s \in \alpha(s)$.
    By upper-semicontinuity of Hom, we have
    $\Hom(\cI_{W_s},F_s) \neq 0$.
    Then \cref{lemma: Hom from I W' to non-split F is zero} implies that $F_s$ is either of the ideal type or $F_s= \cI_{W_s} \oplus i_{s,*}(G_s)$ for some $i_s: Y \hookrightarrow X_s$ and rank one torsion-free sheaf $G_s$.
    In either case, \Cref{lemma: Hom from I W' to non-split F is zero} and \Cref{lemma: Hom from I W' to ideal F is at most dimension 1} imply that $\Hom(\cI_{W_x},F_x)$ is $1$-dimensional for $x = \eta$ and $x=s$.
    By \cref{lemma: global hom from fiberwise homs}, there exists a fiberwise nonzero morphism $\cI_W \to F$; in particular, it must be injective over $s$.
    It now follows from \cite[\href{https://stacks.math.columbia.edu/tag/00ME}{00ME}]{stacks-project} that the morphism $\cI_W \to F$ is injective and the quotient $F/\cI_W$ is $R$-flat.
    Therefore, by the uniqueness of limits in the Quot scheme, the injection $\cI_W \to F$ induces the given filtration $\Fil_\bullet(F)$.
    If $F_s$ is of the ideal type (resp. split torsion type), we get that $\gr_\bullet(F_s) \cong \cI_{W_s} \oplus (\cI_Z / \cI_{W_s})$ with $\cI_Z / \cI_{W_s}$ a pushforward of a rank $1$ torsion-free sheaf from $W_s \in \alpha(s)$ (resp. $\gr_\bullet(F_s) \cong \cI_{W_s} \oplus i_{s,*}(G_s)$),
    hence we conclude that $\gr_\bullet(F_s) \in \cU_{\alpha,\ell}(s)$.

    \paragraph{\textbf{Case: $\Fil_1(F_\eta) = \gr_1(F_\eta) \cong i_{\eta,*}(G_\eta)$.}}
    In this case, $F_\eta$ is a torsion sheaf, and therefore, by openness of torsion-freeness, $F_s$ must be of torsion type as well.
    Take the extension $0\to T\to F\to E\to 0$ of the torsion sequence $0\to i_{\eta,*}(G_\eta) \to F_\eta \to \cI_{W_\eta} \to 0$ as in \cref{construction: limit in Quot}.
    By uniqueness of the extension of a filtration, we get a natural identification $T = \Fil_1(F)$.
    By \cref{lemma: E eta is an ideal sheaf}, $0\to T_s \to F_s \to E_s\to 0$ is the torsion filtration of $F_s$.
    But then $\gr_\bullet(F_s) \cong T_s \oplus E_s$ is in $\cU_{\alpha, \ell}(s)$, as desired.
\end{proof}

% -
\subsection{\texorpdfstring{$S$}{S}-completeness}

In this subsection we show that $\nps_{\alpha, \ell}$ is $S$-complete in the sense of \cite[Def. 3.38]{AHLH}. In order to show S-completeness of $\cU_{\alpha,\ell}$, we check a different property, called axis-completeness, introduced in \cite{zhang2024axis_completeness}. 
We now recall the corresponding definitions.
We introduce the notation $\ST_R := \left[ \Spec(R[s,t]/(st-\pi)) / \bG_m \right]$, for a discrete valuation ring $R$ with uniformizer $\pi$, and $\cX := \left[ \Spec(\bZ[s,t]/st) / \bG_m \right]$.
In both cases, $\bG_m$ acts on $s$ with weight $1$ and on $t$ with weight $-1$.
We let $\mathfrak{o}$ denote the unique closed point of $\ST_R$, and
if $K$ is a field, we let $\mathfrak{o}$ denote the unique closed point of $\cX_K$.
\begin{defn}[{\cite[Def. 3.38]{AHLH}, \cite[Def. 4.14]{zhang2024axis_completeness}}] \label{defn: S and axis completeness}
A morphism $f : \sY \to \sZ$ of locally of finite type algebraic stacks over a field $k$ is \emph{S-complete} (resp. \emph{axis-complete}) if it satisfies the unique right lifting property with respect to 
$\ST_R \setminus \mathfrak{o} \to \ST_R$ for any discrete valuation ring $R$
(resp. $\cX_K\setminus \mathfrak{o} \to \cX_K$ for any field $K$), analogously to the filling condition in \cref{defn: Theta reductivity}. We say that a stack $\sY$ is $S$-complete if the structure morphism $\sY \to \Spec(k)$ is $S$-complete.
\end{defn}

\begin{lemma}[{cf. \cite[Thm. 4.15]{zhang2024axis_completeness}}]
\label{lemma: axis implies S}
    Let $\sY$ be an S-complete algebraic stack over a field $k$.
    Let $\sU \subset \sY$ be an open substack.
    If the inclusion morphism $\sU \to \sY$ is axis-complete (e.g. if $\sU$ is axis-complete), then $\sU$ is S-complete.
\end{lemma}

\begin{proof}
    The first part of the proof in \cite[Thm. 4.15]{zhang2024axis_completeness} works verbatim in this case.
\end{proof}

\begin{prop} \label{prop: S completeness}
    Let $\alpha$ be as in \cref{setup: alpha}. Suppose that $(\alpha, \ell)$ satisfies condition (LR) as in \cref{defn: condition dagger}.
    Then the stack $\nps_{\alpha, \ell}$ is $S$-complete.
\end{prop}

\begin{proof}
    By \cite[Lem. 7.16, Prop. 3.42]{AHLH}, the stack $\SCoh(X)$ is S-complete,
    hence $\SCoh(X)_{\cO_X}$, being its closed substack, is S-complete.
    By \cref{lemma: axis implies S}, it is now enough to prove that $\cU_{\alpha,\ell} \to \SCoh(X)_{\cO_X}$ is axis-complete.
    By \cite[Lem. 4.16]{zhang2024axis_completeness}, we should check the following: given two geometric points $F, F' \in \cU_{\alpha,\ell}(K)$ for some field $K$, given two finite opposite filtrations $\Fil_\bullet(F)$ and $\Fil^\bullet(F')$ such that $\gr(i) := \gr_i(F) \cong \gr^i(F')$ for all $i \in \bZ$, we have that the (isomorphic) associated graded objects $\gr:=\gr_\bullet(F)\cong\gr^\bullet(F')$ are in $\cU_{\alpha,\ell}(K)$.
    See the following diagram for illustration, where dotted lines indicate the steps of the filtration that have isomorphic associated graded objects.
    % https://q.uiver.app/#q=WzAsOCxbMCwwLCJGID0gXFxGaWxfezB9KEYpIl0sWzIsMCwiXFxGaWxfezF9KEYpIl0sWzIsMiwiXFxGaWxeezB9KEYnKSJdLFswLDIsIjAgPSBcXEZpbF57LTF9KEYnKSJdLFs0LDAsIlxcRmlsXzIoRikiXSxbNiwwLCJcXGNkb3RzIl0sWzQsMiwiXFxGaWxeMShGJykiXSxbNiwyLCJcXGNkb3RzIl0sWzMsMiwiIiwwLHsic3R5bGUiOnsidGFpbCI6eyJuYW1lIjoiaG9vayIsInNpZGUiOiJ0b3AifX19XSxbMSwwLCIiLDEseyJzdHlsZSI6eyJ0YWlsIjp7Im5hbWUiOiJob29rIiwic2lkZSI6ImJvdHRvbSJ9fX1dLFs0LDEsIiIsMSx7InN0eWxlIjp7InRhaWwiOnsibmFtZSI6Imhvb2siLCJzaWRlIjoiYm90dG9tIn19fV0sWzIsNiwiIiwwLHsic3R5bGUiOnsidGFpbCI6eyJuYW1lIjoiaG9vayIsInNpZGUiOiJ0b3AifX19XSxbNSw0LCIiLDEseyJzdHlsZSI6eyJ0YWlsIjp7Im5hbWUiOiJob29rIiwic2lkZSI6ImJvdHRvbSJ9fX1dLFs2LDcsIiIsMCx7InN0eWxlIjp7InRhaWwiOnsibmFtZSI6Imhvb2siLCJzaWRlIjoidG9wIn19fV0sWzksOCwiXFxncigwKSIsMSx7InNob3J0ZW4iOnsic291cmNlIjoyMCwidGFyZ2V0IjoyMH0sImxldmVsIjoxLCJzdHlsZSI6eyJib2R5Ijp7Im5hbWUiOiJkb3R0ZWQifSwiaGVhZCI6eyJuYW1lIjoibm9uZSJ9fX1dLFsxMCwxMSwiXFxncigxKSIsMSx7InNob3J0ZW4iOnsic291cmNlIjoyMCwidGFyZ2V0IjoyMH0sImxldmVsIjoxLCJzdHlsZSI6eyJib2R5Ijp7Im5hbWUiOiJkb3R0ZWQifSwiaGVhZCI6eyJuYW1lIjoibm9uZSJ9fX1dXQ==
    \begin{equation}
    \label{eq: diagram for axis-completeness}
    \begin{tikzcd}[ampersand replacement = \&]
    	{F = \Fil_{0}(F)} \&\& {\Fil_{1}(F)} \&\& {\Fil_2(F)} \&\& \cdots \\
    	\\
    	{0 = \Fil^{-1}(F')} \&\& {\Fil^{0}(F')} \&\& {\Fil^1(F')} \&\& \cdots
    	\arrow[""{name=0, anchor=center, inner sep=0}, hook', from=1-3, to=1-1]
    	\arrow[""{name=1, anchor=center, inner sep=0}, hook', from=1-5, to=1-3]
    	\arrow[""{name=4, anchor=center, inner sep=0}, hook', from=1-7, to=1-5]
    	\arrow[""{name=2, anchor=center, inner sep=0}, hook, from=3-1, to=3-3]
    	\arrow[""{name=3, anchor=center, inner sep=0}, hook, from=3-3, to=3-5]
    	\arrow[""{name=5, anchor=center, inner sep=0}, hook, from=3-5, to=3-7]
    	\arrow["{\gr(0)}"{description}, between={0.2}{0.8}, dotted, no head, from=0, to=2]
    	\arrow["{\gr(1)}"{description}, between={0.2}{0.8}, dotted, no head, from=1, to=3]
    	\arrow["{\gr(2)}"{description}, between={0.2}{0.8}, dotted, no head, from=4, to=5]
    \end{tikzcd}
    \end{equation}

    If the filtration $\Fil_\bullet(F)$ is shifted trivial, then so is $\Fil^\bullet(F')$, and we have $\gr\cong F\cong F' \in \cU_{\alpha,\ell}(K)$.
    We are thus left with the case when the filtration $\Fil_\bullet(F)$ is not shifted trivial, and we will treat separately different types of $F$.
    After appropriate shifts, we can assume that $\gr_\bullet(F)$ is concentrated in degrees $0,\dots,m$ for some $m>0$.

    \paragraph{\textbf{Case: $F'$ is of the ideal type.}}
    We can then write $F' = \cI_Z$ for some length $\ell$ subscheme of $X_K$.
    By \eqref{eq: diagram for axis-completeness}, we have an inclusion $\gr(0) = F / \Fil_1(F) \to F' = \cI_Z$, hence $\gr(0)$ must be torsion-free.
    Since $F$ is a rank $1$ sheaf and its quotient $F / \Fil_1(F)$ is torsion-free, $F$ is of torsion type, and its torsion sequence as in \cref{defn: stack U_alpha} is forced to be $\Tors(F) = \Fil_1(F) = i_*(G)$, with $[i: Y\to X_K] \in \alpha_j(K)$ and a rank $1$ torsion-free $G$ on $Y$, and $\TF(F) = \gr(0) = \cI_W$, for some $W\in \alpha_j(K)$.
    The diagram \eqref{eq: diagram for axis-completeness} now becomes:
    \begin{equation}
    \label{eq: diagram for axis-completeness with gr(0) = I W}
    \begin{tikzcd}[ampersand replacement = \&]
    	{F} \&\& {i_*(G)} \&\& {\Fil_2(F)} \&\& \cdots \\
    	\\
    	{0} \&\& {\cI_W} \&\& {\Fil^1(F')} \&\& \cdots
    	\arrow[""{name=0, anchor=center, inner sep=0}, hook', from=1-3, to=1-1]
    	\arrow[""{name=1, anchor=center, inner sep=0}, hook', from=1-5, to=1-3]
    	\arrow[""{name=4, anchor=center, inner sep=0}, hook', from=1-7, to=1-5]
    	\arrow[""{name=2, anchor=center, inner sep=0}, hook, from=3-1, to=3-3]
    	\arrow[""{name=3, anchor=center, inner sep=0}, hook, from=3-3, to=3-5]
    	\arrow[""{name=5, anchor=center, inner sep=0}, hook, from=3-5, to=3-7]
    	\arrow["{\cI_W}"{description}, between={0.2}{0.8}, dotted, no head, from=0, to=2]
    	\arrow["{\gr(1)}"{description}, between={0.2}{0.8}, dotted, no head, from=1, to=3]
    	\arrow["{\gr(2)}"{description}, between={0.2}{0.8}, dotted, no head, from=4, to=5]
    \end{tikzcd}
    \end{equation}
    We notice that $\Fil^1(F') / \cI_W \subset \cI_Z / \cI_W \subset \cO_{X_K} / \cI_W \cong \cO_W$.
    The next step in the filtration then yields $i_*(G) / \Fil_2(F) \cong \Fil^1(F') / \cI_W \subset \cO_W$, hence $\gr(1) = i_*(G) / \Fil_2(F)$ is the pushforward of a torsion-free sheaf from the integral subscheme $W \subset X_K$.
    However, $\gr(1)$ is also supported on the integral subscheme $Y \subset X_K$, which has the same dimension as $W$, so we must have $[W\subset X_K]=[Y \subset X_K] \in \alpha_j(K)$.
    Therefore, $\gr(1) \subset i_* \cO_Y$, and hence $\gr(1)$ is the pushforward from $Y$ of a rank $1$ torsion-free sheaf.
    But recall that $G$ is a rank $1$ torsion-free sheaf on $Y$ as well, so $\gr(1) = i_*(G)/\Fil_2(F)$ is forced to be $i_*(G)$.
    We conclude that $\Fil_2(F)=0$, and hence $\gr = \gr(0)\oplus\gr(1) \cong \cI_Y \oplus i_*(G)$ is a sheaf in $\cU_{\alpha,\ell}(K)$, as desired.

    \paragraph{\textbf{Case: both $F$ and $F'$ are of torsion type.}}
    The isomorphism $F / \Fil_1(F) \cong \Fil_0(F')$ from \eqref{eq: diagram for axis-completeness} induces a nonzero morphism $f : F \to F'$, which is not an isomorphism by non-triviality of filtrations.
    We get an induced morphism between torsion sequences of $F$ and $F'$:
    % https://q.uiver.app/#q=WzAsMTAsWzAsMCwiMCJdLFsyLDAsImlfKihHKSJdLFs0LDAsIkYiXSxbNiwwLCJcXGNJX1ciXSxbOCwwLCIwIl0sWzAsMiwiMCJdLFsyLDIsImknXyooRycpIl0sWzQsMiwiRiciXSxbNiwyLCJcXGNJX3tXJ30iXSxbOCwyLCIwIl0sWzAsMV0sWzEsMl0sWzIsM10sWzMsNF0sWzUsNl0sWzYsN10sWzcsOF0sWzgsOV0sWzEsNiwiZyIsMV0sWzIsNywiZiIsMV0sWzMsOCwiaCIsMV1d
    \begin{equation}
    \label{eq: diagram of torsion sequences for axis-completeness}
    \begin{tikzcd}[ampersand replacement=\&]
    	0 \&\& {i_*(G)} \&\& F \&\& {\cI_W} \&\& 0 \\
    	\\
    	0 \&\& {i'_*(G')} \&\& {F'} \&\& {\cI_{W'}} \&\& 0
    	\arrow[from=1-1, to=1-3]
    	\arrow[from=1-3, to=1-5]
    	\arrow["g"{description}, from=1-3, to=3-3]
    	\arrow[from=1-5, to=1-7]
    	\arrow["f"{description}, from=1-5, to=3-5]
    	\arrow[from=1-7, to=1-9]
    	\arrow["h"{description}, from=1-7, to=3-7]
    	\arrow[from=3-1, to=3-3]
    	\arrow[from=3-3, to=3-5]
    	\arrow[from=3-5, to=3-7]
    	\arrow[from=3-7, to=3-9]
    \end{tikzcd}
    \end{equation}
    We let $Y \in \alpha_j(K)$ (resp. $Y' \in \alpha_{j'}(K)$) denote the support of $G$ (resp. $G'$), then $W \in \alpha_j(K)$ (resp. $W' \in \alpha_{j'}(K)$).

    \begin{itemize}
        \item
        Subcase: $g$ is zero.
        Since $f$ is not zero and $\Hom(\cI_W, i'_*(G'))=0$ by condition (LR), we must have that $h$ is nonzero.
        Hence $h$ is an isomorphism by condition (LR).
        Since $g=0$, the morphism $F / i_*(G) \cong \cI_W \xrightarrow{h} \cI_{W'}$ factors through $\cI_W \to F'$.
        By \cref{lemma: Hom from I W' to non-split F is zero}, $F'$ is of split torsion type, i.e. $F' \cong i'_*(G') \oplus \cI_{W'}$.
        We notice that both compositions $F \to \gr(0) \to F'$ and $F \to \cI_W \cong \cI_{W'} \subset F'$ are equal to $f$, so their images are the same, and hence $\gr(0) \cong \Fil^0(F') = \cI_W$.
        Therefore, the diagram \eqref{eq: diagram for axis-completeness} looks again like \eqref{eq: diagram for axis-completeness with gr(0) = I W}, and the same argument as in that case yields $\gr(2) = 0$ and $\gr \cong \cI_W \oplus i_*(G) \in \cU_{\alpha,\ell}(K)$.
        \item 
        Subcase: $g$ is nonzero.
        In this case $g:i_*(G) \to i'_*(G')$ would correspond by adjunction to a nonzero morphism $(i')^*i_*(G) \to G'$. Note that $(i')^*i_*(G)$ is supported on $Y \cap Y'$. Since $G'$ is torsion-free on the integral scheme $Y'$, this forces $Y\cap Y' = Y'$, or equivalently $Y \supset Y'$. The containment $Y \supset Y'$ corresponds to an inclusion $\cI_{Y} \hookrightarrow \cI_{Y'}$. By condition (LR), this inclusion of ideal sheaves is forced to be an isomorphism, and hence $Y=Y'$.
        Therefore, we have that $Y,Y'$ belong to the same component $\alpha_j = \alpha_{j'}$, and the same holds for $W, W'$. This implies an equality of Hilbert polynomials $P_{i_*(G)} = P_{(i')_*(G')}$.
        We now have a nonzero morphism $G \to G'$ between rank $1$ sheaves on the integral scheme $Y'$ with the same Hilbert polynomial, and hence $g$ is forced to be an isomorphism.
        We noted that $f$ is not an isomorphism, so $h$ cannot be an isomorphism either.
        We conclude from condition (LR) that $h$ is $0$.
        
        From now on, we identify $i_*(G) = i_*(G')$ via $g$.
        Since the composition $F \to F' \to \cI_{W'}$ is zero this time, $f$ factors through $F \to i_*(G)$, and hence $F$ is of the split torsion type: $F \cong i_*(G) \oplus \cI_W$.
        Similar to the argument in the previous subcase, we note that $F \to i_*(G) \subset F'$ and $F \to \gr(0) \to F'$ are both factorizations of $f$, hence their images are equal: $\gr(0) \cong \gr^0(F') = i_*(G)$.
        Now diagram \eqref{eq: diagram for axis-completeness} becomes:
        \begin{equation*}
        \begin{tikzcd}[ampersand replacement = \&]
        	{F} \&\& {\cI_W} \&\& {\Fil_2(F)} \&\& \cdots \\
        	\\
        	{0} \&\& {i_*(G)} \&\& {\Fil^1(F')} \&\& \cdots
        	\arrow[""{name=0, anchor=center, inner sep=0}, hook', from=1-3, to=1-1]
        	\arrow[""{name=1, anchor=center, inner sep=0}, hook', from=1-5, to=1-3]
        	\arrow[""{name=4, anchor=center, inner sep=0}, hook', from=1-7, to=1-5]
        	\arrow[""{name=2, anchor=center, inner sep=0}, hook, from=3-1, to=3-3]
        	\arrow[""{name=3, anchor=center, inner sep=0}, hook, from=3-3, to=3-5]
        	\arrow[""{name=5, anchor=center, inner sep=0}, hook, from=3-5, to=3-7]
        	\arrow["{i_*(G)}"{description}, between={0.2}{0.8}, dotted, no head, from=0, to=2]
        	\arrow["{\gr(1)}"{description}, between={0.2}{0.8}, dotted, no head, from=1, to=3]
        	\arrow["{\gr(2)}"{description}, between={0.2}{0.8}, dotted, no head, from=4, to=5]
        \end{tikzcd}
        \end{equation*}
        We have an isomorphism $\cI_W / \Fil_2(F) \cong \Fil^1(F') / i_*(G)$ from the next step of the filtration.
        The latter quotient is a subsheaf of $F'/i_*(G) \cong \cI_{W'}$, hence $\cI_W / \Fil_2(F)$ is rank $1$ and torsion-free.
        This in turn implies that $\Fil_2(F) = 0$ and $\gr(1) \cong \cI_W$, and therefore $\gr \cong i_*(G) \oplus \cI_W \in \cU_{\alpha,\ell}(K)$. \qedhere
    \end{itemize}
\end{proof}

% -
\subsection{Existence part of the valuative criterion for properness}

\begin{prop} \label{prop: existence valuative criterion}
    The stack $\nps_{\alpha, \ell}$ satisfies the existence part of the valuative criterion of properness with respect to discrete valuation rings.
\end{prop}

\begin{proof}
    Let $R$ be a discrete valuation ring with generic point $\eta$ and special point $s$, and choose $F: \eta \to \nps_{\alpha, \ell}$ corresponding to a sheaf on $X_{\eta}$. If $F$ is of ideal type, then we may write $F \cong \cI_{Z}$ for some length-$\ell$ subscheme $Z \subset X_{\eta}$. If we denote by $\widetilde{Z} \subset X_R$ the flat closure of $Z \subset X_{\eta}$ in $X_R$, then $\widetilde{F}:= \cI_{\widetilde{Z}}$ corresponds to a family $\widetilde{F}: \Spec(R) \to \nps_{\alpha, \ell}$ with generic fiber $\widetilde{F}_{\eta} \cong F$, as desired. 

    We are left to consider the case when $F$ is of torsion type. This means that there exists a short exact sequence
    \[ 0 \to i_*(G) \to F \to \cI_W \to 0\]
    for some $[i: Y \hookrightarrow X_{\eta}], [W \subset X_{\eta}] \in \alpha_j(\eta)$. Let $[\widetilde{i}: \widetilde{Y} \hookrightarrow X_R] \in \alpha_j(R)$ and $[\widetilde{W} \subset X_R] \in \alpha_j(R)$ denote the flat closures of $Y$ and $W$ respectively. Using the existence part of the valuative criterion for properness for the relative moduli of pure sheaves for the family $\widetilde{Y} \to \Spec(R)$ (see e.g. \cite[Prop. 6.16]{torsion-freepaper}), we obtain an $R$-flat family of pure sheaves $\widetilde{G}$ on $\widetilde{Y}$ with generic fiber $\widetilde{G}_{\eta} \cong G$. The local constancy of Hilbert polynomials of the $R$-fibers of $\widetilde{G}$ implies that the special fiber $\widetilde{G}_s$ has rank one. The extension $F$ of $\cI_W$ by $i_*(G)$ yields an element of $\gamma \in \Ext^1(\cI_W, i_*(G))$ which is well-defined up to scaling by a nonzero element of the fraction field $\text{Frac}(R)$. By \cite[\href{https://stacks.math.columbia.edu/tag/0A1J}{0A1J}]{stacks-project} applied to the complexes $\cI_{\widetilde{W}}$, $\widetilde{i}_*(\widetilde{G})$ and the morphism $\pi: X_R \to \Spec(R)$, there is a finite complex $C^{\bullet}$ of free $R$-modules such that for all $R$-algebras $A$ we have an isomorphism of $A$-modules $\Ext^1\left(\cI_{\widetilde{W}}|_{X_A}, \widetilde{i}_*(\widetilde{G})|_{X_A}\right) \cong H^1(C^{\bullet} \otimes_R A)$. The element $\gamma \in \Ext^1(\cI_W, i_*(G)) \cong H^1(C^{\bullet} \otimes_R \text{Frac}(R))$ can be lifted to a cocycle $\widetilde{\gamma} \in C^1 \otimes_R \text{Frac}(R)$. Up to multiplying by a sufficiently large power of the uniformizer of $R$, we can always arrange that $\widetilde{\gamma}$ lies in $C^1 \subset C^1 \otimes_R \text{Frac}(R)$, and that it is a cocycle for the complex $C^{\bullet}$. Then $\widetilde{\gamma}$ represents a cohomology class in $H^1(C^{\bullet}) \cong \Ext^1(\cI_{\widetilde{W}}, \widetilde{i}_*(\widetilde{G}))$, which corresponds to an extension
    \[ 0 \to \widetilde{i}_*(\widetilde{G}) \to \widetilde{F} \to \cI_{\widetilde{W}} \to 0.\]
    with $\widetilde{F}_{\eta} \cong F$. We may view $\widetilde{F}$ as a morphism $\widetilde{F}: \Spec(R) \to \nps_{\alpha, \ell}$, and we have thus obtained the required extension of the point $F \in \nps_{\alpha, \ell}(\eta)$.
\end{proof}

% -
\subsection{Good moduli space}

In this subsection, we collect all our preliminary work in order to construct modifications $\Hilb^{\ell}(X)_{\alpha}$ of the Hilbert scheme $\Hilb^{\ell}(X)$. 

\begin{thm} \label{thm: existence of proper}
    Suppose that $\alpha$ as in \cref{setup: alpha} is of finite type and $(\alpha, \ell)$ satisfies condition (LR) as in \cref{defn: condition dagger}. Then, the stack $\nps_{\alpha, \ell}$ admits a proper good moduli space $\Hilb^{\ell}(X)_{\alpha}$.
\end{thm}

\begin{proof}
    In view of \cref{prop: boundedness}, \cref{prop: local reductivity}, \cref{prop: theta reductivity}, \cref{prop: S completeness} and \cref{prop: existence valuative criterion}, the stack $\nps_{\alpha, \ell}$ satisfies the hypotheses of \cite[Thm. 5.4]{AHLH}, which guarantee the existence of a proper good moduli space.
\end{proof}

Our next goal is to identify an open subset where the modified Hilbert scheme $\Hilb^{\ell}(X)_{\alpha}$ agrees with the usual $\Hilb^{\ell}(X)$.

\begin{notn}
\label{notation: nc-alpha}
    Given $\alpha$ of finite type as in \cref{setup: alpha}, we denote by $\Hilb^{\ell}(X)^{\rnc\alpha} \subset \Hilb^{\ell}(X)$ the open subscheme whose geometric points parameterize subschemes $Z \subset X_K$ for an algebraically closed field $K \supset k$ such that for all $[Y \subset X_K] \in \alpha(K)$ we have that $Z \not\subset Y$ scheme-theoretically.
    The notation $\rnc\alpha$ in the superscript is a shorthand for ``not contained in classes of $\alpha$''.
\end{notn}

\begin{prop} \label{prop: locus where the modified Hilbert scheme is an iso}
    Suppose that $\alpha$ as in \cref{setup: alpha} is of finite type and $(\alpha, \ell)$ satisfies condition (LR) as in \cref{defn: condition dagger}. Then, the morphism of moduli spaces $f: \Hilb^{\ell}(X) \to \Hilb^{\ell}(X)_{\alpha}$ induced by the inclusion $\SCoh^{\rtf}(X)_{\cO_X, P_{\ell}} \hookrightarrow \nps_{\alpha, \ell}$ (see \cref{prop: hilbert scheme as moduli space of sheaves}) is such that the composition $\Hilb^{\ell}(X)^{\rnc\alpha} \hookrightarrow \Hilb^{\ell}(X) \xrightarrow{f} \Hilb^{\ell}(X)_{\alpha}$ is an open immersion.
\end{prop}
\begin{proof}
    We denote by $\SCoh^{\rtf}(X)_{\cO_X, P_{\ell}}^{\rnc\alpha} \subset \SCoh^{\rtf}(X)_{\cO_X, P_{\ell}} \subset \nps_{\alpha, \ell}$ the open preimage of $\Hilb^{\ell}(X)^{\rnc\alpha}$ via the good moduli space morphism $\SCoh^{\rtf}(X)_{\cO_X, P_{\ell}} \to \Hilb^{\ell}(X)$. By \cite[Prop. 2.6]{edidin-rydh-canonical-reduction}, it suffices to show that $\SCoh^{\rtf}(X)_{\cO_X, P_{\ell}}^{\rnc\alpha}$ is contained in the stable locus of both $\SCoh^{\rtf}(X)_{\cO_X, P_{\ell}}$ and $ \nps_{\alpha, \ell}$ (as in \cite[Def 2.5]{edidin-rydh-canonical-reduction}). Note that it is enough to show that $\SCoh^{\rtf}(X)_{\cO_X, P_{\ell}}^{\rnc\alpha}$ is contained in the stable locus of $\nps_{\alpha, \ell}$; it will then follow automatically that $\SCoh^{\rtf}(X)_{\cO_X, P_{\ell}}^{\rnc\alpha}$ is also contained in the stable locus of the open substack $\SCoh^{\rtf}(X)_{\cO_X, P_{\ell}} \subset \nps_{\alpha, \ell}$. In view of \cite[Thm. 4.16(iv)]{alper-good-moduli}, it suffices to show that for all algebraically closed field extensions $K \supset k$ and all points $x \in \SCoh^{\rtf}(X)_{\cO_X, P_{\ell}}^{\rnc\alpha}(K)$, we have that $x$ is closed in the fiber of the moduli space morphism $\nps_{\alpha, \ell} \to \Hilb^{\ell}(X)_{\alpha}$ containing it. After base-changing to $K$ and using \cite[Prop. 4.7(i)]{alper-good-moduli}, we may assume that $K =k$, and so $x$ is defined over the ground field. Then it is sufficient to show that $x$ is a closed point of $\nps_{\alpha, \ell}$, and this follows from \cref{prop: closed points}.
\end{proof}

Before we dive into specific geometric situations, we point out that the morphism $f : \Hilb^\ell(X) \to \Hilb^\ell(X)_\alpha$, constructed in \cref{thm: main theorem on U_alpha}, is not guaranteed to be surjective (see \Cref{example: lack of surjectivity} below).
In our applications, we sidestep this issue by taking the scheme-theoretic image.
However, we now present a simple condition that guarantees surjectivity.

\begin{prop}
\label{prop: condition for surjectivity from Hilbert scheme}
    Let $\alpha$ be as in \Cref{setup: alpha}, and suppose that $\alpha$ is of finite type and $(\alpha,\ell)$ satisfies condition (LR).
    Assume that for all field extensions $K \supset k$ and points $[Y \subset X_K] \in \alpha(K)$, the subvariety $Y \subset X_K$ is smooth with $\Pic(Y) \cong \bZ$.
    Assume in addition that for every pair of subvarieties $Y$ and $W$ from the same connected component $\alpha_j(K)$, we have $\cI_Y \cong \cI_W$. (For example, this happens when $\alpha$ is discrete.)
    Then $f : \Hilb^\ell(X) \to \Hilb^\ell(X)_\alpha$ is surjective.
\end{prop}
\begin{proof}
    It suffices to prove surjectivity after base changing to an algebraically closed field, so for the rest of the proof, we may assume without loss of generality that $k$ is algebraically closed. It is then enough to check surjectivity at the level of $k$-valued points. 
    The $k$-valued points of $\Hilb^\ell(X)_\alpha$ correspond to the $k$-valued points of the stack $\cU_{\alpha,\ell}$ up to closure equivalence \cite[Thm. 4.16(iv)]{alper-good-moduli}.
    By \cref{prop: closed points}, a $k$-point of $\Hilb^\ell(X)_\alpha$ is represented by either an ideal sheaf or a torsion sheaf of split type.
    In the former case, it evidently admits a preimage in $\Hilb^\ell(X)$.
    For the latter case, by the assumption on $\alpha$, we can write $F = i_*(G) \oplus \cI_Y$ for some $[i: Y \hookrightarrow X] \in \alpha(k)$ and a rank one torsion-free sheaf $G$ with Hilbert polynomial $P_G = P_{\cO_Y} - \ell$.
    We will find an ideal sheaf $\cI_Z$ which degenerates to $F$, hence contains $F$ in its closure, and therefore its class in $\Hilb^\ell(X)$ maps to the class of $F$ in $\Hilb^\ell(X)_\alpha$.

    \paragraph{\textbf{Step 1: If $\dim(Y) \geq 2$, then $\det(G) \cong \cO_Y$.}}
    Consider $[\det (G)]$ as a divisor class on $Y$, and let $H$ be the restriction of the hyperplane class on $X$ to $Y$.
    We can write $[\det(G)]=nD$ for a divisor corresponding to the ample generator of $\Pic(Y) \cong \bZ$.
    Let $Z'\subset Y$ be a finite subscheme of length $\ell$, then we know that $P_G = P_{\cO_Y} - \ell = P_{\cI_{Z'\subset Y}}$.
    We then get the following for any $m\in \bZ$:
    \begin{align*}
        0 &= P_{G}(m) - P_{\cI_{Z'\subset Y}}(m)
        \\&= \int_Y (\ch(G) - \ch(\cI_{Z'\subset Y})).e^{mH}.\Td_Y
        , \text{ by Hirzebruch-Riemann-Roch.}
    \end{align*}
    We let $\hot$ denote the degree at least $2$ parts of $\ch(G)$ and $\ch(\cI_{Z'\subset Y})$, and we let $\lot$ denote the smaller degree terms of a polynomial in $m$.
    With this, we continue the calculation above:
    \begin{align*}
        0 &= \int_Y (1 + nD + \hot - 1 - \hot).e^{mH}.\Td_Y
        \\&= \int_Y (nD + \hot).e^{mH}.\Td_Y
        \\&= \frac{n}{(\dim Y - 1)!} D.H^{\dim Y - 1} \cdot m^{\dim Y - 1} + \lot
    \end{align*}
    Since $D$ and $H$ are both ample divisors on $Y$, the intersection $D.H^{\dim Y - 1}$ is nonzero.
    Therefore, in order for the polynomial to vanish, we must have $n=0$, or in other words, $\det G \cong \cO_Y$. 

    \paragraph{\textbf{Step 2: $G$ is isomorphic to an ideal sheaf of a length $\ell$ subscheme in $Y$.}}
    If $\dim Y \geq 2$, then by the previous step, we have $\det G \cong \cO_Y$, but this is also the reflexive hull of $G$ by \cref{lemma: det(rank 1 tf) = double dual}, and so we get an embedding $G\subset \cO_Y$, which means that $G$ is an ideal sheaf.
    Since $P_G = P_{\cO_Y} - \ell$, we conclude that $G$ is an ideal sheaf of a length $\ell$ subscheme in $Y$.

    If $\dim Y = 1$, then our assumptions imply that $Y \cong \mathbb P^1$.
    Since $G$ is rank $1$ torsion-free on a smooth curve, it must be a line bundle.
    The form of its Hilbert polynomial implies that $G \cong \cO_{\mathbb P^1} (-\ell)$, as desired.

    \paragraph{\textbf{Step 3: $F$ is a degeneration of an ideal sheaf of a length $\ell$ subscheme of $X$.}}
    By the previous step, there exists a length $\ell$ subscheme $Z \subset Y$ such that $G \cong \cI_{Z\subset Y}$.
    Then there is an embedding $\cI_Y \subset \cI_Z$ of rank one sheaves on $X$, which yields the short exact sequence
    \[
    0 \to \cI_Y \to \cI_Z \to i_*(\cI_{Z\subset Y}) \to 0.
    \]
    Degenerating this sequence to the split extension yields the result.
\end{proof}

The following examples justify the choice of hypotheses in \Cref{prop: condition for surjectivity from Hilbert scheme}.
\begin{example} \label{example: lack of surjectivity}
    Consider $X := \P^1_k \times C$ over an algebraically closed field, where $C$ is a curve of genus at least $1$, and let $\pi_1 : X \to \P^1_k$, $\pi_2 : X \to C$ denote the projections.
    If one picks $\alpha \cong C$ parameterizing the fibers of $\pi_2$ and denotes by $i : \P^1_s \to X$ the embedding of the fiber $\P^1_s$ over $s \in C$, then $i_*(\cO_{\P^1_k}(-\ell)) \oplus \cI_{\P^1_t}$ is not in the image of $f$ for $s \neq t$.
    If one picks $\alpha \cong \P^1_k$ parameterizing the fibers of $\pi_1$ and denotes by $i : C_s \to X$ the embedding of the fiber $C_s$ over $s \in \P^1$, then one can pick a line bundle $\cL$ on $C_s$ of degree $-\ell$ which is not isomorphic to an ideal sheaf of $\ell$ points by the assumption on the genus.
    Then $i_*(\cL) \oplus \cI_{C_s}$ is not in the image of $f$.
\end{example}

% ---
\section{Applications to contractions} \label{section: applications}

In this section, we spell out in detail the examples that have already been mentioned in the introduction (Abel-Jacobi type contractions, contracting rational curves), and provide some additional examples.
If we have $Y \subset X$ which is a $\P^1$-fibration corresponding to a connected component inside $\Hilb(X)$ and whose normal bundle is fiberwise bounded by $N\in \bZ$ (see \cref{defn: fiberwise bounded normal bundle}), then for $\ell > \max(N,0)$, we can fiberwise contract the $\P^\ell$-fibration in $\Hilb^\ell(X)$ induced by $Y$ (\cref{thm: contraction of projective space bundles}).

% -
\subsection{Abel-Jacobi contractions for Hilbert schemes of surfaces} \label{subsection: Abel-Jacobi contractions}

Let $C$ be a proper geometrically integral curve over $k$. For any given integer $\chi$, the stack $\SCoh^{\rtf}(C)_{(1,\chi)}$ of rank $1$ torsion-free sheaves of Euler characteristic $\chi$ is a $\mathbb{G}_m$-gerbe over the compactified Jacobian $\overline{\Pic}(C)_{\chi}$, which is a projective scheme.
The corresponding Abel-Jacobi morphism $AJ_{\ell}: \Hilb^{\ell}(C) \to \overline{\Pic}(C)_{\chi(\cO_C)-\ell}$ is induced by the map $\Hilb^{\ell}(C) \to \SCoh^{\rtf}(C)_{(1,\chi)}$ that sends a length $\ell$ subscheme $Z \subset C$ to its ideal sheaf $\cI_Z$ (see \cite[(5.16)]{altman-kleiman-compactifying-picard} and \cite[pg. 102]{altman-kleiman-compactifying-picard} setting their $F$ to $\cO_C$).
If the curve $C$ is Gorenstein, then for a sufficiently large $\ell$, the Abel-Jacobi morphism exhibits $\Hilb^{\ell}(C)$ as a projective bundle over $\overline{\Pic}(C)_{\chi(\cO_C)-\ell}$ \cite[Thm. 8.6]{altman-kleiman-compactifying-picard}. These constructions apply verbatim in the setting of families of integral curves: given a flat proper morphism $\cC \to B$ of finite type schemes with geometrically integral fibers of dimension 1, there are relative versions $\Hilb^{\ell}(\cC/B) \to B$ and $\overline{\Pic}(\cC/B)_{\chi(\cO) - \ell} \to B$ of the Hilbert scheme and compactified Jacobian, and a relative version of the Abel-Jacobi morphism $AJ_{\ell}: \Hilb^{\ell}(\cC/B) \to \overline{\Pic}(\cC/B)_{\chi(\cO)-\ell}$.

Let $X$ be a smooth projective surface, and let $\alpha \subset \Hilb(X)$ be a finite union of connected components of the Hilbert scheme such that every geometric point $[C \subset X_K] \in \alpha$ is an integral subcurve. We denote by $I(\alpha)$ the maximum among all intersections numbers $C_1 \cdot C_2$ as we range over all field extensions $K$ and all pairs of points $[C_1 \subset X_K], [C_2 \subset X_K] \in \alpha(K)$.
\begin{lemma} \label{lemma: radicial morphism}
    Let $X$ be a smooth projective surface, and let $\alpha \subset \Hilb(X)$ be a finite union of connected components of the Hilbert scheme such that the corresponding universal family $\cC \to \alpha$ parameterizes geometrically integral curves. For all $\ell> I(\alpha)$, the induced morphism $\psi: \Hilb^{\ell}(\cC/\alpha) \to \Hilb^{\ell}(X)$ given on $T$-points by
    \[
        (a : T \to \alpha, [Z \subset \cC_a]) \mapsto [Z \subset \cC_a \subset X_T]
    \]
    is finite and radicial.
\end{lemma}
\begin{proof}
    Since $\psi$ is proper, it suffices to show that it is injective at the level of geometric points. In other words, we need to show that for every field extension $K$ and any length $\ell$ subscheme $[Z \subset X_K] \in \Hilb^{\ell}(X)(K)$, there is at most one $[C \subset X_K] \in \alpha(K)$ such that $Z \subset C$.
But indeed, if there were two distinct $[C_1 \subset X_K], [C_2 \subset X_K] \in \alpha(K)$ containing $Z$, then we would get $C_1 \cdot C_2 \geq \ell > I(\alpha)$, a contradiction.
\end{proof}

Note that the image of $\psi: \Hilb^{\ell}(\cC/\alpha) \to \Hilb^{\ell}(X)$ in \Cref{lemma: radicial morphism} is closed and is equal to the union $\bigcup_{t \in \alpha} \Hilb^{\ell}(\cC_t)$ of the Hilbert schemes of all subcurves parameterized by $\alpha$.

\begin{thm} \label{prop: abel jacobi contractions}
    Let $X$ be a smooth projective surface, and let $\alpha \subset \Hilb(X)$ be a finite union of connected components parameterizing integral subcurves.
    Denote by $\cC \to \alpha$ the universal family of integral subcurves over $\alpha$.
    Then, for all positive integers $\ell>I(\alpha)$, there is a surjective morphism $f: \Hilb^{\ell}(X) \to Y_{\ell}$ to a proper algebraic space $Y_{\ell}$ satisfying the following: 
    \begin{enumerate}
        \item
        \label{item: birationality with Hilbert scheme for surfaces}
        The morphism $f$ is an isomorphism when restricted to the open complement of the image of $\psi: \Hilb^{\ell}(\cC/\alpha) \to \Hilb^{\ell}(X)$.
        \item
        \label{item: locus contracted by AJ}
        If we denote by $H_{\ell} \subset \overline{\Pic}(\cC/\alpha)_{\chi(\cO)-\ell}$ 
        the scheme-theoretic image of the Abel-Jacobi morphism $AJ_{\ell}: \Hilb^{\ell}(\cC/\alpha) \to \overline{\Pic}(\cC/\alpha)_{\chi(\cO)-\ell}$, then there is closed immersion $\theta: H_{\ell} \hookrightarrow Y_{\ell}$ such that the following diagram commutes
        \[
        \begin{tikzcd}
         \Hilb^{\ell}(\cC/\alpha) \ar[d, "AJ_{\ell}"]
          \ar[r, "\psi"] & \Hilb^{\ell}(X) \ar[d, "f"]\\ H_{\ell} \ar[r, "\theta"] & Y_{\ell}
        \end{tikzcd}
        \]
    \end{enumerate}
\end{thm}

\begin{proof}
We first show that the pair $(\alpha, \ell)$ satisfies property (LR). Part (1) in (LR) is automatically satisfied by \Cref{remark: condition 1 automatic dimension}, so we just need to check part (2) of (LR).
Let $K \supset k$ be an algebraically closed extension, and fix two points $a, b \in \alpha(K)$ and a rank one torsion-free sheaf $G$ on $\cC_{b}$ with Hilbert polynomial $P_{\cO_{\cC_{b}}}-\ell$. After base-changing to $K$, we may assume without loss of generality that $K=k$. Set $i: \cC_{b} \hookrightarrow X$. Using adjunction, we get $\Hom(\cI_{\cC_{a}}, i_*(G)) = \Hom(i^*(\cI_{\cC_a}), G)$. The pullback $i^*(\cI_{\cC_a})$ is a line bundle with Hilbert polynomial $P_{\cO_{\cC_{b}}}-I(\alpha)$. Using the assumption $\ell>I(\alpha)$, it follows that we have an inequality $P_{i^*(\cI_{\cC_a})} > P_G$ of Hilbert polynomials of rank one torsion-free sheaves on the integral curve $\cC_{b}$. This implies that $\Hom(i^*(\cI_{\cC_a}), G)=0$, thus concluding the proof that $(\alpha, \ell)$ satisfies (LR).

In view of \cref{thm: main theorem on U_alpha} and \cref{prop: locus where the modified Hilbert scheme is an iso}, there is a morphism $f: \Hilb^{\ell}(X) \to \Hilb^{\ell}(X)_{\alpha}$ to a proper algebraic space $\Hilb^{\ell}(X)_{\alpha}$, which is an isomorphism when restricted to $\Hilb^{\ell}(X)^{\rnc\alpha} := \Hilb^{\ell}(X) \setminus \Hilb^{\ell}(\cC/\alpha)$ (cf. \cref{notation: nc-alpha}).
This proves part \ref{item: birationality with Hilbert scheme for surfaces}.

For \ref{item: locus contracted by AJ}, we set $Y_{\ell}$ to be the closed scheme-theoretic image of $\Hilb^{\ell}(X)$ in $\Hilb^{\ell}(X)_{\alpha}$.
To conclude the proof, we are left to find a closed immersion $\theta: H_{\ell} \hookrightarrow Y_{\ell}$ such that the diagram in part \ref{item: locus contracted by AJ} is commutative.
For this, it is sufficient to find a closed immersion $\theta: \overline{\Pic}(\cC/\alpha)_{\chi(\cO_{\cc})-\ell} \hookrightarrow \Hilb^{\ell}(X)_{\alpha}$ such that the following diagram is commutative:
\begin{equation}
\label{equation 1: proof of abel jacobi contraction}
    \begin{tikzcd}[ampersand replacement = \&]
    \Hilb^{\ell}(\cC/\alpha) \ar[d, "AJ_{\ell}"]
    \ar[r, "\psi"] \& \Hilb^{\ell}(X) \ar[d, "f"]\\ \overline{\Pic}(\cC/\alpha)_{\chi(\cO_{\cc})-\ell} \ar[r, "\theta"] \& \Hilb^{\ell}(X)_{\alpha}
    \end{tikzcd}
\end{equation}
Let us denote by $\overline{\SPic}(\cC/\alpha) \to \alpha$ the relative stack of rank $1$ torsion-free sheaves with Euler characteristic $\chi(\cO)-\ell$ for the family of integral curves $\cC \to \alpha$ (here we omit the subscript $\chi(\cO_{\cC})-\ell$ for simplicity).
We have that $\overline{\SPic}(\cC/\alpha) \to \overline{\Pic}(\cC/\alpha)_{\chi(\cO_{\cC})-\ell}$ is a $\mathbb{G}_m$-gerbe.
Consider the stack $\overline{\SPic}(\cC/\alpha) \times B\mathbb{G}_m$ 
that sends a $k$-scheme $T$ to the groupoid of tuples $(a: T \to \alpha, G, L)$, where, denoting by $\cC_T$ the pullback of $\cC \to \alpha$ along $a : T \to \alpha$, we take $G$ to be a $T$-flat family of torsion-free sheaves on $\cC_T$, and $L$ is a line bundle on $T$.
The morphisms in this groupoid are isomorphisms of $G$ and $L$.
We define the morphism $\widetilde{\theta}: \overline{\SPic}(\cC/\alpha) \times B\mathbb{G}_m \to \nps_{\alpha, \ell}$ at the level of $T$-points:
\[
    \widetilde{\theta} : (a,G,L) \mapsto i_*(G) \oplus \cI_{\cC_a} \otimes \pi^*(L),
\]
where $i: \cC_a \hookrightarrow X_T$ denotes the closed embedding with ideal sheaf $\cI_{\cC_a}$ and $\pi: X_T \to T$ denotes the structure morphism. Note that this morphism is representable: indeed for any $T$-points as above the automorphisms of both $(a,G,L)$ and $i_*(G) \oplus \cI_{\cC_a} \otimes \pi^*(L)$ are isomorphic to $\mathbb{G}_m^2(T)$, 
and the functor is easily seen to induce an isomorphism of automorphisms groups.
Furthermore, we claim that $\widetilde{\theta}$ is a monomorphism. Indeed, let $T$ be a $k$-scheme, and let $(a, G, L), (b,F, M) \in \overline{\SPic}(\cC/\alpha) \times B\mathbb{G}_m(T)$. Let us denote by $i_a: \cC_a \hookrightarrow X_T$ and $i_b: \cC_b \hookrightarrow X_T$ the corresponding closed immersions.
Suppose that $\widetilde{x} = (i_a)_*(G) \oplus \cI_{\cC_a} \otimes \pi^*(L)$ and $\widetilde{y} = (i_b)_*(F) \oplus \cI_{\cC_b} \otimes \pi^*(M)$ are isomorphic. Note that $(i_a)_*(G) \oplus \cI_{\cC_a} \otimes \pi^*(L) \twoheadrightarrow \cI_{\cC_a} \otimes \pi^*(L)$ and $(i_b)_*(F) \oplus \cI_{\cC_b} \otimes \pi^*(M) \twoheadrightarrow \cI_{\cC_b} \otimes \pi^*(M)$ are uniquely determined as the unique largest quotients which are $T$-flat families of torsion-free sheaves on $X_T \to T$.
Hence, it follows that we have individual isomorphisms $(i_a)_*(G) \cong (i_b)_*(F)$ and $\cI_{\cC_a} \otimes \pi^*(L) \cong \cI_{\cC_b} \otimes \pi^*(M)$.
From the first isomorphism it follows that we have equality of scheme-theoretic supports $\cC_a = \cC_b \subset X_T$ of $(i_a)_*(G) \cong (i_b)_*(F)$ and also $G \cong F$ in $\cC_a = \cC_b$.
On the other hand, by taking determinant of the isomorphism $\cI_{\cC_a} \otimes \pi^*(L) \cong \cI_{\cC_b} \otimes \pi^*(M)$ and tensoring with $\det(\cI_{\cC_a})^{\vee} = \det(\cI_{\cC_b})^{\vee}$, we conclude that $\pi^*(L) \cong \pi^*(M)$. It follows then from the projection formula that $L \cong \pi_*(\pi^*(L)) \cong \pi_*(\pi^*(M)) \cong M$.
Hence $\widetilde{x} \cong \widetilde{y}$, concluding the proof that $\widetilde{\theta}$ is a monomorphism.

We claim that $\widetilde{\theta}$ is a closed immersion.
In order to see this, by \cite[\href{https://stacks.math.columbia.edu/tag/04XV}{04XV}]{stacks-project} it suffices to show that it is universally closed. We may factor $
\widetilde{\theta}$ as a composition $\overline{\SPic}(\cC/\alpha) \times B\mathbb{G}_m \xrightarrow{h} \bigsqcup_j\overline{\SPic}(\cC/\alpha_j) \times B\mathbb{G}_m \times \alpha_j \xrightarrow{\tau} \nps_{\alpha, \ell}$, where $h$ is given at the level of points by
\[ (a,G,L) \mapsto (a,G,L,a) ,\]
and $\tau$ is given at the level of $T$-points by 
\[ (a,G,L,b) \mapsto (i_a)_*(G) \oplus \cI_{\cC_b} \otimes \pi^*(L).\]
Since the morphism $h$ is a closed immersion, we are reduced to showing that $\tau$ is universally closed. Thus, we need to show that for all schemes $T$ and all families of sheaves $F: T \to \nps_{\alpha, \ell}$, the image of $T \times_{\nps_{\alpha, \ell}} \bigsqcup_j\overline{\SPic}(\cC/\alpha_j) \times B\mathbb{G}_m \times \alpha_j \to T$ is closed.
However, note that this image coincides exactly with the points of $T$ where the fiber $F_t$ is of split torsion type.
Alternatively, this can be described as the points $t \in T$ such that $F_t$ has endomorphisms of dimension at least 2 (and in fact equal to 2) over the residue field.
By upper semicontinuity of dimension of endomorphisms for the $T$-flat family $F$ (which follows from \cite[Cor. 7.7.8.2]{egaiii_2}), we conclude that the image is closed, thus concluding the proof that $\widetilde{\theta}$ is a closed immersion. 

By \cite[Lem. 4.12, 4.14]{alper-good-moduli}, the immersion $\widetilde{\theta}$ induces a closed immersion of good moduli spaces $\theta: \overline{\Pic}(\cC/\alpha)_{\chi(\cO_{\cc})-\ell} \hookrightarrow \Hilb^{\ell}(X)_{\alpha}$.
We are left to show that the relevant diagram \eqref{equation 1: proof of abel jacobi contraction} commutes. Consider the morphism $\widetilde{AJ}_{\ell}: \Hilb^{\ell}(\cC/\alpha) \to \overline{\SPic}(\cC/\alpha) \times B\mathbb{G}_m$ given at the level of $T$-points by $(a, [Z \subset \cC_a]) \mapsto (a, \cI_{Z\subset \cC_a}, \cO_T)$.
Note also that the finite radicial morphism
$\psi : \Hilb^{\ell}(\cC/\alpha) \to \Hilb^{\ell}(X)$ can be factored as $\Hilb^{\ell}(\cC/\alpha) \xrightarrow{\widetilde{\psi}} \SCoh^{\rtf}(X)_{\cO_X, P_{\ell}} \to \Hilb^{\ell}(X)$, where $\widetilde{\psi}$ sends $(a, [Z \subset \cC_a])$ to $\cI_{Z \subset X_T}$.
Set $\gamma: \Hilb^{\ell}(\cC/\alpha) \xrightarrow{\widetilde{\psi}} \SCoh^{\rtf}(X)_{\cO_X, P_{\ell}} \hookrightarrow \nps_{\alpha, \ell}$
and $\xi: \Hilb^{\ell}(\cC/\alpha) \xrightarrow{\widetilde{AJ}_{\ell}} \overline{\SPic}(\cC/\alpha) \times B\mathbb{G}_m \xrightarrow{\widetilde{\theta}} \nps_{\alpha, \ell}$.
If we set $\mu: \nps_{\alpha, \ell} \to \Hilb^{\ell}(X)_{\alpha}$ to be the good moduli space morphism, then the commutativity of the diagram \eqref{equation 1: proof of abel jacobi contraction} amounts to showing that $\mu \circ \gamma = \mu \circ \xi$.

We may think of $\gamma$ and $\xi$ as families of sheaves on $X_{\Hilb^{\ell}(\cC/\alpha)} \to \Hilb^{\ell}(\cC/\alpha)$. More precisely, if $(a_{\univ}, [Z_{\univ} \subset \cC\times_{\alpha} \Hilb^{\ell}(\cC/\alpha)])$ denotes the universal family on $\Hilb^{\ell}(\cC/\alpha)$, then $\gamma$ corresponds to the sheaf $F_{\gamma} := \cI_{Z_{\univ}\subset X_{\Hilb^{\ell}(\cC/\alpha)}}$ on $X_{\Hilb^{\ell}(\cC/\alpha)}$ and $\xi$ corresponds to the sheaf 
\[F_{\xi}:= i_*(\cI_{Z_{\univ} \subset \cC\times_{\alpha} \Hilb^{\ell}(\cC/\alpha)}) \oplus \cI_{\cC\times_{\alpha} \Hilb^{\ell}(\cC/\alpha) \subset X_{\Hilb^{\ell}(\cC/\alpha)}}.\]
The inclusion of ideal sheaves $\cI_{\cC\times_{\alpha} \Hilb^{\ell}(\cC/\alpha) \subset X_{\Hilb^{\ell}(\cC/\alpha)}} \hookrightarrow \cI_{Z_{\univ}\subset X_{\Hilb^{\ell}(\cC/\alpha)}}$ induces a short exact sequence of sheaves on  $X_{\Hilb^{\ell}(\cC/\alpha)}$
\[ 0 \to \cI_{\cC\times_{\alpha} \Hilb^{\ell}(\cC/\alpha) \subset X_{\Hilb^{\ell}(\cC/\alpha)}} \to \cI_{Z_{\univ}\subset X_{\Hilb^{\ell}(\cC/\alpha)}} \to i_*(\cI_{Z_{\univ} \subset \cC\times_{\alpha} \Hilb^{\ell}(\cC/\alpha)}) \to 0.\]
We may think of this as a two-step filtration of $F_{\gamma}$ with the associated graded sheaf isomorphic to $F_{\xi}$. By assigning integer weights to the steps of the filtration, we may apply the Rees construction \cite[Prop. 1.0.1]{halpernleistner2018structure} to obtain a morphism of stacks $\widetilde{\gamma}: (\mathbb{A}^1_k/\mathbb{G}_m) \times \Hilb^{\ell}(\cC/\alpha) \to \nps_{\alpha, \ell}$ such that $\widetilde{\gamma}(1) \cong \gamma$ and $\widetilde{\gamma}(0) \cong \xi$.
Note that the stack $\mathbb{A}^1_k/\mathbb{G}_m$ has good moduli space $\Spec(k)$, and so the stack $(\mathbb{A}^1_k/\mathbb{G}_m) \times \Hilb^{\ell}(\cC/\alpha)$ has good moduli space $\Hilb^{\ell}(\cC/\alpha)$.
By the universal property of good moduli spaces \cite[Thm. 6.6]{alper-good-moduli}, the composition $(\mathbb{A}^1_k/\mathbb{G}_m) \times \Hilb^{\ell}(\cC/\alpha) \xrightarrow{\widetilde{\gamma}} \nps_{\alpha, \ell} \xrightarrow{\mu} \Hilb^{\ell}(X)_{\alpha}$ factors through $(\mathbb{A}^1_k/\mathbb{G}_m) \times \Hilb^{\ell}(\cC/\alpha) \to \Hilb^{\ell}(\cC/\alpha)$.
Therefore, by restricting to $0 \times\Hilb^{\ell}(\cC/\alpha)$ and $1 \times \Hilb^{\ell}(\cC/\alpha)$ we get $\mu \circ \widetilde{\gamma}(1) = \mu \circ \widetilde{\gamma}(0)$.
Using $\widetilde{\gamma}(1) \cong \gamma$ and $\widetilde{\gamma}(0) \cong \xi$, we obtain the desired equality $\mu \circ \gamma = \mu \circ \xi$.
\end{proof}

\subsection{Examples of contractions in Hilbert schemes of surfaces}

\begin{example}[Contractions of $\Hilb^{\ell}(\mathbb{P}^2_k)$] \label{example: contractions of hilbert scheme of projective plane}
    Let $\alpha \cong \mathbb{P}^2_k \subset \Hilb(\mathbb{P}^2_k)$ denote the component of the Hilbert scheme parameterizing lines in $\mathbb{P}^2_k$. Using the terminology of \cref{subsection: Abel-Jacobi contractions}, we have $I(\alpha) =1$. In this case we have that $\Hilb^{\ell}(\cC/\alpha) \to \alpha$ is a $\mathbb{P}^\ell$-bundle, and $\overline{\Pic}(\cC/\alpha)_{\chi(\cO)-\ell} = \Pic(\cC/\alpha)_{1-\ell} \cong \alpha$. The corresponding morphism $f: \Hilb^{\ell}(\mathbb{P}^2_k) \to Y_{\ell}$ guaranteed by \cref{prop: abel jacobi contractions} contracts the $\mathbb{P}^{\ell}$-fibers of the embedded $\mathbb{P}^{\ell}$-bundle $\Hilb^{\ell}(\cC/\alpha) \hookrightarrow \Hilb^{\ell}(X)$. 
    
    Two points $Z,Z' \in \Hilb^{\ell}(\mathbb{P}^2_k)$ have the same image in $Y_{\ell}$ if and only if there is a line $\mathbb{P}^1 \subset \mathbb{P}^2$ containing both $Z$ and $Z'$.
    When $\ell=2$ we have $Y_{\ell} \cong \alpha \cong \mathbb{P}^2_k$ and the morphism $\Hilb^2(\mathbb{P}^2_k) \to \mathbb{P}^2_k$ sends a length-$2$ subscheme $Z \subset \mathbb{P}^2_k$ to the unique line in $\mathbb{P}^2$ containing it.
    For $\ell >2$, the morphism $f: \Hilb^{\ell}(\mathbb{P}^2_k) \to Y_{\ell}$ is birational.

    We note that in this example the algebraic space $Y_{\ell}$ is known to be projective, and $\nps_{\alpha, \ell}$ is a stack of Bridgeland semistable complexes, see the rank one walls described in \cite[\S 9]{MMP_hilbert_scheme_bridgeland} and the examples in \cite[\S10]{MMP_hilbert_scheme_bridgeland}.
\end{example}

We thank Yoonjoo Kim for bringing the following example to our attention.
\begin{example}[A Mukai flop] \label{example: mukai flop contractions}
    Let $X \to \mathbb{P}^1_k$ be an elliptic $K3$ surface with two singular fibers $X_a$, $X_b$ of type $I_3$ as considered in \cite[Ex. 1.7(ii)]{namikawa_deformation_symplectic}. The curve $X_a$ is a cycle consisting of three copies of $\mathbb{P}^1_k$, say $C_1, C_2, C_3$. The inclusions $[C_i \subset X]$ are isolated points in the Hilbert scheme $\Hilb(X)$, we denote by $\alpha_1, \alpha_2 , \alpha_3 \subset \Hilb(X)$ the corresponding connected components and set $\alpha = \bigsqcup_{i=1}^3 \alpha_i$. For $\ell \geq 2$, the three subschemes $\Hilb^{\ell}(C_i) \cong \mathbb{P}^{\ell}_k \subset \Hilb^{\ell}(X)$ are pairwise disjoint. For any $\ell \geq 2$, the pair $(\alpha, \ell)$ satisfies (LR), and the corresponding morphism $\Hilb^{\ell}(X) \to \Hilb^{\ell}(X)_{\alpha}$ is a contraction of $\bigsqcup_{i=1}^3 \Hilb^{\ell}(C_i)$. As explained in \cite[Ex. 1.7(ii)]{namikawa_deformation_symplectic}, the proper algebraic space $\Hilb^{2}(X)_{\alpha}$ is not a scheme. We note that $\Hilb^{2}(X)_{\alpha}$ cannot be obtained as a moduli of Bridgeland semistable objects on a smooth projective variety. Indeed, all such moduli spaces admit a strictly nef line bundle (see \cite{bayer_macri_line_bundle} and \cite[Lem. 5.2]{takajakka_bridgeland}), whereas it is proven in \cite[Ex. 1.7(ii)]{namikawa_deformation_symplectic} that $\Hilb^{2}(X)_{\alpha}$ contains a curve class which has trivial intersection with any divisor class. We thank Emanuele Macr\`i for communicating this argument to us.
\end{example}

\begin{example}[Contractions of Hilbert schemes of elliptic K3 surfaces] Let $X$ be a K3 surface equipped with a flat morphism $\pi: X \to \mathbb{P}^1_k$ such that all fibers are geometrically integral. We may think of $X \to \mathbb{P}^1_k$ as a flat family of subcurves of $X$ parameterized by $\mathbb{P}^1_k$, and hence we get a closed immersion $\alpha:= \mathbb{P}^1_k \hookrightarrow \Hilb(X)$. Using the projection formula, we get that $H^0(X, \pi^*(\cO_{\mathbb{P}^1_k}(1))$ has dimension $2$, and hence $\mathbb{P}^1_k = \alpha$ is an open and closed component of the Hilbert scheme of divisors in $\Hilb(X)$.

In this case, using the language of \cref{subsection: Abel-Jacobi contractions}, we have $I(\alpha)=0$. For all $\ell>0$, we have that $\overline{\Pic}_{\chi(\cO)-\ell} = \overline{\Pic}_{-\ell} \to \mathbb{P}^1$ is a family of integral curves, and the corresponding Abel-Jacobi morphism $AJ_{\ell}: \Hilb^{\ell}(X/\mathbb{P}^1) \to \overline{\Pic}_{-\ell}$ is a $\mathbb{P}^{\ell-1}$-bundle \cite[Thm. 8.6]{altman-kleiman-compactifying-picard}. The corresponding morphism $\Hilb^{\ell}(X) \to Y_{\ell}$ guaranteed by \cref{prop: abel jacobi contractions} contracts the fibers of $\Hilb^{\ell}(X/\mathbb{P}^1) \to \overline{\Pic}_{-\ell}$ of the closed subscheme $\Hilb^{\ell}(X/\mathbb{P}^1) \hookrightarrow \Hilb^{\ell}(X)$. 

Consider the composition $\pi_{\ell}: \Hilb^{\ell}(X) \to X^{\ell}/\!/S_{\ell} \to (\mathbb{P}^1_k)^{\ell}/\!/S_{\ell} \cong \mathbb{P}^{\ell}_k$, where the left-most morphism is the Hilbert-Chow morphism. There is a natural closed embedding $\Delta: \mathbb{P}^1_k \to (\mathbb{P}^1_k)^{\ell}/\!/S_{\ell} \cong \mathbb{P}^{\ell}_k$, and we get a commutative diagram of spaces:
           \[
\begin{tikzcd}
 \Hilb^{\ell}(X/\mathbb{P}^1) \ar[d, "AJ_{\ell}"]
  \ar[r] & \Hilb^{\ell}(X) \ar[d, "f"] \ar[ddr, "\pi_{\ell}"] &\\ \overline{\Pic}_{-\ell} \ar[d] \ar[r]  & Y_{\ell} & \\ \mathbb{P}^1_k \ar[rr, "\Delta"]
  & & \mathbb{P}^{\ell}_k.
\end{tikzcd}
\]
Since all the fibers of the birational modification $f$ are contained inside fibers of $\pi_{\ell}$, it follows that, after perhaps replacing $Y_{\ell}$ with its normalization, we obtain a factorization $\pi_{\ell}: \Hilb^{\ell}(X) \xrightarrow{f} Y_{\ell} \to \mathbb{P}^{\ell}_k$.
\end{example}

The following example was suggested to us by Rob Lazardsfeld.
\begin{example}[Contractions of Hilbert schemes of general hypesurfaces]
    Let $k$ be an algebraically closed field of characteristic $0$, and fix $d \geq 4$.
    Let $X$ be a general hypersurface of degree $d$ inside $\mathbb{P}^3_k$.
    By the Noether-Lefschetz theorem, the Picard group $\Pic(X)$ is generated by the restriction of $\cO_{\mathbb{P}^3_k}(1)$.
    We denote by $\alpha \subset \Hilb(X)$ the component of the Hilbert scheme parameterizing hyperplane section in $X$ (i.e. curves $C \subset X$ with ideal sheaf $\cO_{\mathbb{P}^3_k}(-1)|_X$).
    We note that $\alpha \cong \mathbb{P}(H^0(\mathbb{P}^3_k, \cO_{\mathbb{P}^3_k}(1))) \cong \mathbb{P}^3_k$ is the dual projective space, and every curve $[C \subset X] \in \alpha$ is integral because the divisor class $\cO_{\mathbb{P}^3_k}(1)|_X$ is not divisible.
    In this case, we have $I(\alpha)=d$ (as in \Cref{subsection: Abel-Jacobi contractions}), and hence for all $\ell >d$ \Cref{prop: abel jacobi contractions} yields a corresponding Abel-Jacobi-type contraction 
    \[
        \begin{tikzcd}
         \Hilb^{\ell}(\cC/\mathbb{P}^3_k) \ar[d, "AJ_{\ell}"]
          \ar[r, "\psi"] & \Hilb^{\ell}(X) \ar[d, "f"]\\ H_{\ell} \ar[r, "\theta"] & Y_{\ell}
        \end{tikzcd},
    \]
    where $\pi: \cC \to \mathbb{P}^3_k \cong \alpha$ is the universal family of complete intersection curves of genus $\binom{d-1}{2}$ over $\alpha$. 
\end{example}

% -
\subsection{Fiberwise contractions of \texorpdfstring{$\mathbb{P}^{\ell}$}{Pell}-fibrations} \label{subsection: contraction of projective fibrations}

Let $Y \hookrightarrow X$ be a smooth closed subscheme of a smooth scheme $X$.
Assume that $Y$ is a $\mathbb{P}^1$-fibration, which means that there is a smooth proper morphism $Y \to B$ such that every geometric fiber is isomorphic to $\mathbb{P}^1$. 

\begin{defn}
\label{defn: fiberwise bounded normal bundle}
    We say that the normal bundle of $Y \hookrightarrow X$ is \emph{fiberwise bounded} by an integer $N \in \mathbb{Z}$ if for all geometric points $b \to B$, the pullback $(\cI_Y|_{Y_b})^\vee$ is isomorphic to $\bigoplus_{j=1}^{\dim(X)-\dim(Y)} \cO_{\mathbb{P}^1_b}(n_j)$ with $n_j\leq N$ for all $j$.
\end{defn}

In this case, we have $\overline{\Pic}(Y/B)_{\chi(\cO_Y)-\ell} = B$, and the relative Hilbert scheme of points $\Hilb^{\ell}(Y/B) \to B$ associated to the family of integral curves $Y \to B$ is a $\mathbb{P}^{\ell}$-fibration over $B$.
In the following \cref{lemma: closed immersion of Hilbert schemes}, we prove that there is a natural closed immersion $\Hilb^{\ell}(Y/B) \hookrightarrow \Hilb^{\ell}(X)$.

\begin{lemma} \label{lemma: closed immersion of Hilbert schemes}
    Let $\cC \to B$ be a flat family of curves and suppose that there is a closed immersion $\cC \hookrightarrow X$ into a projective $k$-scheme $X$.
    Then for all positive integers $\ell$, the induced morphism $\Hilb^{\ell}(\cC/B) \to \Hilb^{\ell}(X)$ given on $T$-points by
    \[
        (b : T \to B, [Z \subset \cC_b]) \mapsto [Z \subset \cC_b \subset X_T]
    \]
    is a closed immersion.
\end{lemma}

\begin{proof}
    We note that the morphism $\Hilb^{\ell}(\cC/B) \to \Hilb^{\ell}(X)$ factors as a composition $\Hilb^{\ell}(\cC/B) \to \Hilb^{\ell}(\cC) \hookrightarrow \Hilb^{\ell}(X)$, where the right-most morphism $\Hilb^{\ell}(\cC) \hookrightarrow \Hilb^{\ell}(X)$ is a closed immersion. Hence, it suffices to show that $\Hilb^{\ell}(\cC/B) \to \Hilb^{\ell}(\cC)$ is a closed immersion. This morphism is immediately seen to be proper, since both the target and source are proper.
    Hence, we just need to show that it is a monomorphism.
    This follows from looking at the functors of points: for any Noetherian $k$-scheme $T$, we have that $\Hilb^{\ell}(\cC)(T)$ is the set of closed immersions $Z \hookrightarrow \cC_T$, where $Z \to T$ is finite and flat of constant degree $\ell$, whereas $\Hilb^{\ell}(\cC/B)(T)$ is the subset of those $Z \hookrightarrow \cC_T$ such that the composition $Z \to \cC_T \to B_T$ factors (uniquely) through the faithfully flat finite morphism $Z \to T$.
\end{proof}

\begin{lemma} \label{lemma: closed immersion for P^1-fibrations}
    Let $X$ be a projective $k$-scheme, and let $Y \hookrightarrow X$ be a smooth subscheme which admits a $\mathbb{P}^1$-fibration $\pi: Y \to B$. Then the morphism $B \to \Hilb(X)$ induced by the $B$-family of subschemes $Y \to B$ is a closed immersion.
\end{lemma}
\begin{proof}
    The morphism $B \to \Hilb(X)$ is injective on geometric points; since $Y \hookrightarrow X$ is a closed subscheme, distinct geometric points $b, b' \in B$ yield distinct subcurves $[Y_b \subset X_b], [Y_{b'} \subset X_{b'}]$ in $\Hilb(X)$. To show that it is a closed immersion, we may work separately on each connected component of the smooth scheme $B$, and hence we may assume that $B$ is connected. Let $P$ be the common Hilbert polynomial of all fibers of $Y \to B$. Then we have that $B \to \Hilb(X)$ factors through the open and closed subscheme $\Hilb(X)_P \subset \Hilb(X)$ parameterizing subschemes of $X$ with Hilbert polynomial $P$. We have a factorization $B = \Hilb(Y/B)_P \to \Hilb(Y)_P \hookrightarrow \Hilb(X)_P$, where the identification $B \xrightarrow{\sim} \Hilb(Y/B)_P$ sends a $B$-scheme $T \to B$ to the identity $T$-family $[Y_T \hookrightarrow Y_T]$ of subschemes with fiberwise Hilbert polynomial $P$. A similar argument as in the proof of \Cref{lemma: closed immersion of Hilbert schemes} shows that $\Hilb(Y/B)_P \hookrightarrow \Hilb(Y)_P$ is a closed immersion. Since $\Hilb(Y)_P \hookrightarrow \Hilb(X)_P$ is evidently a closed immersion, we conclude that the composition $B = \Hilb(Y/B)_P \hookrightarrow \Hilb(Y)_P \hookrightarrow \Hilb(X)_P$ is a closed immersion.
\end{proof}

\begin{thm}
\label{thm: contraction of projective space bundles}
    Let $X$ be a smooth projective geometrically connected scheme over $k$ of dimension at least $2$. Let $Y \hookrightarrow X$ be a smooth closed subscheme of $X$ which admits a $\mathbb{P}^1$-fibration $\pi: Y \to B$ and such that the following hold:
    \begin{enumerate}
        \item[(A)] the normal bundle of $Y$ is fiberwise bounded by a given integer $N$ in the sense of \cref{defn: fiberwise bounded normal bundle}, and
        \item[(B)] the closed immersion $B \hookrightarrow \Hilb(X)$ induced by the family of subcurves $Y \to B$ (as in \Cref{lemma: closed immersion for P^1-fibrations}) has open image in $\Hilb(X)$. (Hence the image is a union of connected components.)
    \end{enumerate}
    Then, for all positive integers $\ell >N$, there is a surjective morphism $g: \Hilb^{\ell}(X) \to M_{\ell}$ to a proper algebraic space $M_{\ell}$ satisfying the following:
    \begin{enumerate}
        \item
        \label{item: contractions of P l: isomorphism on the open complement}
        $g$ restricts to an isomorphism on the open complement $\Hilb^{\ell}(X) \setminus \Hilb^{\ell}(Y/B)$. 
        \item 
        \label{item: contractions of P l: contracted locus}
        There is a closed immersion $\theta: B \hookrightarrow M_{\ell}$ such that the following diagram commutes
        \[
            \begin{tikzcd}
             \Hilb^{\ell}(Y/B) \ar[d]
              \ar[r] & \Hilb^{\ell}(X) \ar[d, "g"]\\ B \ar[r, "\theta"] & M_{\ell}
            \end{tikzcd}
        \]
    \end{enumerate}
\end{thm}
\begin{proof}

Let $\alpha \subset \Hilb(X)$ denote the open and closed image of $B \to \Hilb(X)$, thought of as a connected component of the Hilbert scheme. We start by proving that $(\alpha, \ell)$ satisfies (LR). Part (1) in (LR) is automatically satisfied by \Cref{remark: condition 1 automatic dimension}, so we just need to check part (2) of (LR).
Let $K \supset k$ be an algebraically closed extension, and fix two points in $\alpha(K)$, which we may think of as points $b_1, b_2 \in B(K)$.
After base-changing to $K$, we may assume without loss of generality that $K=k$ and $Y_{b_1} \cong Y_{b_2} \cong \mathbb{P}^1_k$.
Set $i: Y_{b_2} \hookrightarrow X$.
There is only one rank one torsion-free sheaf on $Y_{b_2}$ with Hilbert polynomial $P_{\cO_{Y_{b_2}}}-\ell$ up to isomorphism, namely, $G= \cO_{\mathbb{P}^1_k}(-\ell)$.
In order to show (LR), we need to check that $\Hom(\cI_{Y_{b_1}}, i_*(\cO_{\mathbb{P}^1_k}(-\ell))) = 0$. 
We distinguish two cases:
\begin{itemize}
    \item
    If $b_1 \neq b_2$, then we get the following sequence of isomorphisms:
    \begin{align*}
        \Hom(\cI_{Y_{b_1}}, i_*(\cO_{\mathbb{P}^1_k}(-\ell))) &= \Hom(i^*(\cI_{Y_{b_1}}), \cO_{\mathbb{P}^1_k}(-\ell)),
        \text{ by adjunction}
        \\&= \Hom(\cO_{\mathbb{P}^1_k}, \cO_{\mathbb{P}^1_k}(-\ell)),
        \text{ because $Y_{b_1}$ and $Y_{b_2}$ are disjoint}
        \\&=0,
        \text{ by the assumption $\ell>0$.}
    \end{align*}
    \item If $b_1 = b_2$, then $i^*(\cI_{Y_{b_1}})$ is the conormal bundle of the immersion $Y_{b_1}\hookrightarrow X$.
    The short exact sequence of normal bundles for $Y_{b_1} \subset Y \subset X$ reads:
    \[ 
    0 \to \cO_{\mathbb{P}^1_k}^{\oplus \dim(B)} \to i^*(\cI_{Y_{b_1}})^\vee \to \bigoplus_{j=1}^{\dim(X)-\dim(Y)} \cO_{\mathbb{P}^1_k}(n_j) \to 0.
    \]
    Applying $\Hom(\blank, \cO_{{\mathbb{P}^1_k}}(-\ell))$ to $\cI_{Y_{b_1}}$ and using the short exact sequence above together with the assumption that $\ell>0$ and $\ell > N \geq n_j$ for all $j$, we get the desired vanishing:
    \begin{align*}
        \Hom(\cI_{Y_{b_1}}, i_*(\cO_{\mathbb{P}^1_k}(-\ell))) &= \Hom(i^*(\cI_{Y_{b_1}}), \cO_{\mathbb{P}^1_k}(-\ell))
        \\&= \mathrm H^0 (Y_{b_1},i^*(\cI_{Y_{b_1}})^\vee (-\ell)) = 0.
    \end{align*}
\end{itemize}

We now turn to proving \ref{item: contractions of P l: isomorphism on the open complement}.
In view of \cref{thm: main theorem on U_alpha} and \cref{prop: locus where the modified Hilbert scheme is an iso}, there is a morphism $f: \Hilb^{\ell}(X) \to \Hilb^{\ell}(X)_{\alpha}$ to a proper algebraic space $\Hilb^{\ell}(X)_{\alpha}$ which is an isomorphism when restricted to the open complement
\[\Hilb^{\ell}(X)^{\rnc\alpha} := \Hilb^{\ell}(X) \setminus \Hilb^{\ell}(\cC/\alpha) = \Hilb^{\ell}(X) \setminus \Hilb^{\ell}(Y/B),\]
cf. \cref{notation: nc-alpha}.
We set $M_{\ell} \subset \Hilb^{\ell}(X)_{\alpha}$ to be the closed scheme-theoretic image of $f$.

To conclude the proof, we are left to find a closed immersion $\theta: B \hookrightarrow M_{\ell}$ such that the diagram in part \ref{item: contractions of P l: contracted locus} is commutative.
Note that there is a natural closed immersion $B \hookrightarrow \alpha$, and that the universal family $\cC \to \alpha$ is a smooth family with geometric fibers isomorphic to $\mathbb{P}^1$. It follows that the relevant compactified Jacobian $\overline{\Pic}(\cC/\alpha)_{\chi(\cO)-\ell}$ defined as in the proof of \Cref{prop: abel jacobi contractions} is actually isomorphic to the base $\alpha$ in this case. Using the inclusion $B \hookrightarrow \alpha = \overline{\Pic}(\cC/\alpha)_{\chi(\cO)-\ell}$, it suffices to find a suitable closed embedding $\theta: \overline{\Pic}(\cC/\alpha)_{\chi(\cO)-\ell} \hookrightarrow \Hilb^{\ell}(X)_{\alpha}$ such that the following diagram commutes
\[
    \begin{tikzcd}[ampersand replacement = \&]
    \Hilb^{\ell}(\cC/\alpha) \ar[d, "AJ_{\ell}"]
    \ar[r, "\psi"] \& \Hilb^{\ell}(X) \ar[d, "f"]\\ \overline{\Pic}(\cC/\alpha)_{\chi(\cO)-\ell} \ar[r, "\theta"] \& \Hilb^{\ell}(X)_{\alpha}.
    \end{tikzcd}
\] 
At this point morphism $\theta$ can be defined exactly as in the proof of \Cref{prop: abel jacobi contractions}, and the same argument shows that it is a closed immersion.
\end{proof}

\begin{coroll}[Contracting $\mathbb{P}^{\ell}$-fibrations with fiberwise negative bundle] \label{cor: contraction fiber bundles negative normal}
        Let $X$ be a smooth projective geometrically connected scheme over $k$ of dimension at least $2$.
        Let $Y \hookrightarrow X$ be a smooth geometrically connected closed subscheme of $X$ which admits a $\mathbb{P}^1$-fibration $\pi: Y \to B$ whose normal bundle is fiberwise bounded by $N=-1$.
        Then, for all positive integers $\ell$, the induced morphism $\Hilb^{\ell}(Y/B) \to \Hilb^{\ell}(X)$ is a closed immersion and there is a surjective morphism $g: \Hilb^{\ell}(X) \to M_{\ell}$ to a proper algebraic space satisfying the following:
    \begin{enumerate}
        \item $g$ restricts to an isomorphism on the open complement $\Hilb^{\ell}(X) \setminus \Hilb^{\ell}(Y/B)$, and 
        \item there is a closed immersion $\theta: B \to M_{\ell}$ such that the following diagram commutes
           \[
\begin{tikzcd}
 \Hilb^{\ell}(Y/B) \ar[d]
  \ar[r] & \Hilb^{\ell}(X) \ar[d, "g"]\\ B \ar[r, "\theta"] & M_{\ell}
\end{tikzcd}
\]
\end{enumerate}
\end{coroll}

\begin{proof}
    By \cref{thm: contraction of projective space bundles}, it suffices to show that the closed immersion $B \hookrightarrow \Hilb(X)$ induced by the family of subcurves $Y \to B$ is open. By \cite[\href{https://stacks.math.columbia.edu/tag/00EH}{00EH}]{stacks-project}, it suffices to show that the conormal sheaf $\cI_B|_B$ of the closed immersion $B \hookrightarrow \Hilb(X)$ is zero. Since $B$ is smooth, we have a short exact sequence
    \[ 0 \to \cI_B|_B \to \Omega^1_{\Hilb(X)/k}|_B \to \Omega^1_{B/k} \to 0\]
    of sheaves in $B$. It suffices then to show that $\Omega^1_{\Hilb(X)/k}|_B \to \Omega^1_{B/k}$ is injective. By \cite[\href{https://stacks.math.columbia.edu/tag/00ME}{00ME}]{stacks-project}, this can be checked fiberwise for every geometric point $b \to B$. After base-changing to the residue field of $b$, we may assume that $b$ is a $k$-point with $Y_b \cong \mathbb{P}^1_k \hookrightarrow X$. After dualizing, we need to show equivalently that the morphism of tangent spaces $T_{B,b} \to T_{\Hilb(X),b}$ is surjective. By assumption, the normal bundle $\cN_{Y_b/X}$ of the inclusion $Y_b \hookrightarrow X$ fits into a short exact sequence 
    \[
        0 \to \cO_{\mathbb{P}^1_k}\otimes_k T_{B,b} \to \cN_{Y_b/X} \to \bigoplus_{j=1}^{\dim(X) - \dim(B)-1} \cO_{\mathbb{P}^1_k}(n_j) \to 0,
    \]
    where $n_j<0$ for all $j$. Recall that the tangent space $T_{\Hilb(X), b}$ is isomorphic to $H^0(\mathbb{P}^1_k, \cN_{Y_b/X})$, and the morphism of tangent spaces $T_{B,b} \to T_{\Hilb(X), b}$ is identified with the induced morphism
    \[ T_{B,b} \cong H^0(\mathbb{P}^1_k, \cO_{\mathbb{P}^1_k}\otimes_k T_{B,b}) \hookrightarrow H^0(\mathbb{P}^1_k, \cN_{Y_b/X}).\]
    The surjectivity of this morphism is then readily seen because $H^0(\mathbb{P}^1_k, \cO_{\mathbb{P}^1_k}(n_j))=0$ for all $j$.
    \qedhere
\end{proof}

\begin{example}[Contractions of $\mathbb{P}^1$-bundles] \label{example: fujiki nakano contractions}
    By setting $\ell=1$ in \cref{cor: contraction fiber bundles negative normal}, it follows that a $\mathbb{P}^1$-bundle $Y \subset X$ with normal bundle fiberwise bounded by $-1$ (\cref{defn: fiberwise bounded normal bundle}) can be contracted in the fiber direction. In the case when $Y$ is a divisor in $X$ (i.e. $\dim(Y) = \dim(X)-1$), this is a known result due to Fujiki-Nakano \cite{fujiki-nakano} and Artin \cite[Cor. 6.11]{artin_algebraization_formal_ii}, who proved more generally a contraction result for divisors $Y \subset X$ that admit $\mathbb{P}^n$-bundle fibrations $Y \to B$.
\end{example}

% -
\subsection{Contractions of Hilbert schemes of non-movable rational curves}
\label{subsec: contractions of rigid P1s}

In this subsection, we explain the special case of \cref{thm: contraction of projective space bundles} when $B =\Spec(k)$.
In this case, we obtain contractions in the following sense.
\begin{defn}[Contractions of subschemes] \label{defn: contraction}
    Let $Y$ be a proper scheme over $k$, and let $Z \subset Y$ be a closed subscheme.
    A \emph{contraction} of $Z$ is a surjective morphism $f: Y \to \overline{Y}$ to a proper algebraic space such that the following are satisfied:
    \begin{enumerate}
        \item The image of the subscheme $Z \subset Y$ via $f$ is a point $p \in \overline{Y}$, and
        \item The restriction $f: Y \setminus Z \to \overline{Y} \setminus p$ is an isomorphism.
    \end{enumerate}
    If such $f$ exists, then we say that $Z \subset Y$ \emph{admits a contraction}.
\end{defn}

\begin{coroll} \label{cor: contraction of hilbert schemes rational curves}
    Let $X$ be a smooth projective geometrically connected scheme over $k$ of dimension at least $2$, and let $C \hookrightarrow X$ be a smooth geometrically connected subcurve of genus $0$. Suppose that $[C \subset X]$ is an isolated point in the Hilbert scheme $\Hilb(X)$ (note that we allow the corresponding component of $\alpha \subset \Hilb(X)$ to be nonreduced, i.e. we allow infinitesimal deformations of $C \subset X$). Let $\cN_{C/X}|_{\mathbb{P}^1_{\overline{k}}} \cong \bigoplus_{j=1}^{\dim(X)-1} \cO_{\mathbb{P}^1_{\overline{k}}}(n_j)$ be the pullback of the normal bundle to the base-change $C_{\overline{k}} \cong \mathbb{P}^1_{\overline{k}}$. Then for all positive integers $\ell > \max_{j} \{ n_j\}$, the morphism $\Hilb^{\ell}(X) \to \Hilb^{\ell}(X)_{\alpha}$ is a contraction of $\Hilb^{\ell}(C) \hookrightarrow \Hilb^{\ell}(X)$ as in \cref{defn: contraction}.
\end{coroll}
\begin{proof}
    This is a special case of \cref{thm: contraction of projective space bundles} when $B \cong \Spec(k)$, where the surjectivity of $\Hilb^{\ell}(X) \to \Hilb^{\ell}(X)_{\alpha}$ is guaranteed by \Cref{prop: condition for surjectivity from Hilbert scheme}.
    \end{proof}

\begin{coroll} 
\label{cor: contraction of P^1 that don't move}
    Let $X$ be a smooth projective geometrically connected scheme of dimension at least $2$, and let $C \hookrightarrow X$ be a smooth geometrically connected subcurve of genus $0$. Suppose that the following hold:
    \begin{enumerate}
        \item[(A)] $[C \subset X]$ is an isolated point in the Hilbert scheme $\Hilb(X)$ (note that we allow the corresponding component of $\alpha \subset \Hilb(X)$ to be nonreduced, i.e. we allow infinitesimal deformations of $C \subset X$), and
        \item[(B)] the pullback of the normal bundle $\cN_{C/X}$ to the base-change $C_{\overline{k}} \cong \mathbb{P}^1_{\overline{k}}$ is isomorphic to $\bigoplus_{j=1}^{\dim(X)-1} \cO_{\mathbb{P}^1_{\overline{k}}}(n_j)$ with $n_j\leq 0$ for all $j$.
    \end{enumerate}
    Then, the morphism $X \to \Hilb^1(X)_{\alpha}$ is a contraction of the curve $C \subset X$ as in \cref{defn: contraction}.
\end{coroll}
\begin{proof}
    This follows from \cref{cor: contraction of hilbert schemes rational curves} for $\ell=1$.
\end{proof}

\begin{example}[Contractions of rational curves in threefolds] \label{example: contraction of rational curves in threefolds}
Let $X$ be a smooth projective threefold, and consider an embedding $i: C \cong \mathbb{P}^1_k \hookrightarrow X$. A result of Reid \cite[Part II]{mmp_threefolds_reid} asserts that if $C$ satisfies both of the following conditions:
\begin{itemize}
    \item we have either $\cN_{C/X} \cong \cO_{\mathbb{P}^1_k}(-1)^{\oplus 2}$ or $\cN_{C/X} \cong \cO_{\mathbb{P}^1_k} \oplus \cO_{\mathbb{P}^1_k}(-2)$, and
    \item $[C \subset X]$ is an isolated point of $\Hilb(X)$,
\end{itemize}
then $C \subset X$ admits a contraction.
We recover this example as a special case of \cref{cor: contraction of P^1 that don't move}.
\end{example}

\begin{example}[Beyond non-positive normal bundle] \label{example: beyond nonpositive normal bundle}
    Let $X$ be a smooth projective threefold, and let $i: C \cong \mathbb{P}^1_k \hookrightarrow X$ be a closed immersion such that $\cN_{C/X} \cong \cO_{\mathbb{P}^1_k}(1) \oplus \cO_{\mathbb{P}^1_k}(-3)$. Assume that $[C \subset X]$ is an isolated point of $\Hilb(X)$. In this case it is not always possible to contract $C \subset X$, see \cite[\S 5]{jimenez_contraction} for an example. However, we may use our results to provide a sufficient condition that guarantees the existence of a contraction.
    
    By \Cref{cor: contraction of hilbert schemes rational curves}, we have that $\mathbb{P}^2_k \cong \Hilb^2(C) \hookrightarrow \Hilb^2(X)$ admits a contraction $\Hilb^2(X) \to \Hilb^2(X)_{\alpha}$.
    Let $\pi: \underline{\mathbb{P}}_X(\Omega^1_X) \to X$ be the projectivization of the cotangent bundle $\Omega^1_X$ parameterizing one-dimensional quotients.
    Then we have a natural closed immersion $\underline{\mathbb{P}}_X(\Omega^1_X) \hookrightarrow \Hilb^2(X)$ which sends a pair $(x, \Omega^1_X|_x \twoheadrightarrow L)$ to the unique length two subscheme of $X$ that is set-theoretically supported on $x$ and has tangent space $L^{\vee} \hookrightarrow T_X|_x$ at $x$.
    Note that the intersection $Z := \Hilb^2(C) \cap \underline{\mathbb{P}}_X(\Omega^1_X) \subset \Hilb^2(X)$ is the image of the section $\sigma: C \to \underline{\mathbb{P}}_X(\Omega^1_X)$ defined by the pair $(C \to X, \Omega^1_X|_C \twoheadrightarrow \Omega^1_C)$.
    Hence $Z$ is naturally identified with $C \subset X$ under the induced morphism $\pi: Z \to X$.
    If we replace $\Hilb^2(X)_{\alpha}$ by the scheme-theoretic image $M_2$ of $\underline{\mathbb{P}}_X(\Omega^1_X) \hookrightarrow \Hilb^2(X) \to \Hilb^2(X)_{\alpha}$, then it follows that the induced morphism $\underline{\mathbb{P}}_X(\Omega^1_X) \to M_2$ is a contraction of $C \cong Z \subset \underline{\mathbb{P}}_X(\Omega^1_X)$. 
    
    We now note that a sufficient condition for $C \subset X$ to admit a contraction is that the section $\sigma: C \to \underline{\mathbb{P}}_X(\Omega^1_X)$ extends to a global section $\widetilde{\sigma}: X \to \underline{\mathbb{P}}_X(\Omega^1_X)$. In other words, a contraction exists if there exists a line bundle quotient $\Omega_X^1 \twoheadrightarrow \widetilde{L}$ that recovers the quotient $\Omega^1_X|_C \twoheadrightarrow \Omega^1_C$ when restricted to $C$. Since the existence of a contraction can be checked formally locally around $C$ \cite[Thm. 3.1]{artin_algebraization_formal_ii}, it is even sufficient to show the existence of such a quotient on the formal completion of $C \subset X$. We note that in this case one may use the line bundle $\widetilde{L}$ to verify the criterion in \cite[Thm. 1]{jimenez_contraction}, so one may view this example as a special case of the main theorem in \cite{jimenez_contraction}.
\end{example}

% ---
\section{Surgery diagrams} \label{section: surgery diagrams}

This section is dedicated to setting up the non-GIT wall-crossing diagram as in \cref{fig: surgery diagram}, proving the second main theorem about existence of good moduli spaces on the other side of the wall (\cref{thm: existence of good moduli space other side of the wall}), and studying how it can be applied to the examples from \cref{section: applications}.

We start by defining the substack $\cU^\circ_{\alpha,\ell}$ of $\cU_{\alpha,\ell}$, proving that it is open in $\cU_{\alpha,\ell}$, and verifying that it admits a proper good moduli space.
We then proceed to explain the corresponding surgery diagrams for the examples of contractions from \cref{section: applications}.
Specializing the surgery of $\P^\ell$-fibrations allows us to construct surgeries of non-movable rational curves, and recover Fujiki-Nakano divisorial contractions and flops of rational curves in threefolds.
We end the section with an application to the DK hypothesis by employing some of the methods from \cite{bridgeland_flops}: under certain conditions, we get equivalences of derived categories for flops of $\P^1$-fibrations (\cref{thm: derived categories of surgeries}).

% -
\subsection{The other side of the wall}

\begin{defn}
\label{defn: other side of the wall}
    Let $\alpha \subset \Hilb(X)$ be as in \Cref{setup: alpha}. We denote by $\nps_{\alpha, \ell}^{\circ} \subset \nps_{\alpha, \ell}$ the substack that sends a $k$-scheme $T$ to the subgroupoid of $T$-families $T \to \nps_{\alpha, \ell}$ of sheaves $F$ on $X_T \to T$ such that for every $t \in T$ the fiber $F_t$ is either:
    \begin{enumerate}
        \item of the ideal type of the form $\cI_Z$ for some length $\ell$ subscheme $[Z \subset X_t] \in \Hilb^{\ell}(X)^{\rnc\alpha}(t)$ (see \Cref{notation: nc-alpha}), or
        \item of the nonsplit torsion type (as in \Cref{remark: split and nonsplit torsion type}).
    \end{enumerate}
\end{defn}

The main theorem of this subsection is the following.
\begin{thm} \label{thm: existence of good moduli space other side of the wall}
    Let $\alpha \subset \Hilb(X)$ be as in \Cref{setup: alpha}. Suppose that $\alpha$ is of finite type, and that $(\alpha, \ell)$ satisfies (LR). Then, the following hold:
    \begin{enumerate}
        \item $\nps_{\alpha, \ell}^{\circ}$ is an open substack of $\nps_{\alpha, \ell}$;
        \item $\nps_{\alpha, \ell}^{\circ}$ is a $\mathbb{G}_m$-gerbe over a proper algebraic space $\Hilb^{\ell}(X)^{\circ}_{\alpha}$. In particular, the induced morphism $\nps_{\alpha, \ell}^{\circ} \to \Hilb^{\ell}(X)^{\circ}_{\alpha}$ is a good moduli space; and
        \item the $\mathbb{G}_m$-gerbe $\nps_{\alpha, \ell}^{\circ} \to \Hilb^{\ell}(X)^{\circ}_{\alpha}$ admits a section. In other words, $\Hilb^{\ell}(X)_{\alpha}^{\circ}$ is a fine moduli space of sheaves.
    \end{enumerate}
\end{thm}

\begin{proof}
    We begin by proving (1). Consider the stack of filtrations $\Filt(\nps_{\alpha, \ell}) := \Map(\Theta_k, \nps_{\alpha, \ell})$ as in \cite[Defn. 1.1.11]{halpernleistner2018structure}.
    By the Rees construction \cite[Prop. 1.0.1]{halpernleistner2018structure}, we may describe the stack $\Filt(\nps_{\alpha, \ell})$ concretely as the algebraic stack that sends a $k$-scheme $T$ to the groupoid of $T$-families of sheaves $F \in \nps_{\alpha, \ell}(T)$ equipped with a $\mathbb{Z}$-indexed filtration $(F_m)_{m \in \mathbb{Z}}$ by subsheaves satisfying the following:
    \begin{enumerate}[label=(\alph*)]
        \item
        $F_{m+1} \subset F_{m}$ for all $m$. Furthermore, we have $F_m = 0$ and $F_{-m} = F$ for all sufficiently large $m \gg0$.
        \item
        The associated graded sheaf $\bigoplus_{m \in \mathbb{Z}} F_m/ F_{m+1}$ is a $T$-point of $\nps_{\alpha, \ell}$.
    \end{enumerate}
    For each of the finitely many components $\alpha_j \subset \alpha$, let $P_j$ denote the common Hilbert polynomial of the ideal sheaves of subschemes in $\alpha_j$.
    Inside of $\Filt(\nps_{\alpha, \ell})$, there is an open and closed substack $\mathscr{W} \subset \Filt(\nps_{\alpha, \ell})$ parameterizing $T$-families of sheaves $F$ equipped with a two-step $\mathbb{Z}$-weighted filtration
    \[\ldots 0 \subset 0  \subset F_{1} \subset F_0 = F \subset F \subset F \ldots \]
    such that the Hilbert polynomial of $F_{1}$ is equal to one of the finitely many Hilbert polynomials $P_j$ of the ideal sheaves of subschemes in $\alpha$.
    The requirement that the associated graded lies in $\nps_{\alpha, \ell}$ actually forces every fiber of $F_{1}$ to be an ideal sheaf $\cI_{W'}$ of a subscheme in $\alpha$.
    Note that $\mathscr{W}$ can be viewed as an open substack of a union of relative Quot stacks
    \[\bigsqcup_j \Quot_{X\times \nps_{\alpha, \ell}/ \nps_{\alpha, \ell}}(F_{\univ})_{P_{\ell} - P_j}\]
    for the universal family of sheaves $F_{\univ}$ on the smooth projective morphism $X\times \nps_{\alpha, \ell} \to \nps_{\alpha, \ell}$.
    Since we are considering finitely many Hilbert polynomials, it follows that $\mathscr{W}$ is of finite type. 
    
    Consider the representable forgetful morphism $\ev_1: \mathscr{W} \to \nps_{\alpha, \ell}$ given by forgetting the filtration (i.e. it sends a filtered sheaf $(F, (F_m)_{m \in \mathbb{Z}})$ to the underlying sheaf $F$).
    We claim that $\ev_1$ is a proper morphism.
    Since the formation of $\Filt(\nps_{\alpha, \ell})$ and the forgetful morphism $\ev_1$ commute with base change on the moduli space $\Hilb^{\ell}(X)_{\alpha}$ \cite[Lem. 1.1.1]{halpernleistner2018structure}, we may check this \'etale locally on $\Hilb^{\ell}(X)_{\alpha}$.
    By the version of the Luna \'etale slice theorem in \cite[Thm. 6.1]{alper-hall-rydh-theetalelocalstructureofalgebraicstacks},
    after \'etale base-change on $\Hilb^{\ell}(X)_{\alpha}$, we may replace the stack $\nps_{\alpha, \ell}$ with a quotient stack of the form $[\Spec(A)/G]$, where $A$ is a finitely generated $k$-algebra and $G$ is a linearly reductive group.
    Hence, we are reduced to proving that for every open and closed substack $\mathscr{V} \subset \Filt([\Spec(A)/G])$ of finite type inside the stack of filtrations, the evaluation morphism $\ev_1: \mathscr{V} \to [\Spec(A)/G]$ is proper.
    This is standard, and follows, for example, from the explicit description of $\Filt([\Spec(A)/G])$ in \cite[Thm. 1.4.7]{halpernleistner2018structure}.

    Let $\ev_1(\mathscr{W}) \subset |\nps_{\alpha, \ell}|$ denote the closed image of the proper morphism $\ev_1: \mathscr{W} \to \nps_{\alpha, \ell}$ inside the topological space of $\nps_{\alpha, \ell}$.
    Recall from our discussion above that the filtrations $F_1 \subset F_0 = F$ appearing in $\mathcal{W}$ are exactly the ones where $F_1$ is isomorphic to the ideal sheaf of a subscheme in $\alpha$.
    Using property (LR), this implies that the open complement $|\nps_{\alpha, \ell}| \setminus \ev_1(\mathscr{W})$ consists exactly of the geometric points $F \in \nps_{\alpha, \ell}$ satisfying conditions (1) or (2) in \Cref{defn: other side of the wall}.
    Hence $\nps_{\alpha, \ell}^{\circ}$ is an open substack of $\nps_{\alpha, \ell}$, thus concluding the proof of part (1).

    For part (2), we note that every family of sheaves contains a copy of $\mathbb{G}_m$ inside its group of automorphisms, given by the scalar automorphisms of the sheaf.
    Equivalently, we may think of this as an inclusion of the relative group scheme $\mathbb{G}_m \times \nps_{\alpha, \ell}^{\circ} \to \nps_{\alpha, \ell}^{\circ}$ into the relatively affine inertia $I_{\nps_{\alpha, \ell}^{\circ}} \to \nps_{\alpha, \ell}^{\circ}$.
    Note that for every algebraically closed field $K \supset k$ and every point $F \in \nps_{\alpha, \ell}^{\circ}(K)$, the sheaf $F$ is simple, and hence its automorphism group scheme is $\mathbb{G}_m$.
    Since the inclusion $\mathbb{G}_m \times \nps_{\alpha, \ell}^{\circ} \hookrightarrow I_{\nps_{\alpha, \ell}^{\circ}}$ is an isomorphism fiberwise,
    it is an isomorphism (by a similar argument as in \cite[Lem. 5.1]{de_Cataldo_herrero_log_hodge}), and then the rigidification $\Hilb^{\ell}(X)_{\alpha}^{\circ} := \nps_{\alpha, \ell}^{\circ}\!\fatslash \, \mathbb{G}_m$ as in \cite[Appendix A]{abramovich2007tame} is an algebraic space.
    Hence we have produced a $\mathbb{G}_m$-gerbe $\nps_{\alpha, \ell}^{\circ} \to \Hilb^{\ell}(X)_{\alpha}^{\circ}$, which is automatically a good moduli space. 
    
    Next, let us show that $\Hilb^{\ell}(X)_{\alpha}^{\circ}$ is separated.
    For this, it suffices to show that $\nps_{\alpha, \ell}^{\circ}$ is $S$-complete as in \Cref{defn: S and axis completeness}.
    Since $\nps_{\alpha, \ell}^{\circ}$ is an open substack of the $S$-complete stack $\nps_{\alpha, \ell}$ (\Cref{prop: S completeness}), then by \Cref{lemma: axis implies S}, we just need to show that the open immersion $\nps_{\alpha, \ell}^{\circ} \hookrightarrow \nps_{\alpha, \ell}$ is axis complete.
    For this, we use a notation and setup similar to that in the proof of \Cref{prop: S completeness}.
    Let $K \supset k$ be a field extension, and fix a morphism $\cX_K \to \nps_{\alpha, \ell}$ such that the image of $\cX_K \setminus \mathfrak{o}$ lies in the open substack $\nps_{\alpha, \ell}^{\circ}$. This corresponds to a pair of sheaves $F,F' \in \nps_{\alpha, \ell}^{\circ}(K)$ equipped with filtrations $\Fil_{\bullet}(F), \Fil^{\bullet}(F')$ as in the proof of \Cref{prop: S completeness} such that $\gr_i(F) \cong \gr^i(F')$. We need to show then that the associated graded sheaf $\gr_{\bullet}(F) \cong \gr^{\bullet}(F') \in \nps_{\alpha, \ell}(K)$ lies in the open substack $\nps_{\alpha, \ell}^{\circ}$.

    The case when the filtrations are shifted trivial is tautological, just as in \cref{prop: S completeness}.
    We will prove that this is the only possibility.
    For the sake of contradiction, assume that the filtrations are not shifted trivial.
    We will analyze the same cases as in \cref{prop: S completeness}.
    \begin{itemize}
        \item Case: $F'$ is of the ideal type, i.e. $F' = \cI_Z$ for some length $\ell$ subscheme of $X_K$.
        The argument in this case of \cref{prop: S completeness} implies that $F$ is of torsion type, and gives an embedding of ideal sheaves $\cI_W \subset \cI_Z$ for some $W \in \alpha(K)$.
        But this then implies that $Z \subset W$ scheme-theoretically, which contradicts \cref{defn: other side of the wall}.
        \item Case: both $F$ and $F'$ are of torsion type.
        The argument in this case of \cref{prop: S completeness} implies that either $F$ or $F'$ is of split torsion type, which again contradicts \cref{defn: other side of the wall}.
    \end{itemize}
    This concludes the proof of separatedness.
    
    We are left to show that $\Hilb^{\ell}(X)_{\alpha}^{\circ}$ is proper.
    Since the surjective morphism $\nps_{\alpha, \ell}^{\circ} \to \Hilb^{\ell}(X)_{\alpha}^{\circ}$ is universally closed \cite[Thm. 4.16(ii)]{alper-good-moduli}, it suffices to prove that $\nps_{\alpha, \ell}^{\circ}$ satisfies the existence part of the valuative criterion for properness.
    Recall that $\nps_{\alpha, \ell}$ satisfies the existence part of the valuative criterion by \Cref{prop: existence valuative criterion}.
    By the semistable reduction theorem for $\Theta$-strata \cite[Thm. 6.3]{AHLH}, in order to show the desired statement for $\nps_{\alpha, \ell}^{\circ} \subset \nps_{\alpha, \ell}$ it is enough to show that the locus $\ev_1: \cW \to \nps_{\alpha, \ell}$ that we are removing is a $\Theta$-stratum.
    We have shown that $\ev_1$ is proper.
    
    We next show that $\ev_1$ is radicial. This amounts to showing that for all algebraically closed fields $K \supset k$ and all sheaves $F \in \nps_{\alpha, \ell}(K)$, there is at most one subsheaf of $F$ that is isomorphic to an ideal sheaf of the form $\cI_{Y}$ for some $[Y \subset X_K] \in \alpha(K)$.
    After base-change, we may replace the ground field with $K$ and assume that $F$ is a $k$-point without loss of generality.
    If $F$ is of torsion type, then the desired uniqueness of such ideal subsheaf follows from \Cref{lemma: Hom from I W' to non-split F is zero}.
    On the other hand, if $F$ is of ideal type, then it is of the form $F \cong \cI_Z$ for some length $\ell$ subscheme $Z \subset X$.
    An inclusion of an ideal sheaf $\cI_{Y} \hookrightarrow \cI_Z$ then corresponds to a subscheme $[i: Y \hookrightarrow X]$ such that $Y \supset Z$.
    Note that in that case $\cI_Z / \cI_{Y} \cong i_*(\cI_{Z \subset Y})$.
    If $Z$ was contained in two distinct such subschemes $Y, Y'$, then the sheaf $\cI_{Z}$ would degenerate via the Rees construction to two non-isomorphic sheaves $\cI_{Y} \oplus i_*(\cI_{Z \subset Y})$ and $\cI_{Y'} \oplus (i')_*(\cI_{Z \subset Y'})$ of split torsion type in $\nps_{\alpha, \ell}(k)$.
    Since these sheaves are closed in $\nps_{\alpha, \ell}$ (\Cref{prop: closed points}), this cannot happen by the property that every $k$-point of the stack $\nps_{\alpha, \ell}$ with good moduli space specializes to a unique point \cite[Thm. 4.16(iv) + Prop. 9.1]{alper-good-moduli}.
    This concludes the proof that $\ev_1 : \mathscr{W} \to \nps_{\alpha, \ell}$ is radicial. 

    Since $\ev_1$ is proper, we have that $\ev_1$ is finite radicial, and so $\ev_1: \cW \to \nps_{\alpha, \ell}$ is a weak $\Theta$-stratum in the sense of \cite[Defn. 2.1.1]{halpernleistner2018structure}. By \cite[Lem. 2.1.7]{halpernleistner2018structure}, in order to show that $\ev_1: \cW \to \nps_{\alpha, \ell}$ is a $\Theta$-stratum it is enough to show that for all field-valued points $w: \Spec(K) \to \cW$ we have $H^0(w^*(\mathbb{L}_{\cW/\nps_{\alpha, \ell}})) = H^1(w^*(\mathbb{L}_{\cW/\nps_{\alpha, \ell}}))=0$ for the fiber of the cotangent complex of $\ev_1$. Choose such $w$ corresponding to a filtration $0 \subset F_1 \subset F$ of sheaves on $X_K$. Since $\ev_1$ is representable, it is automatic that $H^1(w^*(\mathbb{L}_{\cW/\nps_{\alpha, \ell}}))=0$. On the other hand, we have that $H^0(w^*(\mathbb{L}_{\cW/\nps_{\alpha, \ell}}))$ is the dual of relative tangent space of $\ev_1: \cW \to \nps_{\alpha, \ell}$ at $w$. In view of our description of $\cW$ in part (1), this agrees with the dual of the tangent space of the corresponding Quot scheme, which is given by $\Hom(F_1, F/F_1)^{\vee}$. We show that this is zero by considering two separate cases.
    \begin{itemize}
        \item Suppose that $F \cong \cI_Z$ is of ideal type. Then, the inclusion $F_1 \cong \cI_Y \subset \cI_Z$ corresponds to a point $[i: Y \hookrightarrow X_K] \in \alpha(K)$ such that $Z \subset Y$. The quotient $F/F_1$ is isomorphic to $i_*(G)$, where $G$ is the ideal sheaf of $Z$ in $Y$. Since $(\alpha, \ell)$ satisfies (LR), it follows that $\Hom(F_1, F/F_1)^{\vee}=0$.
        \item Suppose that $F \cong \cI_W \oplus i_*(G)$ is of split torsion type. Then, we have that $F_1 \cong \cI_W \subset \cI_W \oplus i_*(G)$ is the inclusion of the first summand. Again, we see that $\Hom(F_1, F/F_1)^{\vee}= \Hom(\cI_W, i_*(G))^{\vee} =0$.
    \end{itemize}
    This concludes the proof of (2).
    
    We are left to show part (3).
    In order to show that the $\mathbb{G}_m$-gerbe $\nps_{\alpha, \ell}^{\circ} \to \Hilb^{\ell}(X)^{\circ}_{\alpha}$ splits, it suffices to find a line bundle $L$ on the stack $\nps_{\alpha, \ell}^{\circ}$ such that the group of scalar automorphisms $\mathbb{G}_m$ acts on the fibers of $L$ with weight $1$.
    Let $F_{\univ}$ denote the universal sheaf on the source of the smooth projective morphism $\pi: X_{\nps_{\alpha, \ell}^{\circ}} \to \nps_{\alpha, \ell}^{\circ}$, and let $\det(F_{\univ})$ be its determinant line bundle.
    Since $\nps_{\alpha, \ell}^{\circ}$ is contained in the closed substack $\SCoh(X)_{\cO_X} \hookrightarrow \SCoh(X)$ of coherent sheaves with fiberwise trivial determinant (see \Cref{subsection: moduli of coherent sheaves}), it follows that $\pi_*(\det(F_{\univ}))$ is a line bundle on $\nps_{\alpha, \ell}^{\circ}$ and the counit $\pi^*\pi_*(\det(F_{\univ})) \to \det(F_{\univ})$ is an isomorphism.
    Since the sheaf $F_{\univ}$ has generic rank $1$ on the fibers of $\pi$, the scaling group of automorphisms $\mathbb{G}_m$ acts with weight $1$ on the fibers of the line bundle $\det(F_{\univ})$.
    Hence the scaling group $\mathbb{G}_m$ acts with weight $1$ on the line bundle $\pi_*(\det(F_{\univ}))$, as required.
\end{proof}

% -
\subsection{Examples of surgery diagrams}

In the setting of \Cref{thm: existence of good moduli space other side of the wall}, the open immersions $\SCoh^{\rtf}(X)_{\cO_X, \ell} \hookrightarrow \nps_{\alpha, \ell}$ and $\nps_{\alpha, \ell}^{\circ} \hookrightarrow \nps_{\alpha, \ell}$ induce a diagram of proper good moduli spaces
    \begin{equation} \label{equation: general surgery diagram}
\begin{tikzcd}[ampersand replacement = \&]
 \Hilb^{\ell}(X) \ar[dr, "p_1"']
   \& \& \Hilb^{\ell}(X)_{\alpha}^{\circ} \ar[dl, "p_2"]\\ \& \Hilb^{\ell}(X)_{\alpha} \&
   \end{tikzcd}
\end{equation}
where all morphisms are isomorphisms over the common open locus $\Hilb^{\ell}(X)^{\rnc\alpha}$.
One may think of this as a surgery diagram, where we remove the fibers contracted by $\Hilb^{\ell}(X) \to \Hilb^{\ell}(X)_{\alpha}$ and glue other fibers to obtain $\Hilb^{\ell}(X)_{\alpha}^{\circ}$. 
In this subsection, we provide several examples of this following the geometric setups studied in \Cref{section: applications}. We note that in some cases we are performing, in the language of algebraic geometry, a flip or flop of the corresponding locus in $\Hilb^{\ell}(X)$ that is contracted.

\begin{example}[Abel-Jacobi surgeries] \label{example: abel jacobi surgeries}
    Let us place ourselves in the situation of \Cref{prop: abel jacobi contractions}. Consider the induced diagram of proper algebraic spaces as in \eqref{equation: general surgery diagram}.
    Let us denote by $Y_{\ell}$ the scheme-theoretic image of $p_1$, as in the proof of \Cref{prop: abel jacobi contractions}. If we set $Y_{\ell}^{\circ} := (p_2)^{-1}(Y_{\ell})$, then we obtain a diagram
    \[
    \begin{tikzcd}
     \Hilb^{\ell}(X) \ar[dr, "p_1"']
       & & Y_{\ell}^{\circ} \ar[dl, "p_2"]\\ & Y_{\ell} &
    \end{tikzcd}
    \]
    where both morphisms are surjective. Consider the closed immersion $\theta: H_{\ell} \hookrightarrow Y_{\ell}$ as in \Cref{prop: abel jacobi contractions}.
    Then both $p_1$ and $p_2$ are isomorphisms over the open complement of $H_{\ell}$. 
    
    A geometric point $h \in H_{\ell}(K)$ corresponds to a sheaf of split torsion type of the form $i_*(G) \oplus \cI_{Y}$ for some integral curve $i: Y \hookrightarrow X_K \in \alpha(K)$.
    The preimage $(p_1)^{-1}(h)$ corresponds to all of the length $\ell$ subschemes $Z \subset Y$ with ideal sheaf $G$.
    This is just the projectivization $\mathbb{P}(\Hom(G, \cO_Y))$ of the vector space $\Hom(G, \cO_Y)$, which is the corresponding fiber of $[G] \in \overline{\Pic}(\cC/\alpha)_{\chi(\cO)-\ell}(K)$ under the Abel-Jacobi morphism $AJ_{\ell} : \Hilb^{\ell}(\cC/\alpha) \to \overline{\Pic}(\cC/\alpha)_{\chi(\cO)-\ell}$.
    On the other side, the fiber $(p_2)^{-1}(h)$ parameterizes all nontrivial extensions
    \[ 0 \to i_*(G) \to F \to  \cI_{Y} \to 0.\]
    Hence, we have that $(p_2)^{-1}(h) = \mathbb{P}(\Ext^1(\cI_Y, i_*(G)))$. In other words, the surgery diagram above replaces the projective space $\mathbb{P}(\Hom(G, \cO_Y))$ with $\mathbb{P}(\Ext^1(\cI_Y, i_*(G)))$ for every $h \in H_{\ell} \subset Y_{\ell}$.
\end{example}

\begin{example}[Mukai flop] \label{example: mukai flop surgery} We can consider the Abel-Jacobi surgery from \Cref{example: abel jacobi surgeries} in the special example of the Mukai flop discussed in \Cref{example: mukai flop contractions}. In this case we recover the other side of the flop as described in \cite[Ex. 1.7(ii)]{namikawa_deformation_symplectic}.
\end{example}

\begin{example}[Surgery of $\mathbb{P}^{\ell}$-fibrations] \label{example: surgery of P ell fibrations}
    Let us now place ourselves in the situation of \Cref{thm: contraction of projective space bundles}. Consider the surgery diagram as in \eqref{equation: general surgery diagram}, and let $M_{\ell} \subset \Hilb^{\ell}(X)_{\alpha}$ be the scheme-theoretic image of $p_1$ as in the proof of \Cref{thm: contraction of projective space bundles}. Setting $M_{\ell}^{\circ} := (p_2)^{-1}(M_{\ell})$, we obtain a diagram of surjective morphisms 
        \[
\begin{tikzcd}
 \Hilb^{\ell}(X) \ar[dr, "p_1"']
   & & M_{\ell}^{\circ} \ar[dl, "p_2"]\\ & M_{\ell} &
\end{tikzcd}
\]
Consider the closed immersion $\theta: B \hookrightarrow M_{\ell}$
corresponding to the contracted locus.
A geometric point $b \in B(K)$ corresponds to a sheaf of split torsion type of the form $i_*(\cO_{\mathbb{P}^1_K}(-\ell)) \oplus \cI_{Y_b}$, where we denote by $i$ the inclusion $i: \mathbb{P}^1_K \cong Y_b \hookrightarrow X_K$. 
Similarly as in \Cref{example: abel jacobi surgeries}, we have that $(p_1)^{-1}(b)$ is the fiber of the corresponding $\mathbb{P}^{\ell}$-fibration $\Hilb^{\ell}(Y/B) \to B$, which is isomorphic to the projectivization $\mathbb{P}(H^0(\mathbb{P}^1_K, \cO_{\mathbb{P}^1_K}(\ell))$ of the vector space $\Hom(\cO(-\ell), \cO) = H^0(\mathbb{P}^1_K, \cO_{\mathbb{P}^1_K}(\ell)))$.
On the other hand, a computation shows that the fiber $(p_2)^{-1}(b)$ is the projectivization $\mathbb{P}(H^1(\mathbb{P}^1_K, \cN_{Y_b/X_K}(-\ell))$ of the vector space $\Ext^1(\cI_{Y_b}, i_*(\cO_{\mathbb{P}^1_K}(-\ell))) \cong H^1(\mathbb{P}^1_K, \cN_{Y_b/X_K}(-\ell))$, where $\cN_{Y_b/X_K}$ denotes the normal bundle of the inclusion $i: Y_b \hookrightarrow X_K$.   
\end{example}

\begin{example}[Fujiki-Nakano divisorial contractions as moduli of sheaves] \label{example: fujiki nakano surgery}
Let us specialize \Cref{example: surgery of P ell fibrations} by assuming that $Y \hookrightarrow X$ is a divisor, and that for every geometric point $b \in B(K)$, the restriction of the normal line bundle $\cN_{Y/X}$ to $Y_b \cong \mathbb{P}^1_K$ is isomorphic to $\cO_{\mathbb{P}^1_K}(-1)$. Then, in the diagram
% https://q.uiver.app/#q=WzAsNCxbMCwwLCJYIl0sWzQsMCwiTV8xXlxcY2lyYyJdLFsyLDAsIlxcd2lkZXRpbGRlIE1fMSJdLFsyLDIsIk1fMSJdLFswLDMsInBfMSIsMl0sWzEsMywicF8yIl0sWzAsMiwiIiwwLHsic3R5bGUiOnsiYm9keSI6eyJuYW1lIjoiZGFzaGVkIn19fV0sWzIsMSwiIiwwLHsic3R5bGUiOnsiYm9keSI6eyJuYW1lIjoiZGFzaGVkIn19fV0sWzIsMywiIiwxLHsic3R5bGUiOnsiYm9keSI6eyJuYW1lIjoiZGFzaGVkIn19fV1d
\[\begin{tikzcd}[ampersand replacement=\&]
	X \&\& {\widetilde M_1} \&\& {M_1^\circ} \\
	\\
	\&\& {M_1}
	\arrow[dashed, from=1-1, to=1-3]
	\arrow["{p_1}"', from=1-1, to=3-3]
	\arrow[dashed, from=1-3, to=1-5]
	\arrow[dashed, from=1-3, to=3-3]
	\arrow["{p_2}", from=1-5, to=3-3]
\end{tikzcd}\]
the fibers of the birational morphism $p_2$ are points.
Hence, the morphism $p_2$ is an isomorphism after possibly passing to normalizations, and we write $\widetilde{M}_{1}$ for the common normalization of $M_1$ and ${M}_{1}^{\circ}$.
By the universal property of normalizations, there is an induced morphism $X \to \widetilde{M}_1$ making the whole diagram commute, and exhibiting $\widetilde{M}_{1}$ as the smooth Fujiki-Nakano contraction (cf. \Cref{example: fujiki nakano contractions}).
Therefore, we have realized the contraction as a fine moduli of sheaves, possibly up to normalization.
\end{example}

\begin{example}[Surgery of non-movable rational curves]
\label{example: surgery of non-movable rational curves}
    We specialize \Cref{example: surgery of P ell fibrations}
    and put ourselves in a situation as in \Cref{cor: contraction of P^1 that don't move}. We assume that we have an embedding $i: C \cong \mathbb{P}^1_k \hookrightarrow X$ such that the normal bundle $\cN_{C/X}$ is non-positive and the curve $C \subset X$ does not move. If we denote by $\alpha \subset \Hilb(X)$ the connected component of the Hilbert scheme whose unique point is $[C \subset X]$, then we obtain a diagram of surjective morphisms of proper algebraic spaces
     \[
     \begin{tikzcd}
         X \ar[dr, "p_1"']
   & & \Hilb^{1}(X)_{\alpha}^{\circ} \ar[dl, "p_2"]\\ & \Hilb^{1}(X)_{\alpha} 
     \end{tikzcd}
     \]
   The curve $\mathbb{P}^1_k \cong C \subset X$ is contracted to a point $h \in \Hilb^{1}(X)_{\alpha}$. Both morphisms $p_1$ and $p_2$ are isomorphisms over the complement of $h$, and the preimage $(p_2)^{-1}(h)$ is isomorphic to the projective space $\mathbb{P}(H^1(\mathbb{P}^1_k, \cN_{C/X}(-1)))$.
\end{example}

\begin{example}[Flop of a rational curve in a threefold] \label{example: surgery for rational curves in threefolds}
    As a further special case of \Cref{example: surgery of non-movable rational curves}, we can assume that $X$ is a threefold containing a rational curve $C= \mathbb{P}^1_k \hookrightarrow X$ with non-positive normal bundle that does not move. Then, the normal bundle is either $\cN_{C/X} \cong \cO_{\mathbb{P}^1_k}(-1) \oplus \cO_{\mathbb{P}^1_k}(-1)$ or $\cN_{C/X} \cong \cO_{\mathbb{P}^1_k} \oplus \cO_{\mathbb{P}^1_k}(-2)$, and we are in the setting of \cite[Part II]{mmp_threefolds_reid} as explained in \Cref{example: contraction of rational curves in threefolds}. In this case $\Hilb^1(X)_{\alpha}^{\circ}$ recovers the flop of the curve $C \subset X$, which is constructed locally in \cite[Part II]{mmp_threefolds_reid}.
\end{example}

% -
\subsection{Derived categories of \texorpdfstring{$\mathbb{P}^1$}{P1}-fibration surgeries} \label{subsection: P^1-bundle surgeries and derived categories}

\begin{setup}
\label{setup for derived categories}
As usual, we let $X$ be a smooth, projective, geometrically connected scheme over $k$.
We place ourselves in the situation of \Cref{thm: contraction of projective space bundles}, and assume that we have a smooth closed subscheme $Y \hookrightarrow X$ equipped with a $\mathbb{P}^1$-fibration $Y \to B$ such that:
\begin{itemize}
    \item the normal bundle $\cN_{Y/X}$ is fiberwise bounded by $N=0$ in the sense of \Cref{defn: fiberwise bounded normal bundle} (equivalently, the normal bundle over every fiber $Y_b$ is non-positive), and
    \item the induced morphism $B \to \Hilb(X)$ has open image.
\end{itemize} 
\end{setup}

Then we obtain a diagram
    \begin{equation} \label{equation: surgery diagram for P^1 fibration}
\begin{tikzcd}[ampersand replacement = \&]
 X \ar[dr, "p_1"']
   \& \& M_1^{\circ} \ar[dl, "p_2"]\\ \& M_1
   \end{tikzcd}
\end{equation}
of surjective morphisms between proper algebraic spaces, as explained in \Cref{example: surgery of P ell fibrations}.

Our next result establishes a relation between the bounded derived categories of coherent sheaves $D^b(M_1^{\circ})$ and $D^b(X)$ under the assumption that $M_1^{\circ}$ is smooth and the fibers of $Y \to B$ have non-positive intersection with the canonical divisor of $X$.
By Kawamata's DK-hypothesis in \cite[Conj. 1.2]{kawamata_dk_hypothesis}, one expects that there should be a fully faithful admissible embedding $D^b(M_1^{\circ}) \hookrightarrow D^b(X)$.
Note that, by adjunction and our assumption on the non-positivity of normal bundles, for all $b \in B$ we have $\deg(\omega_{X}|_{Y_b}) \geq -2$. We first discuss the two most degenerate cases corresponding to $\deg(\omega_{X}|_{Y_b}) = -2$ and $\deg(\omega_X|_{Y_b}) = -1$.

\begin{example}[DK-hypothesis for $\mathbb{P}^1$-fibrations] \label{example: DK-hypothesis for P^1 fibrations}
    Let $X$ and $Y \subset X$ be as in \cref{setup for derived categories}, and suppose that $M_1^{\circ}$ is smooth.
    Assume that there exists some $b \in B$ such that $\deg(\omega_X|_{Y_b})=-2$.
    This implies by the adjunction formula and the non-positivity of $\cN_{Y_b/X_b}$ that the normal bundle $\cN_{Y_b/X_b}$ is a trivial vector bundle.
    Standard deformation theory then implies that $\Hilb(X)$ is smooth of dimension $\dim(X)-1$ at $[Y_b \subset X_b] \in \Hilb(X)$.
    This forces $\dim(B) = \dim(X)-1$ and $\dim(Y) = \dim(X)$, and hence $Y=X$.
    Hence, in this case $X \to B$ is itself a $\mathbb{P}^1$-fibration, and we have that $B \hookrightarrow M_1$ is an isomorphism.
    It follows from the discussion in \Cref{example: surgery of P ell fibrations} that in this case the fibers of $p_2: M_1^{\circ} \to M_1 = B$ are all empty, and therefore $M_1^{\circ}$ is the empty scheme.
    We conclude that in this case $D^b(M_1^{\circ})$ tautologically embeds fully faithfully as an admissible subcategory of $D^b(X)$.
\end{example}

\begin{example}[DK-hypothesis for divisorial contractions] \label{example: DK-hypothesis divisorial contractions}
    Let $X$ and $Y \subset X$ be as in \cref{setup for derived categories}, and suppose that $M_1^{\circ}$ is smooth.
    Suppose that there exists some $b \in B$ such that $\deg(\omega_X|_{Y_b})=-1$.
    If we denote by $B' \subset B$ the connected component of $B$ containing $b$ and we denote by $Y' \to B'$ the restriction of $Y$ over $B'$, then it follows that all fibers (not just geometric) of $Y' \to B'$ are isomorphic to $\mathbb{P}^1$, the codimension of $Y'$ is $1$, and for all $b' \in B'$ we have $\cN_{Y'/X}|_{b'} \cong \cO_{\mathbb{P}^1_{b'}}(-1)$.
    Using the smoothness of $M_1^{\circ}$, it follows that around a neighborhood of $B' \subset M_1$ we have that $M_1^{\circ}$ coincides with the Fujiki-Nakano contraction of $Y' \to B'$ as in \Cref{example: fujiki nakano surgery}.
    Hence, around a neighborhood of $B'$, the scheme $X$ is a blowup of $M_1^{\circ}$ at $B' \subset M_1^{\circ}$.
    If we denote by $\overline{X}$ the smooth blow down of $Y' \to B'$, then we have that $M_1^{\circ}$ is obtained as the corresponding moduli of sheaves on $\overline{X}$ where we consider the $\mathbb{P}^1$-fibration $(Y \setminus Y') \to (B \setminus B')$ with $(Y \setminus Y') \subset \overline{X}$.
    % https://q.uiver.app/#q=WzAsOSxbMiwwLCJYIl0sWzIsMSwiXFxvdmVybGluZSBYIl0sWzIsMiwiTV8xIl0sWzMsMSwiTV8xXlxcY2lyYyJdLFsxLDAsIlkiXSxbMSwyLCJCIl0sWzAsMSwiWVxcc2V0bWludXMgWSciXSxbMCwyLCJCXFxzZXRtaW51cyBCJyJdLFsxLDFdLFswLDFdLFsxLDJdLFszLDJdLFs0LDAsIiIsMix7InN0eWxlIjp7InRhaWwiOnsibmFtZSI6Imhvb2siLCJzaWRlIjoidG9wIn19fV0sWzUsMiwiIiwyLHsic3R5bGUiOnsidGFpbCI6eyJuYW1lIjoiaG9vayIsInNpZGUiOiJ0b3AifX19XSxbNyw1LCIiLDIseyJzdHlsZSI6eyJ0YWlsIjp7Im5hbWUiOiJob29rIiwic2lkZSI6InRvcCJ9fX1dLFs2LDddLFs0LDVdLFs2LDgsIiIsMSx7InN0eWxlIjp7InRhaWwiOnsibmFtZSI6Imhvb2siLCJzaWRlIjoidG9wIn0sImhlYWQiOnsibmFtZSI6Im5vbmUifX19XSxbOCwxXV0=
    \[\begin{tikzcd}[ampersand replacement=\&]
    	\& Y \& X \\
    	{Y\setminus Y'} \& {} \& {\overline X} \& {M_1^\circ} \\
    	{B\setminus B'} \& B \& {M_1}
    	\arrow[hook, from=1-2, to=1-3]
    	\arrow[from=1-2, to=3-2]
    	\arrow[from=1-3, to=2-3]
    	\arrow[hook, no head, from=2-1, to=2-2]
    	\arrow[from=2-1, to=3-1]
    	\arrow[from=2-2, to=2-3]
    	\arrow[from=2-3, to=3-3]
    	\arrow[from=2-4, to=3-3]
    	\arrow[hook, from=3-1, to=3-2]
    	\arrow[hook, from=3-2, to=3-3]
    \end{tikzcd}\]
    Since we have a fully faithful admissible embedding $D^b(\overline{X}) \hookrightarrow D^b(X)$ by Orlov's description of the derived category of a blowup \cite{orlov_monoidal}, in order to verify the DK-hypothesis we may iteratively replace $X$ with smooth blow-downs $\overline{X}$ to reduce to the case when $\deg(\omega_X|_{Y_b}) \neq -1$ for all $b \in B$.
\end{example}

We conclude the verification of the DK-hypothesis in this setting by dealing with the case when $\deg(\omega_X|_{Y_b}) =0$ for all $b \in B$. This corresponds to a flop of the family of rational curves $Y \to B$. We shall make use of the following result from \cite{bridgeland_maciocia_elliptic}.
\begin{lemma} \label{lemma: intersection theorem lemma}
    Let $X$ and $X'$ be smooth proper connected algebraic spaces of dimension $n$ over an algebraically closed field, and let $H$ be an object in $D^b(X \times X')$. For any closed point $s \in X'$, let us denote by $H_s \in D^b(X)$ the corresponding fiber. Suppose that the following hold:
    \begin{enumerate}
        \item We have $\Hom(H_s, H_t[j]) =0$ for all closed points $s\neq t \in X'$ and $j \notin [1,n-1]$.
        \item The closed subset $Z \subset X' \times X'$ consisting of pairs $(s,t)$ such that $\RHom(H_s,H_t) \neq 0$ has dimension at most $n+1$.
    \end{enumerate}
    Then, for all closed points $s \neq t \in X'$ we have $\RHom(H_s, H_t) =0$.
\end{lemma}
\begin{proof}
    The beginning of the proof in \cite[Thm. 6.1]{bridgeland_maciocia_elliptic} applies to this context with only minor modifications: in Eq. (3) in \textit{loc.~cit.}, one needs to consider $\mathcal S^\vee \otimes \cO_{(y_1,y_2)}$ instead of $\left( \mathcal S \otimes \cO_{(y_1,y_2)} \right)^\vee$; and we note that $\Ext^n_X(H_{y_2},H_{y_1}) = 0$ for $y_1\neq y_2$ by our assumption and not by Serre duality.
\end{proof}

We are now ready to show the main result of this subsection, for which we adapt the techniques in \cite{bridgeland_flops}.
\begin{thm} \label{thm: derived categories of surgeries}
    Let $X$ be a smooth, projective, geometrically connected scheme over a field $k$.
    Let $Y \hookrightarrow X$ be a smooth $\mathbb{P}^1$-fibration $Y \to B$ as in \cref{setup for derived categories}, and suppose that $B$ is connected and $M_1^{\circ}$ is smooth.
     Assume that the restriction $\omega_X|_{Y_b}$ of the canonical bundle $\omega_X$ satisfies $\deg(\omega_X|_{Y_b}) \leq 0$ for all $b \in B$.
    Then, the following hold:
    \begin{enumerate}
        \item The bounded derived category $D^b(M_1^{\circ})$ embeds as an admissible subcategory of $D^b(X)$.
        \item If in addition we have $\deg(\omega_X|_{Y_b})=0$ for all $b \in B$, then $D^b(M_1^{\circ})$ is equivalent to $D^b(X)$.
    \end{enumerate}
\end{thm}

\begin{proof}
    Since $B$ is connected,  $\deg(\omega_X|_{Y_b})$ is the same integer for all $b \in B$.
    The cases when this integer is negative are already dealt with by \Cref{example: DK-hypothesis for P^1 fibrations} and \Cref{example: DK-hypothesis divisorial contractions}.
    Therefore, we may assume without loss of generality that $\deg(\omega_X|_{Y_b}) =0$ for all $b \in B$.
    In this case, for all geometric points $b \to B$, the normal bundle $\cN_{Y_b/X}$ is isomorphic to either $\cO_{\mathbb{P}^1_b}^{\oplus \dim(X)-3} \oplus \cO_{\mathbb{P}^1_b}(-1)^{\oplus 2}$ or $\cO_{\mathbb{P}^1_b}^{\oplus \dim(X)-2} \oplus \cO_{\mathbb{P}^1_b}(-2)$.
    It follows that $Y \neq X$, that $M_1^{\circ}$ is birational to $X$, and that the fibers of $p_2$ in diagram \eqref{equation: surgery diagram for P^1 fibration} are of dimension at most one (by the description of the nontrivial fibers in \Cref{example: surgery of P ell fibrations}).
    
    Since $M_1^{\circ}$ is a fine moduli of sheaves, there is a universal $M_1^{\circ}$-flat coherent sheaf $F_{\univ}$ on $X \times M_1^{\circ}$. Let $\pi: X \times M_1^{\circ} \to M_1^{\circ}$ denote the second projection. Note that for all field extensions $K \supset k$ and all points $y \in M_1^{\circ}(K)$, the fiber $(F_{\univ})_y$ is of one of the following forms:
    \begin{enumerate}
        \item[(a)] $(F_{\univ})_y \cong \cI_x$ for a closed point $x \in X_K(K)$ in the complement of $Y_K \subset X_K$, or
        \item[(b)] $(F_{\univ})_y$ is isomorphic to the middle term $F$ of a nonsplit extension
        \[ 0 \to i_*(\cO_{\mathbb{P}^1_K}(-1)) \to F \to \cI_{Y_b} \to 0, \]
        where $b \in B(K)$ is such that $Y_b \cong \mathbb{P}^1_K$ and $i: Y_b \hookrightarrow X_K$ denotes the corresponding closed immersion.
    \end{enumerate}
    In particular, it follows by inspection that $\Hom((F_{\univ})_y, \cO_{X_K}) = K$ for all such $y$. 
    
    We claim that $\pi_*(\cH om(F_{\univ}, \cO_{X\times M_1^{\circ}}))$ is a line bundle $L$ on $M_1^{\circ}$, and that we have a natural identification $\Hom(F_{\univ}, \pi^*(Q)) = \Hom(L, Q)$ for all quasicoherent sheaves $Q$ on $M_1^{\circ}$.
    Indeed, this can be checked \'etale locally on $M_1^{\circ}$, and after replacing $M_1^{\circ}$ with an affine atlas we may apply \cite[Cor. 7.7.8]{egaiii_2} to conclude that there exists a coherent sheaf $N$ on $M_1^{\circ}$ such that for all quasicoherent sheaves $Q$ on $M_1^{\circ}$ there is a natural identification $\Hom(F_{\univ}, \pi^*(Q)) = \Hom(N, Q)$. The fact that $\Hom((F_{\univ})_y, \cO_{X_K}) = K$ for all $y \in M_1^{\circ}(K)$ and the reducedness of $M_1^{\circ}$ implies, similarly as in the proof of \Cref{lemma: global hom from fiberwise homs}, that $N$ is locally free of rank $1$, and that we must have $L:= N^{\vee} \cong \pi_*(\cH om(F_{\univ}, \cO_{X\times M_1^{\circ}}))$, as desired.

    The counit $\pi^*(L) = \pi^*\pi_*(\cH om(F_{\univ}, \cO_{X\times M_1^{\circ}})) \to \cH om(F_{\univ}, \cO_{X\times M_1^{\circ}})$ corresponds to a homomorphism $\varphi: F_{\univ} \to \pi^*(L)^{\vee}$ such that fiberwise we have $\varphi_y \neq 0$ for all $y \in M_1^{\circ}$. We denote by $H \in D^b( X \times M_1^{\circ})$ the cone of the morphism $\varphi$. By construction, for all field extensions $K \supset k$ and all $y \in M_1^{\circ}(K)$, the fiber $H_y$ is of one of the following forms:
    \begin{enumerate}
        \item[(a')] If $(F_{\univ})_y \cong \cI_x$ for some closed point $j: x \hookrightarrow X_K$ as in (a) above, then $H_y \cong j_*(\cO_x)$ is the structure sheaf of $x$.
        \item[(b')] If $(F_{\univ})_y$ is a nonsplit extension as in (b) above, determined by a morphism $\alpha: \cI_{Y_b} \to i_*(\cO_{\mathbb{P}^1_K}(-1))[1]$ in the derived category, then $H_y$ fits into a nonsplit exact triangle
        \[ i_*(\cO_{\mathbb{P}^1_K}(-1))[1] \to H_y \to i_*(\cO_{\mathbb{P}^1_K}) \]
        which is determined by the composition $i_*(\cO_{\mathbb{P}^1_K}) \xrightarrow{\beta} \cI_{Y_b}[1] \xrightarrow{\alpha[1]} i_*(\cO_{\mathbb{P}^1_K}(-1))[2]$, where $\beta$ fits into the exact triangle
        \[   \cO_{X_K} \to  i_*(\cO_{\mathbb{P}^1_K}) \xrightarrow{\beta} \cI_{Y_b}[1] \]
        obtained by rotating the exact triangle corresponding to the short exact sequence $0 \to \cI_{Y_b} \to \cO_{X_K} \to i_*(\cO_{\mathbb{P}^1_K}) \to 0$.  
    \end{enumerate}
    We note that, in case (b') above, precomposition with $\beta$ induces an isomorphism 
    \[\Hom(\cI_{Y_b}[1], i_*(\cO_{\mathbb{P}^1_K}(-1))[2]) \xrightarrow{\sim} \Hom(i_*(\cO_{\mathbb{P}^1_K}), i_*(\cO_{\mathbb{P}^1_K}(-1))[2]).\]
    Therefore, if two points $y,y' \in M_1^{\circ}(K)$ yield non-isomorphic extensions $(F_{\univ})_y \not\cong (F_{\univ})_{y'}$, then $H_y \not\cong H_{y'}$.

    We denote by $\Phi: D^b(M_1^{\circ}) \to D^b(X)$ the Fourier-Mukai functor whose kernel is $H \in D^b(X \times M_1^{\circ})$. To conclude the proof of the theorem, we shall prove the following:

    \paragraph{\textbf{Claim.}} $\Phi$ induces a fully faithful embedding $D^b(M_1^{\circ}) \to D^b(X)$.

    Let us first explain why the claim implies the theorem. Assume that the claim holds. By our assumption that $\deg(\omega_X|_{Y_b})=0$ for all $b \in B$, it follows that $\omega_X|_{Y_b}$ is trivial.
    From our description of the fibers $H_y$ above, it follows that $H_y \otimes \omega_X =H_y$ for all $y$. Applying \cite[Thm. 2.4]{bridgeland_king_reid} to the spanning class $\Omega$ of structure sheaves of closed points of $M_1^{\circ}$, it follows that $\Phi$ is an equivalence. So Part (2) of the theorem follows.

    We are left to show the claim that $\Phi$ is fully faithful.
    We recall from \cite[Thm. A]{hall2024generalizedbondalorlovfaithfulnesscriterion} (cf. \cite[Thm. 1.1]{bondalorlov1995sod} in the setting of smooth projective varieties) that the Fourier-Mukai functor $\Phi$ is a fully faithful embedding if and only if the following two conditions hold:
    \begin{itemize}
        \item $\RHom(\Phi(\cO_s),\Phi(\cO_t)) = 0$ for any pair of distinct closed points $s,t \in M^\circ_1$, and
        \item $\Hom(\Phi(\cO_t),\Phi(\cO_t)) = \kappa(t)$, and $\Hom(\Phi(\cO_t),\Phi(\cO_t)[j]) = 0$ for any closed point $t \in M^\circ_1$ and any $j$ outside of the interval $[0,\dim(M^\circ_1)]= [0,\dim(X)]$.
    \end{itemize}
    We note that it is enough to check these conditions when we base-change our whole setup to the algebraic closure of $k$, and hence we may assume without loss of generality that the ground field $k$ is algebraically closed.

    By inspecting the classification of fibers $H_y$ (see (a') and (b') above), the only nontrivial case is when $H_s$ and $H_t$ are both nonsplit extensions 
    \begin{equation} \label{equation: defining triangles}
        i_*(\cO_{\mathbb{P}^1_k}(-1))[1] \to H_s \to i_*(\cO_{\mathbb{P}^1_k}) \; \;  \text{and} \; \; i_*(\cO_{\mathbb{P}^1_k}(-1))[1] \to H_t \to i_*(\cO_{\mathbb{P}^1_k})
    \end{equation}
    of type (b') with the same embedding $i: Y_b\cong \mathbb{P}^1_k \hookrightarrow X_K$.
    Note that this implies that $p_2(s) = p_2(t)$ in diagram \eqref{equation: surgery diagram for P^1 fibration}, so in particular $H_s$ and $H_t$ are orthogonal whenever $p_2(s) \neq p_2(t)$.
    
    We begin by checking the required cohomology vanishing when $s \neq t$, which by our discussion above implies $H_s \not \cong H_{t}$. 
    Let $j \leq 0$, and consider a morphism $\psi: H_s \to H_t[j]$. 
    Since $j-1<0$, we have $\Hom(i_*(\cO_{\mathbb{P}^1_k}(-1)), i_*(\cO_{\mathbb{P}^1_k})[j-1])=0$, and hence the composition $i_*(\cO_{\mathbb{P}^1_k}(-1))[1] \to H_s \xrightarrow{\psi} H_t[j] \to i_*(\cO_{\mathbb{P}^1_k})[j]$ is trivial, where left-most and right-most maps are from the exact triangles in \eqref{equation: defining triangles} above. Therefore, the morphism $i_*(\cO_{\mathbb{P}^1_k}(-1))[1] \to H_s \xrightarrow{\psi} H_t[j]$ factors through a morphism $\psi_1: i_*(\cO_{\mathbb{P}^1_k}(-1))[1] \to i_*(\cO_{\mathbb{P}^1_k}(-1))[1+j]$.  The axiom TR3 of triangulated categories then implies that $\psi$ fits into a morphism between the triangles in \eqref{equation: defining triangles}:
   \[
    \begin{tikzcd}
      i_*(\cO_{\mathbb{P}^1_k}(-1))[1] \ar[d, "\psi_1"]
      \ar[r] & H_s \ar[d, "\psi"] \ar[r] & i_*(\cO_{\mathbb{P}^1_k}) \ar[d, "\psi_2"] \\ i_*(\cO_{\mathbb{P}^1_k}(-1))[1+j] 
      \ar[r] & H_t[j]  \ar[r] & i_*(\cO_{\mathbb{P}^1_k})[j].
    \end{tikzcd}
    \]
    Furthermore, a direct computation (using the non-positivity of the normal bundle of $i: Y_b \hookrightarrow X$ for $j=0$) shows that $\Hom(i_*(\cO_{\mathbb{P}^1_k}), i_*(\cO_{\mathbb{P}^1_k}(-1))[1+j]) =0$, and so it follows that $\psi$ is uniquely determined by $\psi_1$ and $\psi_2$. If $j<0$, we note that $\psi_1$ and $\psi_2$ are forced to be $0$, and so we conclude that $\Hom(H_s, H_t[j]) =0$. On the other hand, if $j=0$, then $\psi_1$ and $\psi_2$ are scalar multiples of the identity, and hence $\Hom(H_s, H_t)$ is at most two-dimensional over $k$. Moreover, since the exact triangles are non-split, we cannot have that one of $\psi_1,\psi_2$ is $0$ while the other one is nonzero, and hence we conclude that $\Hom(H_s, H_t)$ is at most one-dimensional. Note that we have not used $s \neq t$ so far. 

    If $s \neq t$, then $H_s \not\cong H_t$ and hence we cannot have both $\psi_1$ and $\psi_2$ being nonzero. Therefore $\Hom(H_s,H_t)=0$ when $s \neq t$. Under our assumption that $\deg(i^*(\omega_X))=0$, we have $i^*(\omega_X) \cong \cO_{\mathbb{P}^1_k}$, and hence by Serre duality we have that
    \begin{gather*}
        \Hom(H_s, H_t[h]) = \Hom(H_t, H_s \otimes \omega_X[\dim(X)-h])^{\vee}= \Hom(H_t, H_s[\dim(X)-h])^{\vee} =0
    \end{gather*}
    if $h>\dim(X)$, or if $s \neq t$ and $h \geq \dim(X)$. We have in particular shown that $\Hom(H_s, H_t[j]) =0$ if $s\neq t$ and $j \notin [1,\dim(X)-1]$. Now using the orthogonality of $H_s$ and $H_t$ whenever $p_2(s) \neq p_2(t)$ jointly with the fact that $p_2$ has fibers of dimension at most one, we may apply \Cref{lemma: intersection theorem lemma} to conclude that $\RHom(H_s, H_t) =0$ whenever $s \neq t$.
    
    We are left to deal with the case when $s=t$. We have already shown that $\Hom(H_s, H_s[j]) =0$ if $j \notin [0,\dim(X)]$ above. Also, we have shown that $\Hom(H_s, H_s)$ is at most one-dimensional, and so we have $\Hom(H_s, H_s) =k$, as desired. \qedhere
\end{proof}

\appendix

\section{Properties of determinants}
\label{app: properties of determinants}

We recall the main properties of determinants that we use freely throughout the article. Let \(X\) be a smooth projective variety over \(k\) and \(S\) be a \(k\)-scheme.
We set \(X_S := X \times S\).
Let \(F\) be an \(S\)-flat quasicoherent sheaf on \(X_S\) locally of finite presentation.
It follows that, locally on \(S\), the sheaf \(F\) admits a finite resolution by locally free sheaves \cite[pg. 36,37]{huybrechts.lehn}.
By \cite[Theorem 2]{knudsenmumford1976}, the determinant \(\det(F)\) is defined for such sheaves.
The determinant \(\det (F)\) is a line bundle, its formation commutes with pullbacks, and it is multiplicative in short exact sequences.

\begin{lemma}[{\cite[Defn. 2.7, Lem. 2.8]{rho-sheaves-paper}}]
\label{lemma: det(rank 1 tf) = double dual}
    Let \(F\) be an $S$-flat sheaf on $X_S$ as above. If for all $s \in S$, the fiber \(F_s\) is torsion-free of rank one, then \(F^{\vee\vee} \cong \det(F)\). \qed
\end{lemma}

\begin{lemma}
\label{lemma: det is an effective Cartier divisor}
    Let \(S\) be the spectrum of a discrete valuation ring \(R\) or of a field \(K\).
    Let \(i: Y \hookrightarrow X_S\) be the inclusion of a closed subscheme with $Y \neq X_S$, and let \(G\) be a coherent sheaf on \(Y\).
    Then there is an $S$-flat effective Cartier divisor \(D\) on \(X_S\) such that \(\det(i_*(G))\cong \cO(D)\).
    More precisely, if \(\codim(Y) = 1\) and \(Y\) is integral, then the divisor is \(D = \rk(G) \cdot Y\),
    and if \(\codim(Y) \geq 2\), then \(\det(i_*(G))\cong \cO_{X_S}\).
\end{lemma}

\begin{proof}
    When \(S = \Spec (K)\), the claim follows from the coefficient formula in \cite[p. 48 and p. 49, Thm. 3 (vi)]{knudsenmumford1976}.
    When \(S = \Spec(R)\), denote by \(\eta\) the generic point of \(S\).
    Then by the above, \(\det(i_*(G))_\eta \cong \det(i_*(G)_{\eta})\) is an effective Cartier divisor $D$ on $X_{\eta}$. Consider the $S$-flat effective Cartier divisor $\overline{D} \subset X_S$. By construction, $\cO(\overline{D})$ and $\det(i_*(G))$ agree when restricted to the generic fiber $X_{\eta}$. This means that we have $\det(i_*(G)) = \cO(\overline{D} + n \cdot X_s)$ for some integer $n$. Since the line bundle $\cO(X_s)$ is trivial, we get $\cO(\overline{D}) \cong \det(i_*(G))$.
\end{proof}

\footnotesize{\bibliography{moduli-sheaves.bib}}
\bibliographystyle{alpha}
\end{document}